\documentclass[1p]{elsarticle}
\usepackage{hyperref,bm,amsmath,amssymb,algorithm,amsfonts,xcolor}
\usepackage[]{algpseudocode}
\newdefinition{rmk}{Remark}

\newcommand\norm[1]{\left\lVert#1\right\rVert}
\newcommand\abs[1]{\left\lvert#1\right\rvert}
\usepackage{subcaption}
\usepackage{color}
\usepackage{tcolorbox}

\usepackage{lipsum}
\makeatletter
\def\ps@pprintTitle{%
	\let\@oddhead\@empty
	\let\@evenhead\@empty
	\def\@oddfoot{ }%
	\let\@evenfoot\@oddfoot}
\makeatother

\begin{document}

\begin{frontmatter}

\title{A localized reduced basis approach for unfitted domain methods on parameterized geometries\tnoteref{mytitlenote}}

\tnotetext[mytitlenote]{This is the peer-reviewed manuscript accepted in Computer Methods in Applied Mechanics and Engineering \href{https://doi.org/10.1016/j.cma.2023.115997}{https://doi.org/10.1016/j.cma.2023.115997}.}
\author[mymainaddress]{Margarita Chasapi\corref{mycorrespondingauthor}}
\cortext[mycorrespondingauthor]{Corresponding author}
\ead{margarita.chasapi@epfl.ch}

\author[mymainaddress]{Pablo Antolin}
\author[mymainaddress,mysecondaryaddress]{Annalisa Buffa}

\address[mymainaddress]{Institute of Mathematics, \'Ecole Polytechnique F\'ed\'erale de Lausanne, Lausanne, Switzerland }
\address[mysecondaryaddress]{Instituto di Matematica Applicata e Tecnologie Informatiche 'E. Magenes' (CNR), Pavia, Italy}

\begin{abstract}
This work introduces a reduced order modeling (ROM) framework for the solution of parameterized second-order linear elliptic partial differential equations formulated on unfitted geometries. The goal is to construct efficient projection-based ROMs, which rely on techniques such as the reduced basis method and discrete empirical interpolation. The presence of geometrical parameters in unfitted domain discretizations entails challenges for the application of standard ROMs. Therefore, in this work we propose a methodology based on i) extension of snapshots on the background mesh and ii) localization strategies to decrease the number of reduced basis functions. The method we obtain is computationally efficient and accurate, while it is agnostic with respect to the underlying discretization choice. We test the applicability of the proposed framework with numerical experiments on two model problems, namely the Poisson and linear elasticity problems. In particular, we study several benchmarks formulated on two-dimensional, trimmed domains discretized with splines and we observe a significant reduction of the online computational cost compared to standard ROMs for the same level of accuracy. Moreover, we show the applicability of our methodology to a three-dimensional geometry of a linear elastic problem. 
\end{abstract}


\end{frontmatter}

\section{Introduction}

In recent years, unfitted domain methods have attracted a lot of attention. The main idea behind these methods is embedding a geometric representation into a simple background domain. A wide class of immersed and embedded methods fall within this category, where the geometry is decoupled from the discretization of the solution. This is the case, for example, in the immersed boundary method \cite{Peskin2002}, the immersed interface \cite{Li2006}, the fictitious domain and finite cell method \cite{parvizian2007finite}, CutFEM \cite{burman2015cutfem}, the Shifted Boundary Method \cite{Main2018,Main2018b}, and others. Some of the challenges involved in immersed methods have been the focus of research activities in the past, such as numerical integration and imposition of boundary conditions. We further refer the reader to the review in \cite{Mittal2005}.

Moreover, great efforts have been devoted in the research area of design-oriented simulation. Isogeometric Analysis (IGA) was introduced in \cite{Hughes2005} and has been successfully applied in several fields of computational science and engineering. The idea behind IGA is to adopt the same representation employed in Computer Aided Design (CAD), such as B-splines and their rational variants (NURBS), for the approximation of the solution in finite element analysis. This paradigm provides a unified framework from design to analysis that is capable of simplifying the clean-up and meshing of geometric models. A review of isogeometric methods can be found in \cite{Hughes2005,Cottrell2007}. One of the main challenges in IGA is dealing with boundary representations (B-rep). Commercial CAD software are currently using B-reps for solid modeling, that is, the solid is only modeled by its boundary. Since the volumetric description is missing, such boundary representations are not analysis-suitable \cite{Akhras2016,Klinkel2020,Chasapi2021}. 
In the context of IGA, immersed methods have become particularly popular since they circumvent the need to construct volumetric parameterizations by simply embedding B-reps into a background domain \cite{Hollig2001,Schillinger2012,Elfeverson2018,Rank2012,Messmer2022}. Furthermore, immersed isogeometric methods bear connections to Boolean operations and trimming. In CAD, spline parameterizations are commonly trimmed in order to represent complex geometric shapes. The result of trimming operations are unfitted meshes, since the parameterization is defined in the original background domain. As trimming is the prevailing technology to represent complex shapes in CAD, its treatment is crucial to achieve a unified design-through-analysis framework and tackle geometries that are relevant from an industrial viewpoint. The reader is further referred to \cite{Marussig2018} for a detailed review and current challenges on trimming. Moreover, we refer to previous works addressing the challenges posed by trimming in the analysis, such as numerical integration \cite{Nagy2015a,Kudela2015,Antolin2019,Divi2020,Antolin2022,Antolin2022b}, conditioning \cite{Marussig2018b} and stabilization techniques \cite{Elfverson2019,Buffa2020}.

One further aspect to be considered is that in many cases partial differential equations (PDEs) need to be solved multiple times for several parametric configurations. This is the case in a real-time and many-query context arising, for example, in design optimization, uncertainty quantification, control, and other applications. Efficient reduced order modeling techniques are essential to establish a suitable offline/online procedure that achieves a computational speedup. To this end, projection-based reduced order models (ROMs) have been successfully employed in a wide range of domains for the solution of parameterized PDEs. 

In the past years, there have been several successful applications of IGA in the context of reduced order modeling techniques \cite{Garotta2020}. The fields of application span from fluid dynamics \cite{Manzoni2015,Salmoiraghi2016}, to parabolic problems \cite{Zhu2017} and cardiac electrophysiology \cite{Fresca2017}. Moreover, ROMs were constructed using IGA on complex, multipatch geometries and isotopological meshes in \cite{Maquart2020}. The combination of reduced basis methods (RB) and IGA has been particularly motivated by their combined advantages to solve PDEs on parameterized geometries \cite{Devaud2017}. The interested reader is referred to \cite{Hesthaven2016,QMN_RBspringer} for a thorough discussion on reduced basis methods as well as to our previous work \cite{Chasapi2022} in the context of isogeometric ROMs and domain decomposition. We also refer to \cite{Gabriel2022} in the context of IGA and tensor train compression for parameterized geometries. On the other hand, the combination of unfitted domain methods with ROMs and in particular IGA is still in its infancy. Model  reduction within an embedded framework was first investigated in \cite{Nouy2011} for uncertain parameterized geometries, where a fictitious domain method was combined with Proper Generalized Decomposition (PGD) and also in \cite{Balajewicz2014} for interfaces evolving in time, where a low-rank approximation  was formulated for snapshot compression. Recently, ROMs were combined with CutFEM \cite{Karatzas2020,Karatzas2021} and the Shifted Boundary Method \cite{Karatzas2019,Karatzas2020b,Zeng2022} for parameterized geometries. These works address aspects related to embedded methods, such as definition of solutions on a common mesh through suitable extension and transportation of snapshots on the background domain. Their main advantage is that they avoid remeshing and overcome the need for reference domain formulations typically used for ROMs on parameterized geometries, where the snapshots are mapped to a reference domain and the transformation depends highly on the given problem at hand. Nevertheless, these works do not resolve the nonaffine dependence of differential operators and transport maps on the geometrical parameters, i.e., do not make use of hyper-reduction techniques, which is essential for an efficient offline/online decomposition. In the context of unfitted finite elements, hyper-reduction was combined with ROMs for PDE-constrained optimization in \cite{Karatzas2022}. However, the combination of reduced basis and isogeometric methods formulated on unfitted geometries has not been thoroughly investigated in previous studies.

In this work, we propose a full reduction framework for nonaffine problems on parameterized unfitted geometries that relies on hyper-reduction techniques to achieve an efficient solution on the fly. The proposed framework is agnostic with respect to the chosen discretization and can be applied in combination with finite element as well as isogeometric methods formulated on unfitted domains.

Our contribution falls within the context of projection-based ROMs, while we employ the Proper Orthogonal Decomposition (POD) to construct reduced basis spaces. We remark that the RB method is based on the assumption of affine dependence of the operators on the parameters. Since we are interested in parameterized geometries, this assumption is not always fulfilled. To recover the affine dependence on the geometrical parameters we rely on the empirical interpolation method (EIM) \cite{Barrault2004} and in particular its discrete variant (DEIM) for vectors and matrices \cite{Chaturantabut2010,Negri2015}. We recall that the solution of PDEs on parameterized embedded domains necessitates proper extension of solutions to a common mesh. In this respect, our approach is inspired by previous studies on snapshot extension techniques \cite{Karatzas2020,Karatzas2019}. The extended solutions for varying geometrical parameters may exhibit discontinuities for different values of the parameters. In fact, it is inefficient to construct a reduced basis approximation with a single, linear subspace since a very large set of global reduced basis functions is required to achieve a sufficient accuracy. It should be noted that even for a large dimension of the basis the approximation properties may be poor and characterized by oscillatory behavior as discussed in \cite{Benner2021} and references therein.

On the other hand, strategies based on the idea of local reduced bases have been proposed in the past to circumvent these shortcomings. Local ROMs based on clustering of solution snapshots were first introduced in \cite{Amsallem2012} and further proposed, for example, in the context of discrete empirical interpolation \cite{Peherstorfer2014}, cardiac electrophysiology \cite{Pagani2018} and bifurcation problems \cite{Hess2019}. In this work, we propose a parameter-based partitioning of snapshots. This approach bears connections to $hp$-reduced basis methods introduced in \cite{Eftang2010a,Eftang2011} for elliptic and parabolic PDEs and later in \cite{Eftang2012} for empirical interpolation. We also refer to \cite{Haasdonk2011, Maday2013} for adaptive local reduced bases. Finally, we illustrate the methodology with numerical experiments employing spline discretizations on trimmed geometries. The features of the proposed ROM framework are summarized as follows:
\begin{itemize}
    \item It is agnostic to the underlying discretization and cutting operations performed on parameterized unfitted domains.
    \item It enables an efficient offline/online decomposition for nonaffine problems with geometrical parameters based on  hyper-reduction. The latter is applied to the algebraic structure of the differential operators and combined with interpolation to ensure a fast and less intrusive treatment of nonaffine dependencies. 
    \item It allows to construct efficient ROMs based on a localization strategy. The online cost is low, as the dimension of the local problems is small and one can easily switch between subspaces in the online phase.
\end{itemize}

We structure this contribution as follows: Section 2 provides a brief review of the main concepts related to unfitted domain discretizations formulated on parameterized geometries. Section 3 presents the generic framework of parameterized linear elliptic PDEs considered throughout this work. In Section 4 we discuss the reduced basis method for PDEs on unfitted geometries, in particular the snapshots extension, reduced space construction and the discrete empirical interpolation method. In Section 5 we present the localization strategy to construct efficient ROMs. Section 6 provides several numerical experiments in order to assess both the computational efficiency and accuracy of the proposed methodology. Finally, the main conclusions are summarized in Section 7.

\section{Unfitted domain discretization of geometrically parameterized
problems}\label{sec:unfitted}

\begin{figure}[h!]
\centering
\includegraphics[width=0.8\textwidth]{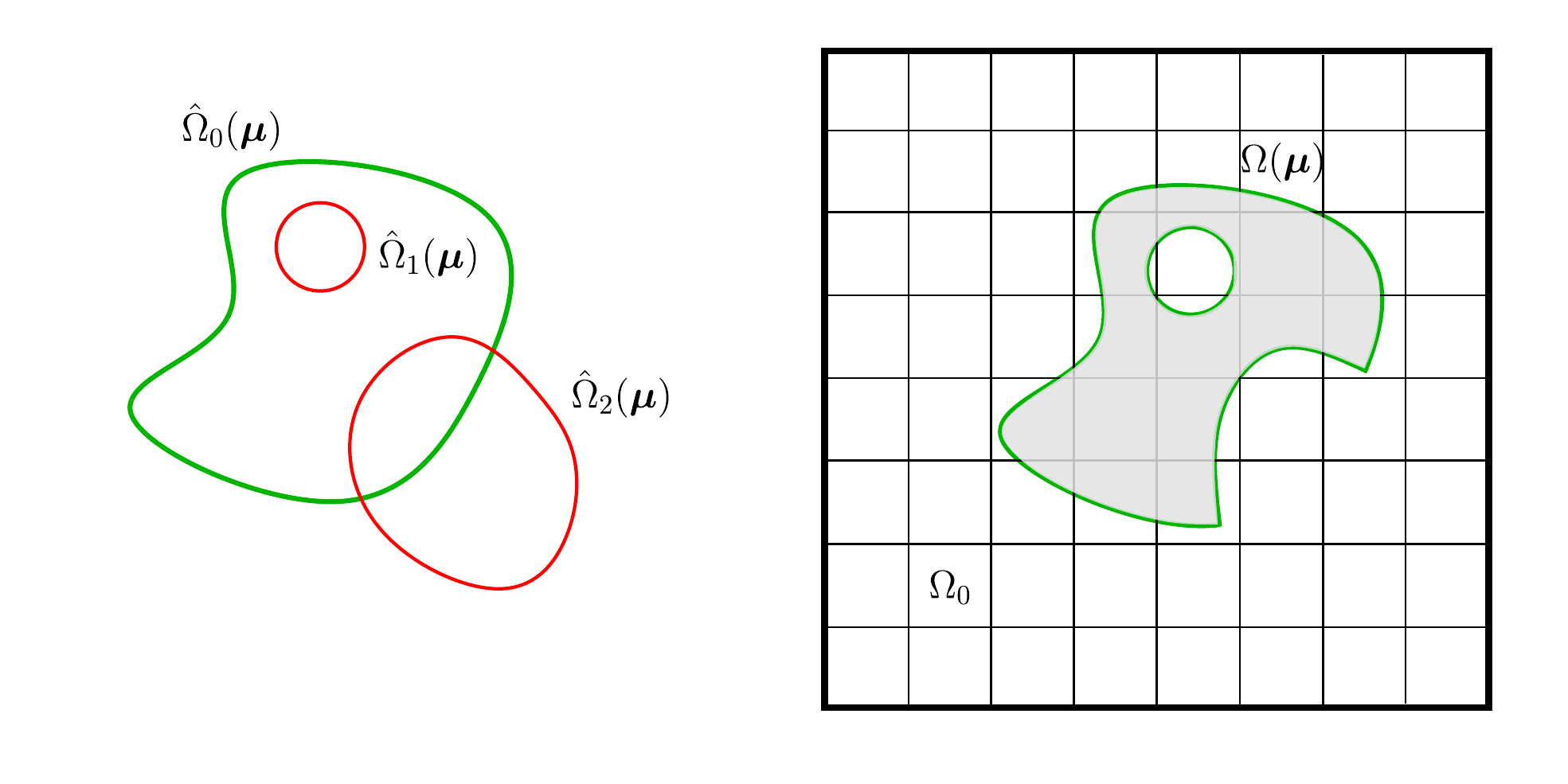}
\caption{Geometrical setting: the geometrically parameterized domain
$\Omega(\bm{\mu})$, built through subtraction operations, is embedded in
the background domain $\Omega_0$.}
\label{fig:geometry}
\end{figure}

In this section we provide a brief overview of some basic concepts related to
unfitted boundary methods, which will constitute the basis of the methods to be
developed along the manuscript.
Let us denote as $\Omega(\bm{\mu})\subset\mathbb{R}^d$ the domain in which we want to solve our PDE problem at hand, where $d$ is the spatial dimension, either 2D or 3D.
This physical domain is described by a set of geometrical parameters
$\bm{\mu}\in\mathcal{P}\subset\mathbb{R}^M$, where $\mathcal{P}$ is the space
of parameters and $M>0$ is the number of parameters.
In this work we assume that $\Omega(\bm{\mu})$ is built by
cutting a master domain $\hat\Omega_0(\bm{\mu})$ with a series of $K>0$ domains
$\left\lbrace\hat\Omega_i(\bm{\mu})\right\rbrace^K_{i=1}$, all of them potentially
dependent on the geometrical parameters $\bm{\mu}$, as (see
Figure \ref{fig:geometry})
\begin{align}\label{eq:cut}
\Omega(\bm{\mu}) = \hat\Omega_0(\bm{\mu}) \setminus
\bigcup^K_{i=1}\overline{\hat\Omega_i(\bm{\mu})}\,.
\end{align}

In order to deal with such generic domain parameterizations in a way that is compatible with reduced order modeling techniques we rely upon unfitted domain methods.
Thus, we assume that $\Omega(\bm{\mu})$ is immersed in a larger \emph{background domain} $\Omega_0$, as
\begin{align}\label{eq:boundingbox}
\Omega(\bm{\mu})\subset\Omega_0\subset\mathbb{R}^d\,\
\forall\bm{\mu}\in\mathcal{P}\,,
\end{align}
that is independent of the geometrical parameters.
$\Omega_0$ is the base upon which we create a discrete functional space
$V_{h,0}$, that is independent of $\bm\mu$, defined in a general way as
\begin{align}
V_{h,0} = \text{span}\left\lbrace \mathcal{B}_i,\
i=1,\dots,\mathcal{N}_{h,0} \right\rbrace,
\end{align}
where $\left\lbrace\mathcal{B}_i\right\rbrace_{i=1}^{\mathcal{N}_{h,0}}$ is a basis of the space,
being $\mathcal{N}_{h,0}=\text{dim}\left(V_{h_0}\right)$. Then, for a given set of geometrical parameters $\bm{\mu}$, we want to solve our problem of interest on the domain $\Omega(\bm{\mu})$, and for that purpose we introduce a smaller space $V_h(\bm{\mu})$, defined as
\begin{align}\label{eq:space}
V_h(\bm{\mu}) = \text{span}\left\lbrace \mathcal{B}_i \in V_{h,0}
\ \vert\ \text{supp}(\mathcal{B}_i) \cap \Omega(\bm{\mu}) \neq\emptyset \right\rbrace,
\end{align}
being $\mathcal{N}_h(\bm{\mu}) = \text{dim}(V_h(\bm{\mu}))$. Clearly it holds $\mathcal{N}_h(\bm{\mu})\leq\mathcal{N}_{h,0}$.
Using such a space for discretizing a differential problem on $\Omega(\bm{\mu})$, the domain partition (i.e. meshing) is decoupled from the solution discretization, which renders this family of methods flexible alternatives to traditional boundary fitted techniques.
Furthermore, it is worth highlighting that the basis functions are defined on
$\Omega_0$, and so, their definition remains unchanged, not depending on
the geometrical parameters. In what follows we will rely upon our previous works for the computation of integrals arising from unfitted domain discretizations \cite{Antolin2019,Wei2021}.

From the definition~\eqref{eq:space}, it is clear that a set of basis functions may be inactive (those whose support does not intersect $\Omega(\bm{\mu})$).
They do not contribute to the solution discretization and therefore they will be just simply not considered.
In addition, that set of functions depends on the particular choice of $\bm\mu$, and may change from problem to problem.
We will further elaborate on this aspect and its implications in the context of efficient reduced order modeling in Section \ref{sec:ROM}.

\begin{rmk}\label{remark1}
The numerical experiments discussed in Section \ref{sec:examples} are based on the use of spline
discretizations (i.e., Isogeometric Analysis~\cite{Hughes2005}).
Therefore, in those examples we assume $\Omega_0$ to be a Cartesian product
domain, as, for instance, the bounding box of all the possible domains
$\Omega(\bm{\mu})$.
We define a Cartesian mesh in such domain (see Figure \ref{fig:geometry}) upon
which the discretization spline space is built.
However, the framework presented above and used in the following sections is completely agnostic with respect
to the discretization chosen for $V_{h,0}$, and it can be applied in
combination with other unfitted domain techniques as, e.g., the Finite Cell
Method~\cite{parvizian2007finite},
CutFEM~\cite{burman2015cutfem}, or the aggregated unfitted domain 
method~\cite{BADIA2018533}.

In addition, due to the fact that in Section \ref{sec:examples} we limit our discussion to the
Poisson and linear elasticity problems, there we assume that $V_{h,0}\subset
H^1(\Omega_0)$. However, the method introduced in this work applies to other
problems and space choices.
\end{rmk}
\begin{rmk}\label{remark2}
As it will be illustrated in the example of Section \ref{sec:ring}, the geometry
$\Omega(\bm{\mu})$ can be further transformed using a mapping
$\bm{F}:\Omega_0\times\mathcal{P}\to\mathbb{R}^s$, with $s\geq d$, in a similar way as
proposed in \cite{Antolin2019}.
\end{rmk}

\section{Parameterized model problem}\label{sec:models}
In the following, we introduce a generic linear elliptic PDE that will serve as model problem for our exposition. We consider the following equation:
\begin{equation}\label{eq1}
\mathcal{L}u = {f} \quad  \text{in} \ \Omega({\bm{\mu}}),
\end{equation}
equipped with proper boundary conditions on the boundary $\partial{\Omega}({\bm{\mu}})$. We suppose that we have the well-posed discrete weak formulation of the parameterized problem in Equation \eqref{eq1} as: find $u_h \in V_h(\bm{\mu})$ such that 
\begin{equation}\label{eq2}
a(u_h,v_h;\bm{\mu}) = f(v_h;\bm{\mu}), \qquad \forall v_h \in V_h(\bm{\mu}),
\end{equation}
where $a(\cdot,\cdot;\bm{\mu}): V_h(\bm{\mu}) \times V_h(\bm{\mu}) \to \mathbb{R}$ is a bilinear, continuous, and coercive form and $f(\cdot;\bm{\mu}): V_h(\bm{\mu}) \to \mathbb{R}$ is a linear and continuous functional associated to a parameterized PDE for every $\bm{\mu}\in \mathcal{P}$. The space $V_h(\bm{\mu})$ is a discrete subspace whose choice depends in general on the boundary conditions. In the case where Dirichlet boundary conditions are imposed on an unfitted part of the boundary, a suitable stabilization technique must be adopted, as for example the one introduced in \cite{Buffa2020,Wei2021}. Note that in the numerical experiments discussed in Section \ref{sec:examples} we will consider homogeneous Dirichlet and Neumann boundary conditions for ease of exposition, while we will impose Dirichlet boundary conditions on the part of the boundary that  coincides with the boundary of the background domain $\partial{\Omega} ({\bm{\mu}})\cap \partial{\Omega_0}$ and not on the unfitted part of $\partial{\Omega}({\bm{\mu}})$. From the algebraic viewpoint, the discrete approximation leads to the following parameterized linear system of dimension $\mathcal{N}_h(\bm{\mu})=\text{dim}(V_h(\bm{\mu}))$
\begin{equation}\label{eq3}
{\bf{A}}(\bm{\mu}){\bf{u}}_h(\bm{\mu}) = {\bf{f}}(\bm{\mu}),
\end{equation}
where ${\bf{A}}$ $\in \mathbb{R}^{\mathcal{N}_h(\bm{\mu})\times \mathcal{N}_h(\bm{\mu})}$ is the stiffness matrix corresponding to the differential operator, ${\bf{f}} \in \mathbb{R}^{\mathcal{N}_h(\bm{\mu})}$ is the vector representing the source term and ${\bf{u}}_h(\bm{\mu})\in \mathbb{R}^{\mathcal{N}_h(\bm{\mu})}$ is the solution vector. In Section \ref{sec:examples} we will consider two model problems to validate the methods, namely the Poisson and linear elasticity problems. From now on, we refer to the problem \eqref{eq2} as high-fidelity or full order model (FOM). We are interested in solving Equation \eqref{eq3} for different values of the parameter vector $\bm{\mu}$ (order of at least hundreds) and analyze different geometrical representations. Driven by this, in the next sections we will turn to reduced order models as a means of tackling parameterized problems in an efficient manner.

\section{Reduced basis method for PDEs on parameterized unfitted geometries}\label{sec:ROM}
In the context of parameterized PDEs, the main idea behind projection-based reduced order models (ROMs) is to approximate the solution of FOMs based on a linear combination of \emph{global} reduced basis functions. The latter can be obtained from selected solutions of the FOM, which are referred to as \emph{snapshots}. In the following we will discuss the key features related to unfitted domain discretizations and briefly review some basic concepts in order to obtain an effective model order reduction.

\subsection{Snapshots extension}\label{sec:extension}
\begin{figure}[!htb]
	\centering
	\includegraphics[width=0.7\textwidth]{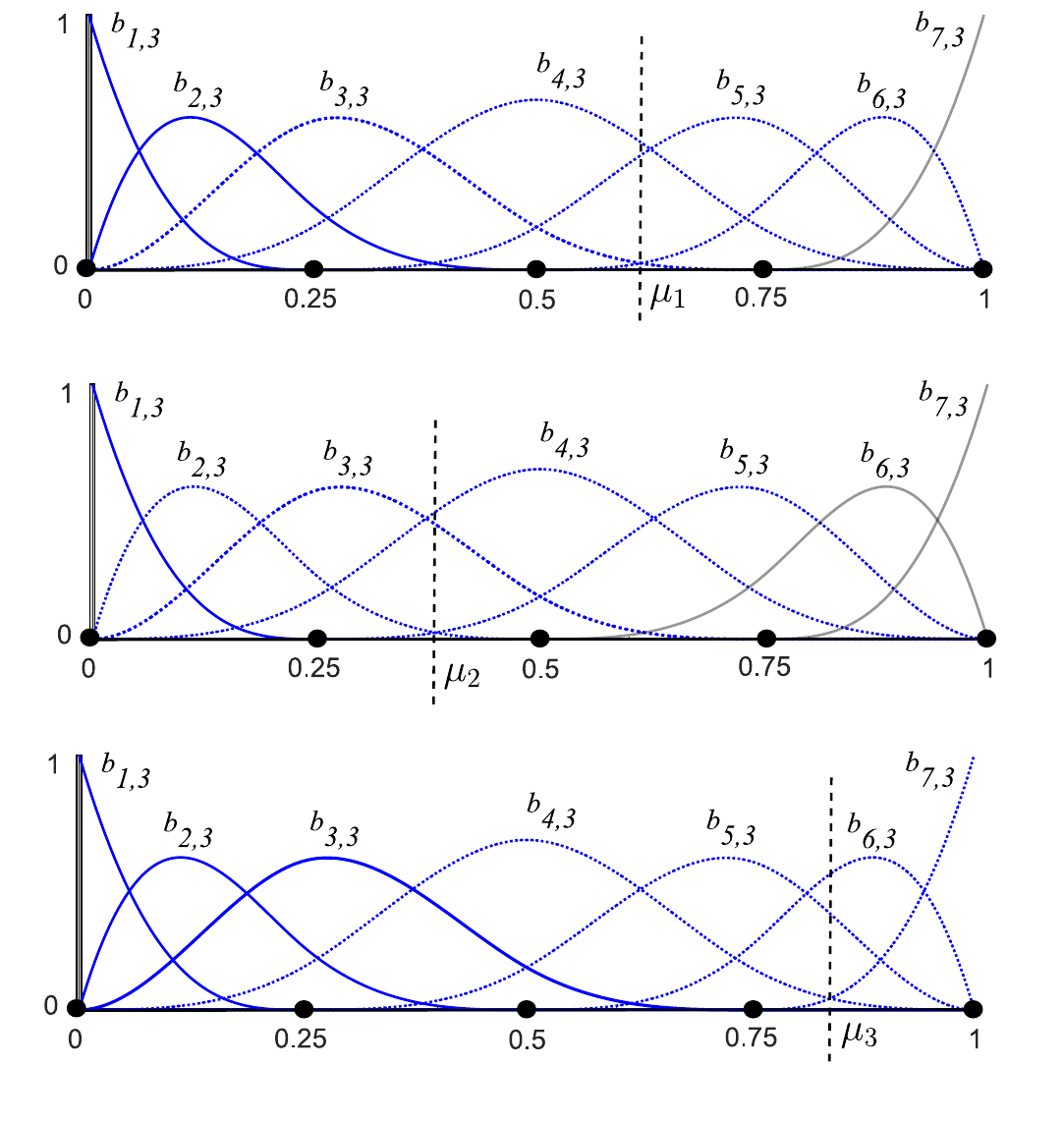} \\
	\caption{Univariate B-spline basis for different parameters $\mu_1,\mu_2,\mu_3$ defining the trimming location. The functions depicted in blue are fully active, in dotted blue are trimmed active and in grey are inactive.}\label{fig:Bsplines}
\end{figure}

The solution of problems on parameterized unfitted domains might vary highly over the parameter space $\mathcal{P}$. Let us consider a spline discretization as an illustrative example. In fact, as the active domain $\Omega(\bm{\mu})$ depends on $\bm{\mu}$, the support of B-spline basis is also $\bm{\mu}$-dependent: the basis functions that are active or inactive may change for different values of the parameters. An example of a trimmed univariate B-spline basis is illustrated in Figure \ref{fig:Bsplines}. Let us now consider a geometrical parameter $\mu$ affecting the location of the elements cut away from the Cartesian mesh upon which the spline discretization is built. The trimmed basis comprises basis functions that are cut and depicted in blue dotted lines in Figure \ref{fig:Bsplines}. The active basis functions with full support inside the domain $(0, \ \mu)$ are depicted in blue color, while the inactive functions, outside of $(0, \ \mu)$ are shown in grey. For different values of $\mu$, different basis functions are fully or partially active.

Thus, when solving problem \eqref{eq3} we seek a solution $u_h(\bm{\mu})$ in a spline space $V_h(\bm{\mu})$ \eqref{eq:space} whose set of active basis functions depends on the parameters $\bm{\mu}$ and the same holds for its dimension $\mathcal{N}_h(\bm{\mu})$. Indeed, depending on the parameters $\bm{\mu}$, the set of active basis functions may change for different snapshots. Practically, this implies that the snapshot solution vectors ${\bf{u}}_h(\bm{\mu})$ obtained from \eqref{eq3} might have different length. This fact will hinder the formation of snapshots matrices for constructing a reduced basis by techniques such as the Proper Orthogonal Decomposition (POD). Therefore, a suitable extension of the solution vectors ${\bf{u}}_h(\bm{\mu})$ has to be performed in order to render all vectors of the same length. Since the background domain 
$\Omega_0$ and the associated space $V_{h,0}$ remain unchanged by the cutting operation \eqref{eq:cut}, it is convenient to extend the snapshot solutions of \eqref{eq3} to the background domain $\Omega_0$. In this paper, we consider a \emph{trivial} extension of the snapshot solutions to zero in the inactive regions of the domain $\Omega_0$. We remark that a thorough discussion and investigation of other possible extensions in the context of projection-based reduced order models are given in \cite{Karatzas2020}. Since the latter work concluded that the trivial extensions only slightly affects the eigenvalues decay, i.e., the dimension of the reduced basis, the \emph{zero extension} is the choice we adopt in this work, although other alternatives are also possible. Therefore, the extended version of the full order problem (FOM) \eqref{eq3} reads
\begin{equation}\label{eq3b}
\widehat{{\bf{A}}}(\bm{\mu})\widehat{{\bf{u}}}_h(\bm{\mu}) = \widehat{{\bf{f}}}(\bm{\mu}).
\end{equation}
The size of this extended problem is $\bm{\mu}$-independent, being $\widehat{\bf{A}}$ $\in \mathbb{R}^{\mathcal{N}_{h,0}\times \mathcal{N}_{h,0}}$ and $\widehat{\bf{u}}_h(\bm{\mu}),\widehat{\bf{f}}, \in \mathbb{R}^{\mathcal{N}_{h,0}}$.

\subsection{Reduced basis problem}\label{sec:rb_problem}
In order to solve the FOM problem \eqref{eq3b} using ROM techniques, we seek for a reduced basis ${\bf{V}} \in \mathbb{R}^{\mathcal{N}_{h,0}\times N}$ where $N$ is the reduced space dimension that is ideally chosen to be
$N \ll \mathcal{N}_{h,0}$. Throughout this work we will consider the POD to construct the reduced basis ${\bf{V}}$, while other techniques, such as the Greedy algorithm \cite{QMN_RBspringer} can also be used. The POD will be briefly reviewed in Section \ref{sec:pod}. The Galerkin reduced basis problem reads: find $u_N \in V_N$ such that
\begin{equation}\label{eq6}
a(u_N,v_N;\bm{\mu}) = f(v_N;\bm{\mu}), \qquad \forall v_N \in V_N,
\end{equation}
where $V_N$ denotes the reduced basis space spanned by $\bf{V}$. Thus, using the reduced basis ${\bf{V}}$, the solution $\widehat{{\bf{u}}}_h(\bm{\mu})$ can be approximated as 
\begin{equation}\label{eq5}
\widehat{{\bf{u}}}_h(\bm{\mu}) \approx {\bf{V}} {\bf{u}}_N(\bm{\mu}),
\end{equation}
where ${\bf{u}}_N(\bm{\mu})\in \mathbb{R}^{N}$ is the solution vector of the reduced problem. To obtain a projection-based ROM from \eqref{eq3b}, the residual is enforced to be orthogonal to the subspace $V_N$ such that
\begin{equation}\label{eq5b}
{\bf{V}}^T (\widehat{{\bf{A}}}(\bm{\mu}){\bf{V}}  {\bf{u}}_N(\bm{\mu})-{\widehat{{\bf{f}}}(\bm{\mu})}) = {\bf{0}}.
\end{equation}
Thus, the reduced basis approximation leads to the linear system 
\begin{equation}\label{eq7}
{\bf{A}}_N(\bm{\mu}) {\bf{u}}_N(\bm{\mu}) = {\bf{f}}_N(\bm{\mu}),
\end{equation} 
where the reduced matrices and vectors are given as
\begin{equation}\label{eq8}
{\bf{A}}_N = {\bf{V}}^T\widehat{{\bf{A}}}(\bm{\mu}){\bf{V}} , \qquad {\bf{f}}_N = {\bf{V}}^T\widehat{{\bf{f}}}(\bm{\mu}).
\end{equation}

\begin{rmk}
By introducing the snapshots extension, we construct a reduced basis with global basis functions that are defined on the $\bm{\mu}$-independent background domain $\Omega_0$. The solution of the reduced problem in \ref{eq7} is likewise defined on the background domain. Its values inside the inactive regions are not relevant and can be discarded during the analysis.
\end{rmk}

The reduced problem \eqref{eq7} has size $N$, which makes it suitable for fast online solution given many different parameters $\bm{\mu} \in \mathcal{P}$. Nevertheless, beyond the size of the problem, \eqref{eq7} still requires the assembly of the FOM matrix $\widehat{{\bf{A}}}(\bm{\mu})$ and vector $\widehat{{\bf{f}}}(\bm{\mu})$ for each parameter $\bm{\mu}$, that is in principle expensive. Therefore, a crucial aspect for the efficiency of the ROM is the assumption that both $\widehat{\bf{A}}(\bm{\mu})$ and  $\widehat{\bf{f}}(\bm{\mu})$ depend affinely on the parameters $\bm{\mu}$. Unfortunately due to the fact that $\bm{\mu}$ encodes geometrical parameters, neither of $\widehat{{\bf{A}}}(\bm{\mu})$ and $\widehat{{\bf{f}}}(\bm{\mu})$ can be affinely decomposed as functions of $\bm{\mu}$ in general. Instead we will approximate them as
\begin{equation}\label{eq9}
\widehat{{\bf{A}}}(\bm{\mu}) \approx \sum_{q=1}^{Q_{a}} \theta_q^{a} (\bm{\mu})\widehat{{\bf{A}}}_q, \qquad \widehat{{\bf{f}}}(\bm{\mu})\approx\sum_{q=1}^{Q_f} \theta_q^f(\bm{\mu})\widehat{{\bf{f}}}_q,
\end{equation}
where $\left\lbrace \theta_q^{a} (\bm{\mu}) \right\rbrace_{q=1}^{Q_{a}}$ and $\left\lbrace\theta_q^f (\bm{\mu})\right\rbrace_{q=1}^{Q_f}$ are $\bm{\mu}$-dependent parameter functions and $\widehat{{\bf{A}}}_q \in \mathbb{R}^{\mathcal{N}_{h,0}\times\mathcal{N}_{h,0}}$ and $\widehat{{\bf{f}}}_q \in \mathbb{R}^{\mathcal{N}_{h,0}}$ are $\bm{\mu}$-independent matrices and vectors, respectively. In order to build the approximation \eqref{eq9}, we rely on hyper-reduction techniques, such as the empirical interpolation method (EIM) \cite{Barrault2004,Maday2009} and its discrete variant (DEIM) for vectors and matrices \cite{Chaturantabut2010,Negri2015} to recover the affine dependence. In particular, our goal is to provide an efficient and non-intrusive procedure that is agnostic to parameter-dependent cutting operations for rapid online evaluation of the coefficients $\theta_q^{a}(\bm{\mu}),\theta_q^{f}(\bm{\mu})$ in Equation \eqref{eq9}. In this work, we will exploit interpolation with radial basis functions (RBFs) \cite{Powell1992} for fast online evaluation of the parameter-dependent functions in Equation \eqref{eq9}. This allows to obtain infinite or piecewise smoothness depending on the chosen type of RBFs. The hyper-reduction procedure will be further discussed in Section \ref{sec:deim}. Introducing the affine approximation \eqref{eq9} into \eqref{eq8}, the reduced matrix ${\bf{A}}_N \in \mathbb{R}^{{N}\times{N}}$ and the right-hand side vector  ${\bf{f}}_N \in \mathbb{R}^{{N}}$ are computed for a given parameter $\bm{\mu}$ as
\begin{equation}\label{eq9a}
{\bf{A}}_N(\bm{\mu}) = \sum_{q=1}^{Q_{a}} \theta_q^{a} (\bm{\mu}){\bf{V}}^T\widehat{{\bf{A}}}_q{\bf{V}}, \qquad {\bf{f}}_N(\bm{\mu}) = \sum_{q=1}^{Q_f} \theta_q^f(\bm{\mu}){\bf{V}}^T\widehat{{\bf{f}}}_q.
\end{equation}
The matrices $\left\lbrace {\bf{V}}^T\widehat{{\bf{A}}}_q {\bf{V}}\right\rbrace_{q=1}^{Q_{a}}$ and vectors $\left\lbrace{\bf{V}}^T\widehat{{\bf{f}}}_q \right\rbrace_{q=1}^{Q_f}$ are ${\bm{\mu}}$-independent and can be pre-computed once and stored during the offline phase. During the online phase, we solve the reduced problem \eqref{eq7} for a given value of $\bm{\mu}$. To this end, we first compute the coefficients $\theta_q^{a}(\bm{\mu}),\theta_q^{f}(\bm{\mu})$ in the affine approximation of \eqref{eq9}, assemble the reduced matrix ${\bf{A}}_N(\bm{\mu})$ and vector ${\bf{f}}_N(\bm{\mu})$ as in \eqref{eq9a} and then solve the linear system \eqref{eq7}. Finally, the solution referred to the space $V_{h,0}$ is reconstructed through \eqref{eq5}. It should be highlighted that ideally $Q_f, Q_{a} \ll \mathcal{N}_{h,0}$ and therefore the online assembly in the form of \eqref{eq9a} is inexpensive.

So far we have introduced two approximations with respect to the original problem \eqref{eq3}: the reduced solution and the affine approximations of $\widehat{\bf{A}}(\bm{\mu})$ and  $\widehat{\bf{f}}(\bm{\mu})$. To build such reduction framework there are two crucial steps, namely, the construction of the reduced basis ${\bf{V}}$ and the creation of the approximate affine decompositions in \eqref{eq9}. In this work, the reduced basis ${\bf{V}}$ is constructed by means of the POD and its generation is detailed in Section \ref{sec:pod}. On the other hand, the affine approximations \eqref{eq9} are performed based on the DEIM and will be elaborated in Section \ref{sec:deim}. Nevertheless, the construction of the reduced basis ${\bf{V}}$ on the extended domain $\Omega_0$ required by unfitted domain discretizations results in a manifold that is highly nonlinear on the parameters ${\bm{\mu}}$. The same holds also for the affine approximations of $\widehat{\bf{A}}(\bm{\mu})$ and  $\widehat{\bf{f}}(\bm{\mu})$. In fact, the approximation of a nonlinear solution manifold with a global linear subspace may be accurate only for a very high number of basis functions, which hinders the construction of efficient ROMs. Therefore, we will provide a localization strategy to construct local bases that can be switched online in an efficient manner. We further elaborate the localized reduced basis method in Section \ref{sec:local}.

\subsection{Proper Orthogonal Decomposition}\label{sec:pod}
In this section we briefly review the Proper Orthogonal Decomposition (POD) technique. The reader is referred to \cite{Hesthaven2016,QMN_RBspringer} for a more detailed exposition. Let us first set our notation for the POD approach based on the singular value decomposition (SVD) algorithm that we will use a few times later. The SVD of a matrix ${\bf{S}} \in \mathbb{R}^{m \times n}$ reads:
\begin{equation}\label{eq15}
{\bf{S}}= \mathbb{U}\boldsymbol{\Sigma}\mathbb{Z}^T,
\end{equation}
where the orthogonal matrices $\mathbb{U} \in \mathbb{R}^{{m}\times{m}}$ and $\mathbb{Z} \in \mathbb{R}^{{n}\times{n}}$ contain the left and right singular vectors of ${\bf{S}}$, respectively, and $\boldsymbol{\Sigma} \in \mathbb{R}^{{m}\times{n}}$ is a diagonal matrix containing the positive singular values of ${\bf{S}}$ sorted in descending order. For $m \ge n$, the correlation matrix is defined as $\mathbb{C}={\bf{S}}^T {\bf{S}} \in \mathbb{R}^{n \times n}$. The following eigenvalue problem can be then derived:
\begin{equation}\label{eq17}
\mathbb{C}\bm{\psi}_i = {\sigma}_i^2 \bm{\psi}_i, \ i=1,\dots,r.
\end{equation}
Here, ${\sigma}_i^2$ are the nonzero eigenvalues of the correlation matrix $\mathbb{C}$ sorted in nondecreasing order and $\bm{\psi}_i \in \mathbb{R}^{n \times n}$ are the associated normalized eigenvectors being $r \le n$ the rank of ${\bf{S}}$. The POD basis of dimension $P$ is then obtained from the first $P$ eigenvectors of the correlation matrix as 
\begin{equation}\label{eq18}
\boldsymbol{\zeta}_j = \frac{1}{\sigma_j}{\bf{S}}\bm{\psi}_j, \quad j=1,\dots,P,
\end{equation}
where $\boldsymbol{\zeta}_j \in \mathbb{R}^{m}$. The POD basis is orthonormal by construction and its dimension $P$ can be chosen such that the projection error induced by the POD, that is, the energy captured by the neglected modes, is smaller than a prescribed tolerance $\epsilon_{POD}$ \cite{QMN_RBspringer}. Therefore, it is sufficient to choose $P$ as the smallest integer such that
\begin{equation}\label{eq19}
 1 - \frac{\sum_{i=1}^{P}\sigma_i^2}{\sum_{i=1}^{r}\sigma_i^2} \le \epsilon_{POD}.
\end{equation}

In order to construct a POD basis ${\bf{V}}$ for the approximation in \eqref{eq5} we assume to have a sufficiently fine and properly selected training sample set $\mathcal{P}_{train}=\{\bm{\mu}_1,...,\bm{\mu}_{N_s}\} \subset \mathcal{P}$ of dimension $N_s=\text{dim}(\mathcal{P}_{train})$. Using this sample set, we form the solution snapshots matrix ${\bf{S}}_u \in \mathbb{R}^{\mathcal{N}_{h,0} \times N_s}$
\begin{equation}\label{eq10}
{\bf{S}}_u= [\widehat{\bf{u}}_1,...,\widehat{\bf{u}}_{N_s}],
\end{equation}
where the vectors $\widehat{\bf{u}}_j \in \mathbb{R}^{\mathcal{N}_{h,0}}$ represent the solutions $\widehat{\bf{u}}_h{(\bm{\mu}_j)}$ extended to the background domain $\Omega_0$ for $j=1,...,N_s$. The reduced basis ${\bf{V}}  = [{\boldsymbol{\zeta}}_1,..,{\boldsymbol{\zeta}}_N] \in \mathbb{R}^{\mathcal{N}_{h,0} \times N}$ is then extracted with the POD as shown above in Equations \eqref{eq15}-\eqref{eq19}.

\begin{rmk}
We remark that also other techniques, such as the Greedy algorithm, can be in principle used to construct the reduced basis. In this work we focus on the POD since our interest lies in addressing its shortcomings related to parameterized unfitted geometries, i.e., the slow decay of the singular values of the SVD. Similarly to advection dominated problems with slowly decreasing Kolmogorov n-widths, the effective model reduction is a challenging task. In fact, the proposed framework bears connections to such type of problems that need to be tackled efficiently from the reduction viewpoint.
\end{rmk}

\subsection{Discrete Empirical Interpolation Method}\label{sec:deim}
The key feature to ensure efficiency of the reduced basis method is the affine parametric assumption discussed previously, which allows to decompose the stiffness matrix and right-hand side vector with respect to the parameters $\bm{\mu}$. As a very first step, an affine approximation in the form of Equation \eqref{eq9} is constructed with the discrete empirical interpolation method (DEIM) for matrices and vectors. The reader is further referred to \cite{Negri2015,Wirtz2014} for a detailed overview of this procedure. 

Similar to the solution, the same extension is performed to form the snapshots matrices in the DEIM procedure. Let us consider the stiffness matrix $\widehat{\bf{A}}(\bm{\mu})$ and right-hand side $\widehat{\bf{f}}(\bm{\mu})$ obtained by extending ${\bf{A}}(\bm{\mu})$ and ${\bf{f}}(\bm{\mu})$ to zero inside inactive regions. Following \cite{Negri2015}, the matrices $\widehat{{\bf{A}}}_q$ for $q=1,\dots,Q_a$ and vectors $\widehat{{\bf{f}}}_q$ for $q=1,\dots,Q_f$ in \eqref{eq9} are obtained by applying the POD procedure described in Section \ref{sec:pod}. In the following, we assume that we have a sufficiently fine training sample set  $\mathcal{P}_{train}^d=\{\bm{\mu}_1,...,\bm{\mu}_{N_{s}^d}\} \subset \mathcal{P}$ of dimension $N_{s}^d=\text{dim}(\mathcal{P}_{train}^d)$ and apply the POD to the vectorization $\widehat{{\bf{k}}}(\bm{\mu}) = \text{vec}(\widehat{{\bf{A}}} (\bm{\mu}))$ and to the vectors $\widehat{\bf{f}}(\bm{\mu})$ for each $ \bm{\mu} \in \mathcal{P}_{train}^d$. We denote the snapshots matrices ${{\bf{S}}_a} \in \mathbb{R}^{{{\mathcal{N}^2_{h,0}}} \times N_{s}^d}$ and ${{\bf{S}}_f} \in \mathbb{R}^{\mathcal{N}_{h,0} \times N_{s}^d}$ upon which the POD bases are built as
\begin{equation}\label{eq20a}
{\bf{S}}_a= [\widehat{\bf{k}}_1,...,\widehat{\bf{k}}_{N_{s}^d}], \qquad {\bf{S}}_f= [\widehat{\bf{f}}_1,...,\widehat{\bf{f}}_{N_{s}^d}],
\end{equation}
where $\widehat{\bf{k}}_j=\text{vec}(\widehat{\bf{A}}{(\bm{\mu}_j))}$ and $\widehat{\bf{f}}_j=\widehat{\bf{f}}{(\bm{\mu}_j)}$, for $j=1,...,N_{s}^d$. It should be highlighted that the procedure depends on the sparsity pattern associated to the background domain $\Omega_0$ and only non-zero entries are considered for the implementation. The reader is further referred to \cite{Negri2015} for a detailed discussion on implementation aspects. After performing the SVD to the matrices ${{\bf{S}}_a}$ and ${{\bf{S}}_f}$, the number of affine terms $Q_a$ and $Q_f$ can be determined by prescribing a tolerance $\epsilon_{POD}^d$ and applying the expression \eqref{eq19}.
Note that hyper-reduction in general requires a higher level of accuracy (i.e., $\epsilon_{POD}^d < \epsilon_{POD}$) to obtain reduced basis approximations that are not impaired by the accuracy of the DEIM approximations \cite{QMN_RBspringer}. 

\

\begin{rmk}
In this work, we first consider the DEIM procedure to obtain affine approximations and then construct a reduced basis space with the POD approach. The solution snapshots in Equation \eqref{eq10} are obtained by solving \eqref{eq3b} and exploiting \eqref{eq9}. Such solution is performed offline and its evaluation is nevertheless inexpensive provided that $Q_f, Q_{a} \ll \mathcal{N}_{h,0}$. We further refer to the discussion in \cite{Negri2015,Stewart1990} regarding the preservation of the non-singularity of the stiffness matrix for $ {\bm{\mu}} \in \mathcal{P}$ by the affine approximation. It is worthwhile highlighting that for the sake of offline savings one can construct the snapshots matrices for both the DEIM approximations \eqref{eq20a} and the reduced basis \eqref{eq10} simultaneously in case the reduced space is built with the POD. More details on this alternative option are given in \cite{Negri2015}. Note that our strategy is intended to be also applicable with other techniques, e.g.~the Greedy algorithm, thus we do not further investigate this option here.
\end{rmk}

Now, a strategy to efficiently compute the parameter-dependent coefficients $\theta_q^{a} (\bm{\mu})$, $q=1,\dots,Q_a$ and $\theta_q^{f} (\bm{\mu})$, $q=1,\dots,Q_f$, in \eqref{eq9} is needed. Following the empirical interpolation procedure \cite{Barrault2004}, we set ourselves at using the \emph{magic points} \cite{Maday2009}. Let us focus on the stiffness matrix first. A Greedy algorithm \cite{QMN_RBspringer} that minimizes the interpolation error over the snapshots is used to select a collection of $Q_a$ matrix entries, which we denote as $\mathcal{J}_a$. These entries fulfill exactly the interpolation constraint for the matrix $ \widehat{{\bf{A}}}(\bm{\mu})$, i.e., for each $(i,j) \in \mathcal{J}_a$
 \begin{equation}\label{eq20}
 \sum_{q=1}^{Q_{a}} \theta_q^{a} (\bm{\mu})[\widehat{{\bf{A}}}_q]_{i,j} = [\widehat{{\bf{A}}}(\bm{\mu})]_{i,j}.
\end{equation}
The interpolation constraint for the vector $\widehat{{\bf{f}}}(\bm{\mu})$ reads, for each $i \in \mathcal{J}_f$
 \begin{equation}\label{eq20_f}
 \sum_{q=1}^{Q_{f}} \theta_q^{f} (\bm{\mu})[\widehat{{\bf{f}}}_q]_{i} = [\widehat{{\bf{f}}}(\bm{\mu})]_{i}.
\end{equation}
Note that the well-posedness of the DEIM procedure follows from \cite{Chaturantabut2010} and \cite[Theorem~10.1]{QMN_RBspringer}. In order for Equations \eqref{eq20}-\eqref{eq20_f} to be efficiently computed online, we need to assure that for a given ${\bm{\mu}} \in \mathcal{P}$ their right-hand side can be rapidly computed on the fly. However, the following aspects should be considered:
\begin{itemize}
    \item This evaluation requires to assemble online the matrix and vector associated to a collection of $Q_a$ matrix and $Q_f$ vector entries given a new value of $\bm{\mu}$. Depending on the parametric complexity, i.e.~the number of functions $Q_a,Q_f$, this operation can be costly.
    \item In the context of finite element methods, the PDE operators are in practice assembled employing a \emph{reduced mesh} that benefits from the local support of basis functions \cite{Negri2015}. However, the goal of our work is to provide a strategy that is independent of the underlying discretization.
    \item In the case of unfitted domains, the magic points selected by DEIM may correspond to active or cut functions. The latter are identified during integration and assembly and may change depending on $\bm{\mu}$. Therefore, efficient implementation of this operation requires several intrusive techniques in the high-fidelity approximation and assembly routines \cite{Karatzas2022}. In what follows, we aim to provide a non-intrusive procedure that is completely agnostic to the cutting operations, i.e., it is independent of the number of cut basis functions  and cut domains.
\end{itemize}

To this end, our approach is inspired by interpolation-based ROMs that have been a subject of research in previous studies \cite{Garotta2020,Georgaka2020}. In particular, the online computation of the coefficients $\theta_q^{a} (\bm{\mu})$ and $\theta_q^{f} (\bm{\mu})$ can be made more efficient by: %\cite{Xiao2017,Hijazi2020,Georgaka2020, Garotta2020}. 
\begin{enumerate}
    \item pre-computing the values of $\left\lbrace \theta_q^{a} (\bm{\mu})\right\rbrace _{q=1}^{Q_{a}}$ and $\left\lbrace \theta_q^{f} (\bm{\mu})\right\rbrace _{q=1}^{Q_{f}}$ in \eqref{eq20}-\eqref{eq20_f} during the offline phase for each ${\bm{\mu}} \in \mathcal{P}_{train}^d$,
    \item using such computations to \emph{train} a fast interpolation method during the offline phase,
    \item evaluating efficiently the interpolants during the online phase for any given ${\bm{\mu}} \in \mathcal{P}$.
\end{enumerate}

In the following we adopt radial basis functions (RBFs) for the interpolation \cite{Powell1992}, due to their capability to interpolate scattered data, although the use of other methods is of course also possible. The procedure is identical for both $\left\lbrace \theta_q^{a} (\bm{\mu})\right\rbrace _{q=1}^{Q_{a}}$ and $\left\lbrace \theta_q^{f} (\bm{\mu})\right\rbrace _{q=1}^{Q_{f}}$, therefore we will consider from now on only the first one to keep the exposition concise. During the online phase, the function $\left\lbrace \theta_q^{a} (\bm{\mu})\right\rbrace _{q=1}^{Q_{a}}$ in \eqref{eq9} is approximated as
\begin{equation}\label{eq20_5}
{\theta}^a_q(\bm{\mu})\approx \sum_{j=1}^{N_{s}^d} \omega_{q,j}^{a} \phi_{q,j} (\norm{\bm{\mu}-\bm{\mu}_j}_2), \qquad q=1,\dots,Q_a.
\end{equation}
where $\phi_{q,j}$ denotes the radial basis function associated to the $j$-th center parameter point $\bm{\mu}_j$ and $\norm{\cdot}_2$ represents the Euclidean norm. There are several alternatives for radial basis functions, such as Gaussian, multi-quadratic, and others. In the numerical experiments discussed in Section \ref{sec:examples} we will use cubic radial basis functions. These feature piecewise, higher-order smoothness without spurious oscillations.  It should be noted that the number of interpolation parameter points coincides with the number of training parameter samples $N_{s}^d$. The unknown weights $\omega_{q,j}^{a}$ are computed during the offline phase such that they fulfill the interpolation constraint exactly for ${\bm{\mu}_k} \in \mathcal{P}_{train}^d$
\begin{equation}\label{eq20_4}
\sum_{j=1}^{N_{s}^d} \omega_{q,j}^{a} \phi_{q,j} (\norm{\bm{\mu}_k-\bm{\mu}_j}_2) = \theta_q^{a}(\bm{\mu}_k), \quad k=1,\dots,N_{s}^d, \ q=1,\dots,Q_a.
\end{equation}
We refer the reader to \cite{Powell1992}, where the unique solvability of the underlying linear system is analyzed. Note that depending on the type of radial basis functions, polynomials may be augmented to the above definition to render the problem uniquely solvable \cite{Buhmann2000}. We remark that the condition number of the matrix associated to the RBF problem grows with the number of interpolation points and preconditioning techniques \cite{Beatson1999} or tuning of shape parameters \cite{Fornberg2007} may be needed for large data sets. However, the localization strategy we propose in Section \ref{sec:local} mitigates this effect to some extent, since we partition the data set to construct local approximations.

We recall that in the context of unfitted domain discretizations, the construction of affine approximations on the extended domain $\Omega_0$ results in a manifold that is nonlinear on the parameters ${\bm{\mu}}$. Thus, constructing one global approximation may lead to a high number of affine terms $Q_a,Q_f$. This impedes the overall efficiency of the ROM. The same holds also for the construction of the reduced basis ${\bf{V}}$ and its dimension $N$ discussed in Section \ref{sec:rb_problem}. Therefore, we will introduce a localization strategy to construct accurate, local approximations while containing the dimension of the reduced problem. 

\section{Localization strategy}\label{sec:local}
In the following we will present a strategy to construct efficient ROMs based on localized reduced bases. Our approach is inspired by problems with moving fronts and discontinuities, where local subspaces are constructed for the DEIM and reduced basis approximation \cite{Amsallem2012,Peherstorfer2014,Pagani2018}. This allows the approximation with multiple, smaller subspaces and switching between different local bases in the online phase. Since online evaluations depend only on the dimension of the local bases, one can construct more efficient ROMs compared to a single, global reduced basis approach. The main steps involved in the proposed strategy are: 

\begin{enumerate}
    \item setup a clustering strategy to partition separately the parameters, i.e.~the associated snapshots, for the DEIM and reduced basis approximations,
    \item for each cluster combination train local DEIM approximations and reduced bases during the offline phase as discussed in Section \ref{sec:ROM},
    \item during the online phase, select the cluster with the smallest distance to a given ${\bm{\mu}} \in \mathcal{P}$ and solve the reduced problem \eqref{eq7}.
\end{enumerate}

\subsection{Parameter-based clustering}
Having in mind problems formulated on parameterized unfitted geometries, the question that arises is how to partition snapshots obtained by extension such that the solution can be approximated by a local subspace of sufficiently small dimension. For this purpose we formulate a parameter-based partitioning strategy. 

Let us first present the main idea behind the proposed strategy. We assume that we have $N_c$ partitions that are centered around fixed points, i.e.~centroids, $\bar{\bm{\mu}}_1,\dots,\bar{\bm{\mu}}_{N_c}$ in the parameter space $\mathcal{P}$. We will discuss later how to obtain those. Then let us recall the cut domains $\hat\Omega_i(\bm{\mu})$ for $i=0,\dots,K$ in \eqref{eq:cut} and define the following distance for a given $\bm{\mu} \in \mathcal{P}$
\begin{align}\label{eq:distance}
\mathcal{D}({\bm{\mu}},k) = \max_i \text{dist}(\partial{\hat{\Omega}}_i(\bar{\bm{\mu}}_k)\ \cap \ \partial{\hat{\Omega}}_i({\bm{\mu}})), \quad k=1,\dots,N_c.
\end{align}
In what follows we set ourselves to assign a given parameter vector ${\bm{\mu}}$ to the partition $k$ that minimizes $\mathcal{D}({\bm{\mu}},k) $. This strategy allows us to form partitions comprising unfitted discretizations with similar active and inactive regions. In order to ensure the efficient computation of the above operation, we assume that $\exists C >  0$ such that
\begin{align}\label{eq:assumption}
 \max_i \text{dist}(\partial{\hat{\Omega}}_i(\bar{\bm{\mu}}_k) \ \cap \ \partial{\hat{\Omega}}_i({\bm{\mu}})) \le C \norm{{\bm{\mu}} - {\bar{\bm{\mu}}}_k}_2^2, \quad  k=1,\dots,N_c, \ i=0,\dots,K.
\end{align}
That is, the maximum distance between boundaries is bounded by the distance between the parameters in the Euclidean norm. To this end, the natural choice is to use the parameters as indicator for grouping together snapshots. We remark that in this work we focus on unfitted domain discretizations, where active basis functions may vary for different values of the parameters ${\bm{\mu}}$. In fact, the proposed strategy can be also adapted to other cases where the discontinuity or variability of the solutions stems, for example, from the underlying physical problem. 

 Let us now discuss how to obtain a partition of the parameter space in $N_c$ subspaces as $\mathcal{P} = \bigcup_{k=1}^{N_c} \mathcal{P}_k$. In practice, we consider the partitioning applied to the discrete counterpart of the parameter space. Moreover, in what follows we opt for an automatic partitioning with the \emph{k-means} clustering algorithm \cite{Likas2003}, although other partitioning techniques are also possible \cite{Eftang2010a,Eftang2011,Haasdonk2011}. Then the $i$-th snapshot, i.e., the $i$-th column of the matrices \eqref{eq10} and \eqref{eq20a} is assigned to a specific cluster $k$ if ${\bm{\mu}}_i \in \mathcal{P}_k$. Note that neighboring snapshots can be added to each cluster to obtain overlapping clusters with smooth transitions from one cluster to another. In the numerical experiments of Section \ref{sec:examples} we will not consider overlaps between clusters, although this is in principle possible \cite{Amsallem2012}. In the following, we will consider a separate partitioning for the DEIM approximation, as this allows to render the dimension of the local bases associated to the DEIM approximations independent of the dimension of the local reduced bases obtained with the POD.

\subsection{Offline phase}
We now present in detail the offline steps involved in the localization strategy. First, we consider the DEIM approximation and the parameter set $\mathcal{P}_{train}^d=\{\bm{\mu}_1,...,\bm{\mu}_{N_{s}^d}\} \subset \mathcal{P}$ introduced in Section \ref{sec:deim}. The first step in the offline phase is to partition the matrix $\mathcal{P}_{train}^d$ into $N_{c}^d$ submatrices corresponding to subregions $\mathcal{P}_k^d \subset \mathcal{P}$ for $k=1,\dots,N_{c}^d$. The \emph{k-means} algorithm starts by choosing random cluster centers (i.e.~centroids) $\{\bar{\bm{\mu}}_{k}^d\}_{k=1}^{N_{c}^d}$. Then the partition is performed such that $\forall \bm{\mu} \in \mathcal{P}_{train}^d$
\begin{equation}\label{eq21_a}
{\mathcal{P}_k^d}= \{{\bm{\mu}} \ \vert\ \text{if} \ \arg \min_{i}  \norm{{\bm{\mu}} - {\bar{\bm{\mu}}}_{i}^d}_2^2=k\}, \quad k=1,...,N_{c}^d.
\end{equation}
The cluster centroids are updated iteratively until the algorithm converges such that
\begin{equation}\label{eq21_b}
{\bar{\bm{\mu}}_k^d} = \frac{1}{\abs{\mathcal{P}_k^d}}\sum_{{\bm{\mu}} \in {\mathcal{P}_k^d}} {\bm{\mu}}, \quad k=1,...,N_{c}^d. \ 
\end{equation}
The k-means clustering minimizes the distance between each parameter vector and the cluster's centroid with respect to the Euclidean norm $\norm{\cdot}_2$. We refer the reader to \cite[Algorithm~5]{Amsallem2012} for a detailed overview of the k-means algorithm. The snapshots are then grouped into the same clusters as their corresponding parameters following the assumption in \eqref{eq:assumption}. Thereafter, the DEIM procedure and RBF interpolation described in Section \ref{sec:deim} is performed separately for each cluster. In the numerical experiments, we will adopt the same number of clusters for the DEIM approximation of the stiffness matrix and right-hand side, although in principle this could be chosen differently. The offline localization procedure for the DEIM approximations is presented in Algorithm \ref{algo1}.

\makeatletter

\def\algbackskip{\hskip-\ALG@thistlm}
\makeatother
 \begin{algorithm}
    \caption{Localized DEIM procedure}\label{algo1}
    \begin{algorithmic}[1]
    \Procedure{[$\{{\bar{\bm{\mu}}_{k}^d}\},\{\widehat{\bf{A}}_q^k\},\{\widehat{\bf{f}}_q^k\},\{\omega^a_{q,j,k}\},\{\omega^f_{q,j,k}$\}] = OFFLINE}{$P_{train}^d, N_{c}^d, \epsilon_{POD}^{d}$}
    \State $\{{\bar{\bm{\mu}}_{k}^d}\}_{k=1}^{N_{c}^d},\{{\mathcal{P}_1^d,.., \mathcal{P}^d_{N_{c}^d}}\} \gets \textit{k-means clustering } (\mathcal{P}_{train}^d,N_{c}^d)$
    \State \emph{Local DEIM basis functions, indices and interpolation weights:}
    \For{$k = 1,...,N_{c}^d$}    
     \For{${\bm{\mu}} \in \mathcal{P}_k^d$} 
    \State \emph{Compute} $\widehat{\bf{A}}{(\bm{\mu})}, \widehat{\bf{f}}{(\bm{\mu})}$ \emph{with full order model}
    \State ${\bf{S}}_a^{k} = [{\bf{S}}_a^{k}, \widehat{\bf{A}}{(\bm{\mu})}]; \ \ 
     {\bf{S}}_f^{k} = [{\bf{S}}_f^{k}, \widehat{\bf{f}}{(\bm{\mu})}]$
    \EndFor
    \State $\widehat{\bf{A}}_q^k \gets
    \textit{POD}({\bf{S}}_a^{k},\epsilon_{POD}^{d}); \ \ \mathcal{J}_{\alpha}^k 
    \gets \textit{DEIM-indices}(\widehat{\bf{A}}_q^k)$ ,\ $q=1,\dots,Q_a$
    \State $\widehat{\bf{f}}_q^k \gets\textit{POD}({\bf{S}}_f^{k},\epsilon_{POD}^{d}); \ \ 
    \mathcal{J}_{f}^k\gets 
    \textit{DEIM-indices}(\widehat{\bf{f}}_q^k)$, ,\ $q=1,\dots,Q_f$
    \State $\omega^a_{q,j,k} \gets \textit{RBF}(P_k^d,{\bf{S}}_a^{k},\widehat{\bf{A}}_q^k,\mathcal{J}_{\alpha}^k)$ \eqref{eq20_4}, \ $q=1,\dots,Q_a$,\ $j=1,\dots,\text{dim}(\mathcal{P}_{train}^d)$
     \State $\omega^f_{q,j,k} \gets \textit{RBF}(P_k^d,{\bf{S}}_f^{k},\widehat{\bf{f}}_q^k,\mathcal{J}_{f}^k)$ \eqref{eq20_4}, \ $q=1,\dots,Q_f$,\ $j=1,\dots,\text{dim}(\mathcal{P}_{train}^d)$
    \EndFor
    \EndProcedure 
    \end{algorithmic}
    \end{algorithm}

\noindent Once the local DEIM approximations are constructed, the next step in the offline phase is to construct local reduced bases. We consider the parameter set $\mathcal{P}_{train}=\{\bm{\mu}_1,...,\bm{\mu}_{N_s}\} \subset \mathcal{P}$ and the solution snapshots matrix defined in Equation \eqref{eq10}. The partitioning is performed in the same manner as before, that is, we seek $N_c$ partitions corresponding to subregions $\mathcal{P}_k \subset \mathcal{P}$, $k=1,..,N_c$. Then, the k-means algorithm initially selects random cluster centroids $\{{\bar{\bm{\mu}}}_k\}_{k=1}^{N_{c}}$ that are updated iteratively following the steps in Equations \eqref{eq21_a} and \eqref{eq21_b}. In order to evaluate \eqref{eq9}, we select the local DEIM approximation by minimizing the distance between a given parameter $\bm{\mu} \in \mathcal{P}_{train}$ and the clusters' centroids such that: 
\begin{equation}\label{eq22}
l = \arg \min_i  \norm{{\bm{\mu}} - {\bar{\bm{\mu}}}_{i}^d }_2^2, \ i=1,\dots,N_{c}^d.
\end{equation}

\makeatletter
\def\algbackskip{\hskip-\ALG@thistlm}
\makeatother
 \begin{algorithm}
    \caption{Localized reduced basis procedure}\label{algo2}
    \begin{algorithmic}[1]
    \Procedure{[${\{\bar{\bm{\mu}}}_k\},\{{\bf{A}}_{N,q}^{k,\bar{k}}\},\{{\bf{f}}_{N,q}^{k,\bar{k}}\}$] = OFFLINE}{$\text{DEIM ARRAYS},{\mathcal{P}_{train}}, N_c, \epsilon_{POD}$}
      \State $\{\bar{{\bm{\mu}}}_k\}_{k=1}^{N_{c}},\{{\mathcal{P}_1,.., \mathcal{P}_{N_{c}}}\} \gets \textit{k-means clustering } (\mathcal{P}_{train},N_{c})$
    \State \emph{Local reduced basis functions and reduced arrays}:
    \For{$k = 1,...,N_c$}    
    \For{${\bm{\mu}} \in \mathcal{P}_k$} 
    \State \emph{Full order arrays}:
    \State $l = \arg \min_i  \norm{{\bm{\mu}} - \bar{{\bm{\mu}}}_{i}^d}_2^2$, \ $i=1,\dots,N_{c}^d$
    \State $\{\{\widehat{\bf{A}}_q^l\}_{q=1}^{Q_a},\{\theta_{q,l}^a(\bm{\mu})\}_{q=1}^{Q_a}\} \gets \textit{assemble} \ \widehat{\bf{A}}{(\bm{\mu})}$ \textit{using} \eqref{eq9} \textit{and}  \eqref{eq20_5}
       \State $\{\{\widehat{\bf{f}}_q^l\}_{q=1}^{Q_f},\{\theta_{q,l}^f(\bm{\mu})\}_{q=1}^{Q_f} \} \gets \textit{assemble} \  \widehat{\bf{f}}{(\bm{\mu})}$ \textit{using} \eqref{eq9} \textit{and}  \eqref{eq20_5}
    \State $\widehat{\bf{u}}_h{(\bm{\mu})} \gets \textit{solve FOM in } \eqref{eq3b}$ 
    \State \emph{Solution snapshots}:
    \State ${\bf{S}}_u^{k} = [{\bf{S}}_u^{k}, \widehat{\bf{u}}_h{(\bm{\mu})}]$
    \EndFor
    \State ${\bf{V}}_k \gets \textit{POD}({\bf{S}}_u^{k},\epsilon_{POD});$
    %\State \emph{Local reduced arrays}:
    \For{$\bar{k} = 1,...,N_{c}^d$}  
    \State $\{\{{\bf{A}}_{N,q}^{k,\bar{k}}\}_{q=1}^{Q_a}, \{{\bf{f}}_{N,q}^{k,\bar{k}}\}_{q=1}^{Q_f}\} \gets \textit{projection of full order arrays onto} \ {\bf{V}}_k $ \eqref{eq9a}
    \EndFor
    \EndFor
    \EndProcedure
    \end{algorithmic}
    \end{algorithm}

\noindent Once the snapshots matrix is constructed, each snapshot is assigned to the same cluster as its respective parameter. Then we construct local reduced bases and project all possible combinations of full order arrays obtained by DEIM onto each local subspace as described in Section \ref{sec:ROM}. Algorithm \ref{algo2} presents the offline localization procedure to construct the reduced bases. Note that the input DEIM arrays refer to the output of Algorithm \ref{algo1}.

Since the number of clusters has to be chosen in advance, the \emph{k-means variance} can be considered to choose the optimal number during the offline phase. In this work, we adopt this criterion for the parameter vectors. The \emph{k-means variance} reads
\begin{equation}\label{eq23}
{\bf{\mathcal{V}}} = \sum_{k=1}^{N_c} \sum_{{\bm{\mu}} \in \mathcal{P}_k} \norm{{\bm{\mu}} - \bar{{\bm{\mu}}}_k}_2^2.
\end{equation}
The same criterion holds also for the DEIM approximations by replacing the sum over $N_c^d$ clusters and evaluating the Euclidean distance to the centroids $\bar{{\bm{\mu}}}_k^d$ for $\bm{\mu} \in \mathcal{P}_k^d$. As the number of clusters increases, the variance is expected to decrease rapidly until it reaches a plateau. One can choose the number of clusters based on this \emph{elbowing} effect of the variance, that is, the smallest integer for which a transition from a steep slope to a plateau occurs. We further refer the reader to \cite{Hess2019} for a thorough discussion on this criterion. It is worthwhile noting that the optimal choice of clusters depends to some extent on the given problem at hand, that is, the targeted accuracy and computational speedup. We now summarize the steps of the offline phase as follows:
\\
\begin{enumerate}
    \item We partition the parameters $\mathcal{P}_{train}^d$ into clusters for a given number of clusters $N_{c}^d$. 
    \item We construct the snapshots matrices ${\bf{S}}_a^k$, ${\bf{S}}_f^k$ with $ k=1,
    \dots,N_{c}^d$ for the DEIM approximations by solving the FOM. Each snapshot (i.e.~column of ${\bf{S}}_a^k$, ${\bf{S}}_f^k$) is assigned to the same cluster as its corresponding parameter.
    \item We construct local DEIM approximations and store the basis functions and interpolation weights for each cluster.
    \item We partition the parameters $\mathcal{P}_{train}$ for a given number of cluster $N_c$. 
     \item We construct the solution snapshots matrix ${\bf{S}}_u^k$ with $ k=1,\dots
    \dots,N_{c}$ by solving the problem \eqref{eq3b} exploiting the affine form of Equation \eqref{eq9}. Each snapshot (i.e.~column of ${\bf{S}}_u^k$) is assigned to the same cluster as its corresponding parameter.
    \item We construct local reduced bases for each cluster applying the POD technique. 
    \item We construct local ROMs for all cluster combinations, that is, by projecting each local DEIM approximation onto each local reduced basis space.
\end{enumerate}

\subsection{Online phase}
In the online phase, for a given parameter $\bm{\mu} \in \mathcal{P}$, we switch between DEIM approximations and local bases such that the distance to the respective cluster centroid is minimized following Equation \eqref{eq22}. In this way, the online evaluation depends only on the dimension of the local bases. Note that in case the reduced basis is constructed with the Greedy algorithm, an additional transformation of the basis is done as discussed in \cite{Maday2013}.  It should be highlighted that since both DEIM and reduced basis approximations are localized, during the offline phase we project and store all possible cluster combinations when pre-computing the matrices $\left\lbrace {\bf{V}}^T\widehat{{\bf{A}}}_q {\bf{V}}\right\rbrace_{q=1}^{Q_{a}}$ and vectors $\left\lbrace{\bf{V}}^T\widehat{{\bf{f}}}_q \right\rbrace_{q=1}^{Q_f}$ in \eqref{eq9a}. Then, during the online phase we pick the operators associated to the selected clusters. Thus, the primary goal in constructing local, low-dimensional reduced bases is to reduce the online computational cost at the price, however, of additional offline effort associated to constructing and storing multiple reduced bases. We remark that the online cost might vary between different clusters depending on the dimension of the associated local basis and the number of affine terms. In the presented algorithms and numerical experiments, the respective dimensions $N,Q_a,Q_f$ refer to the number of local functions in the selected cluster. The online phase is given in Algorithm \ref{algo3}. Note that the input ROM arrays refer to the output of Algorithm \ref{algo2} and RBF arrays to the output of the RBF interpolation in Algorithm \ref{algo1}. 

\makeatletter
\def\algbackskip{\hskip-\ALG@thistlm}
\makeatother
 \begin{algorithm}
    \caption{Online phase}\label{algo3}
    \begin{algorithmic}[1]
   \Procedure{[${\bf{u}}_N$] = ONLINE}{\text{ROM ARRAYS},\ \text{RBF ARRAYS}, $\{{\bar{\bm{\mu}}_{k}^d}\}_{k=1}^{N_c^d},{\bm{\mu}}$}
    \State $l = \arg \min_i  \norm{{\bm{\mu}} - \bar{{\bm{\mu}}}_{i}^d}_2^2, \ i=1,\dots,N_{c}^d$
    \State $m = \arg \min_j  \norm{{\bm{\mu}} - \bar{{\bm{\mu}}}_{j}}_2^2 \ \ j=1,\dots,N_{c}$
    \State \emph{Reduced order arrays}:
    \State \emph{compute} $\theta_{q,l}^{a}({\bm{\mu}}),\ q=1,\dots,Q_a \ \eqref{eq20_5}$
    \State \emph{compute} $\theta_{q,l}^{f}({\bm{\mu}}),\ q=1,\dots,Q_f \ \eqref{eq20_5}$
    \State ${\bf{A}}_{N}({\bm{\mu}})=\sum_{q=1}^{Q_a} \theta_{q,l}^{a} (\bm{\mu}){\bf{A}}_{N,q}^{l,m}; \quad {\bf{f}}_{N}({\bm{\mu}})=\sum_{q=1}^{Q_f} \theta_{q,l}^{f} (\bm{\mu}){\bf{f}}_{N,q}^{l,m}$
    \State ${\bf{u}}_N \gets \textit{solve reduced linear system in} \ \eqref{eq7}$ 
    \EndProcedure
    \end{algorithmic}
    \end{algorithm}

\section{Numerical experiments}\label{sec:examples}
In this section we present some numerical experiments for the Poisson and linear elasticity problems to assess the performance of the proposed methodology in constructing efficient reduced order models for PDEs defined on parameterized unfitted geometries. As discussed in Section \ref{sec:unfitted}, we make use of spline discretizations that are built upon a Cartesian mesh (see Remark \ref{remark1}). However, the method is agnostic to the underlying discretization and perfectly suitable for other choices. We further refer the reader to \cite{Hughes2005,Cottrell2007} and references therein for a detailed review on splines and isogeometric analysis in general as well as to our previous works on trimming using isogeometric analysis \cite{Coradello2020,Coradello2021,Buffa2022}. The results have been obtained using the open-source Octave/Matlab isogeometric package \emph{GeoPDEs} \cite{Vazquez2016} in combination with the open-source library \emph{redbKIT} \cite{redbKIT} and the re-parameterization tool for integration of trimmed geometries presented in \cite{Antolin2019,Wei2021}. It is worthwhile remarking that we adopt a simple diagonal pre-conditioning to limit the consequences of trimming on the condition number while a more detailed discussion is provided in \cite{Deprenter2017}. Unless stated otherwise, we approximate the parameter-dependent coefficients in \eqref{eq9} using cubic RBFs to compute \eqref{eq20_5} and employ Latin hypercube sampling \cite{McKay1979} to select the parameters for our training and test sets. Table \ref{tab:parameters} summarizes the notation defined in the previous sections and used in the numerical experiments.

\begin{table}[h!]
 	\centering
 	\caption{Overview of parameters employed for the numerical experiments}\label{tab:parameters}
 	\begin{tabular}{ll} \hline 
 		$N$  & dimension of reduced basis  \\ 
 		$Q_{a}$  & number of affine terms for stiffness matrix  \\ 
 		$Q_f$ & number of affine terms for right-hand side vector  \\
 		$\mathcal{N}_{h,0}$  &  degrees of freedom of the background domain \\
 		$N_{s}^d$  & dimension of training sample for DEIM approximations  \\ 
 		$N_s$ & dimension of training sample for POD  \\ 
 		$N_c$  & number of clusters for reduced basis  \\ 
 		$N_{c}^d$  & number of clusters for DEIM approximations  \\ 
 		$N_{t}$  & dimension of test sample for error analysis\\\hline
 	\end{tabular}
 \end{table}

\subsection{ The Poisson problem}
Let us first consider the Poisson equation on a parameterized domain. The continuous formulation of the problem reads in strong form: for any $\bm{\mu} \in \mathcal{P}$, find $u \in H_{0,{\Gamma_D}}^1(\Omega(\bm{\mu}))$ such that
\begin{equation}\label{eq1a}
\begin{cases}
- \Delta u &= {f} \ \quad \text{in} \ \Omega({\bm{\mu}}) \\
 \quad \ \ u &= 0 \ \quad \text{on} \ \Gamma_D({\bm{\mu}}) \\
 \quad \ \displaystyle \frac{\partial{u}}{\partial{\boldsymbol{n}}} &= 0 \ \quad \text{on} \ \Gamma_{N}({\bm{\mu}}), 
\end{cases}
\end{equation}
where $\Gamma_D({\bm{\mu}})\subset \partial{\Omega({\bm{\mu}})}\cap\partial{\Omega_0}$ denotes the Dirichlet part of the boundary. We define $H_{0,{\Gamma_D}}^1(\Omega(\bm{\mu})) \subset H^1(\Omega(\bm{\mu}))$ as the subspace of  $H^1(\Omega(\bm{\mu}))$ such that functions vanish on the Dirichlet boundary. The Neumann part of the boundary is $\Gamma_N({\bm{\mu}})$ and it holds that $\overline{\Gamma}_D({\bm{\mu}})\cup\overline{\Gamma}_N({\bm{\mu}}) = \partial{\Omega}({\bm{\mu}})$ and $\Gamma_D({\bm{\mu}}) \cap \Gamma_N({\bm{\mu}})= \varnothing$.  Furthermore, ${f} \in L^2(\Omega(\bm{\mu}))$ is the source term and  $\boldsymbol{n}$ the outward unit normal to the boundary $\partial{\Omega}({\bm{\mu}})$. For simplicity of exposition, we assumed above homogeneous Dirichlet and Neumann boundary conditions without loss of generality. We can now write the discrete weak formulation of the parameterized problem as: find $u_h \in V_h$ such that 
\begin{equation}\label{eq2a}
a(u_h,v_h;\bm{\mu}) = f(v_h;\bm{\mu}), \qquad \forall v_h \in V_h,
\end{equation}
where $V_h \subset H_{0,\Gamma_D}^1(\Omega(\bm{\mu}))$ is the finite-dimensional subspace spanned by a B-spline basis. 
The associated parameterized bilinear form $a(\cdot,\cdot;\bm{\mu})$ and the linear functional $f(\cdot;\bm{\mu})$ read:
\begin{equation}\label{eq3a}
\begin{aligned}
a(u_h,v_h;\bm{\mu}) &= \int_{{\Omega}(\bm{\mu})} \nabla{u}_h \cdot \nabla{v}_h \, \textrm{d}\Omega, \\
f(v_h;\bm{\mu}) &= \int_{{\Omega}(\bm{\mu})} {f} v_h \, \textrm{d}\Omega.
\end{aligned}
\end{equation}

\subsubsection{Square with circular hole: 1D geometrical parameterization}
\begin{figure}[!htb]
	\centering
	\includegraphics[width=1.0\textwidth]{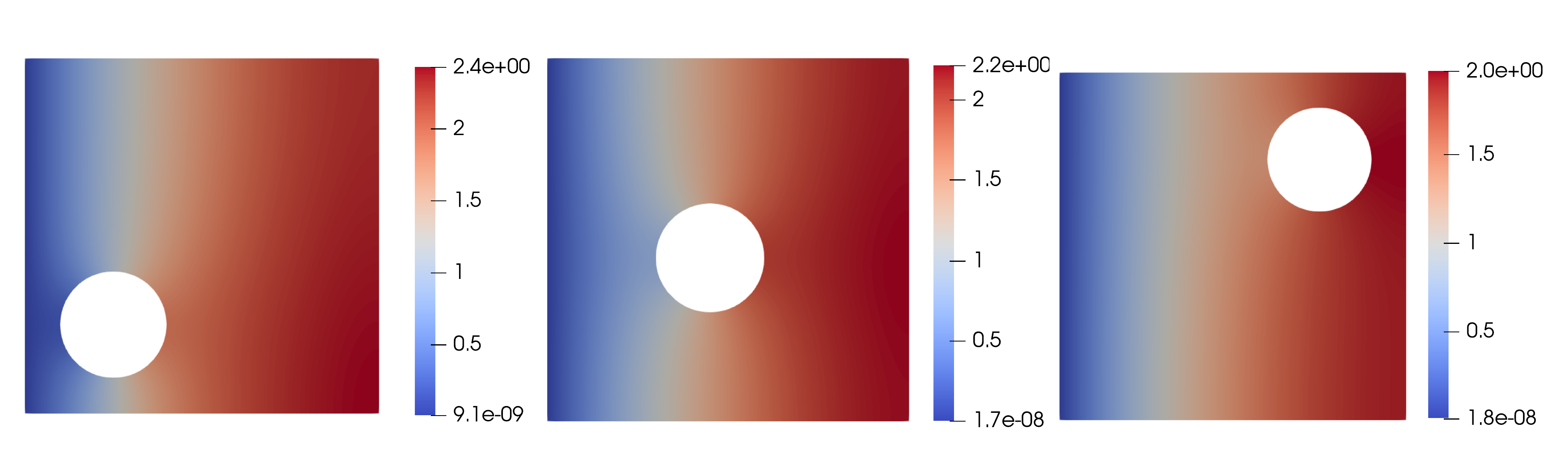} \\
	\caption{Example 6.1.1: Exemplary solution snapshots for  $\mu=[0.5,0.9,1.5]$.}\label{fig:FOM_sol}
\end{figure}

The first example we consider is a two-dimensional problem with a single geometrical parameter. The model is defined on a rectangular domain $\Omega_0 = (0, 2)^2$, which is trimmed by a circular curve of radius $R=0.3$. The trimmed domain $\Omega(\mu)$ is parameter-dependent, where $\mu \in [0.5, 1.5]$ is a parameter representing the coordinates of the center of the circular hole. The hole is centered at $(\mu,\mu)$ and therefore moves along one diagonal of the square $\Omega_0$.
Homogeneous Dirichlet boundary conditions are imposed on the left boundary of the domain and a constant source term is set as ${f}=1$. The geometry is discretized with cubic $C^2$-continuous B-splines using a mesh with $32$ elements per direction over a Cartesian grid, resulting in $\mathcal{N}_{h,0}=1225$ degrees of freedom. We remark that the radius of the hole is fixed and the trimming causes the number of active basis functions to change slightly for the problem at hand. However, different basis functions are active depending on the location of the circular hole. The solution of the FOM is depicted in Figure \ref{fig:FOM_sol} for three different values of the parameter $\mu$. Although we have chosen a simplified setup, the solution of the problem varies significantly for different values of the parameter depending on the location of the trimmed region. Therefore, constructing an efficient ROM for this problem poses challenges to traditional reduced basis methods.

\begin{figure}[!h]
     \begin{subfigure}[b]{0.49\textwidth}
         \centering
         \includegraphics[width=0.9\textwidth]{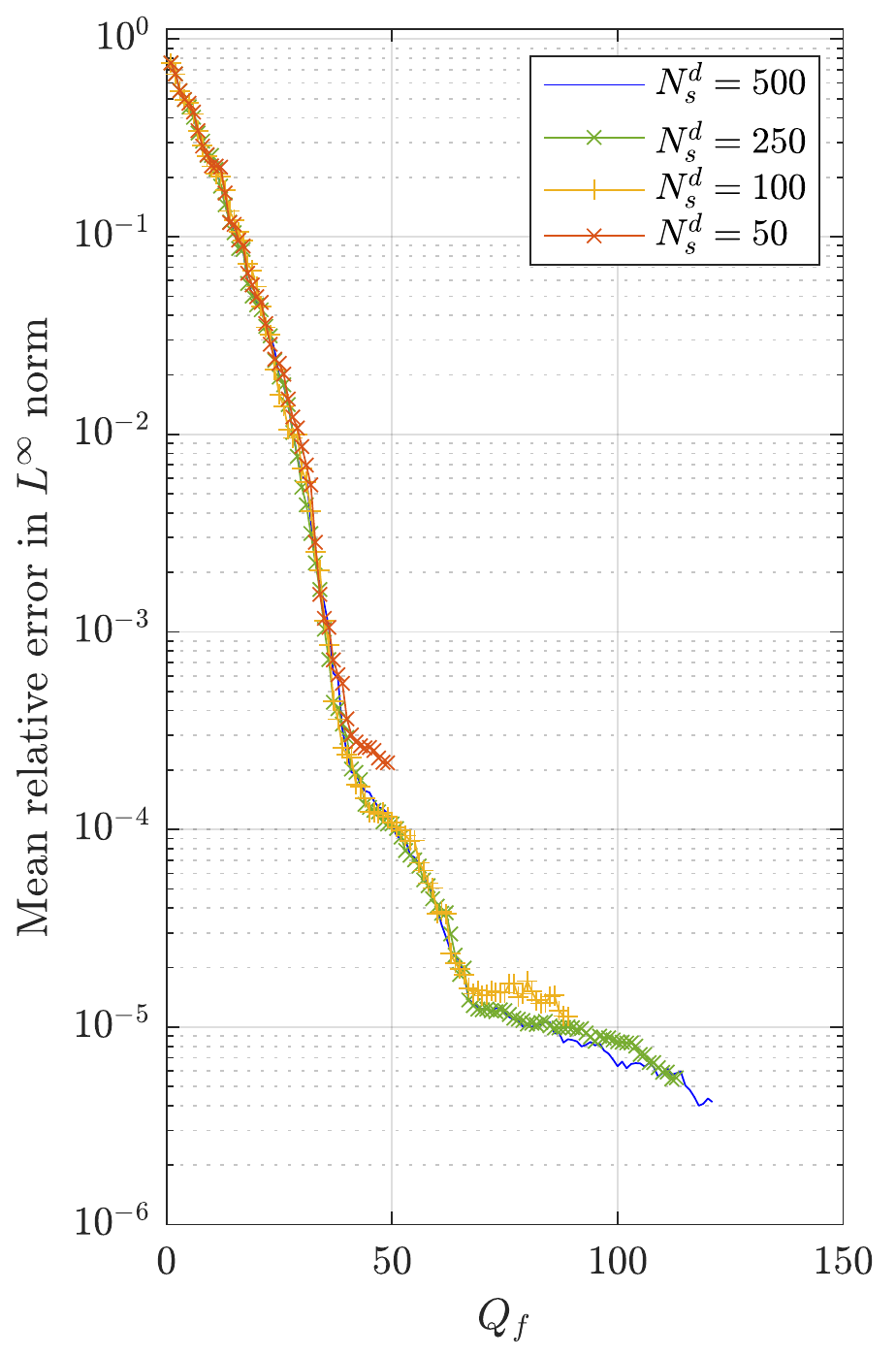}
         \caption{Right-hand side}
         \label{fig:DEIM_error_1}
     \end{subfigure}
     \begin{subfigure}[b]{0.49\textwidth}
         \centering
         \includegraphics[width=0.9\textwidth]{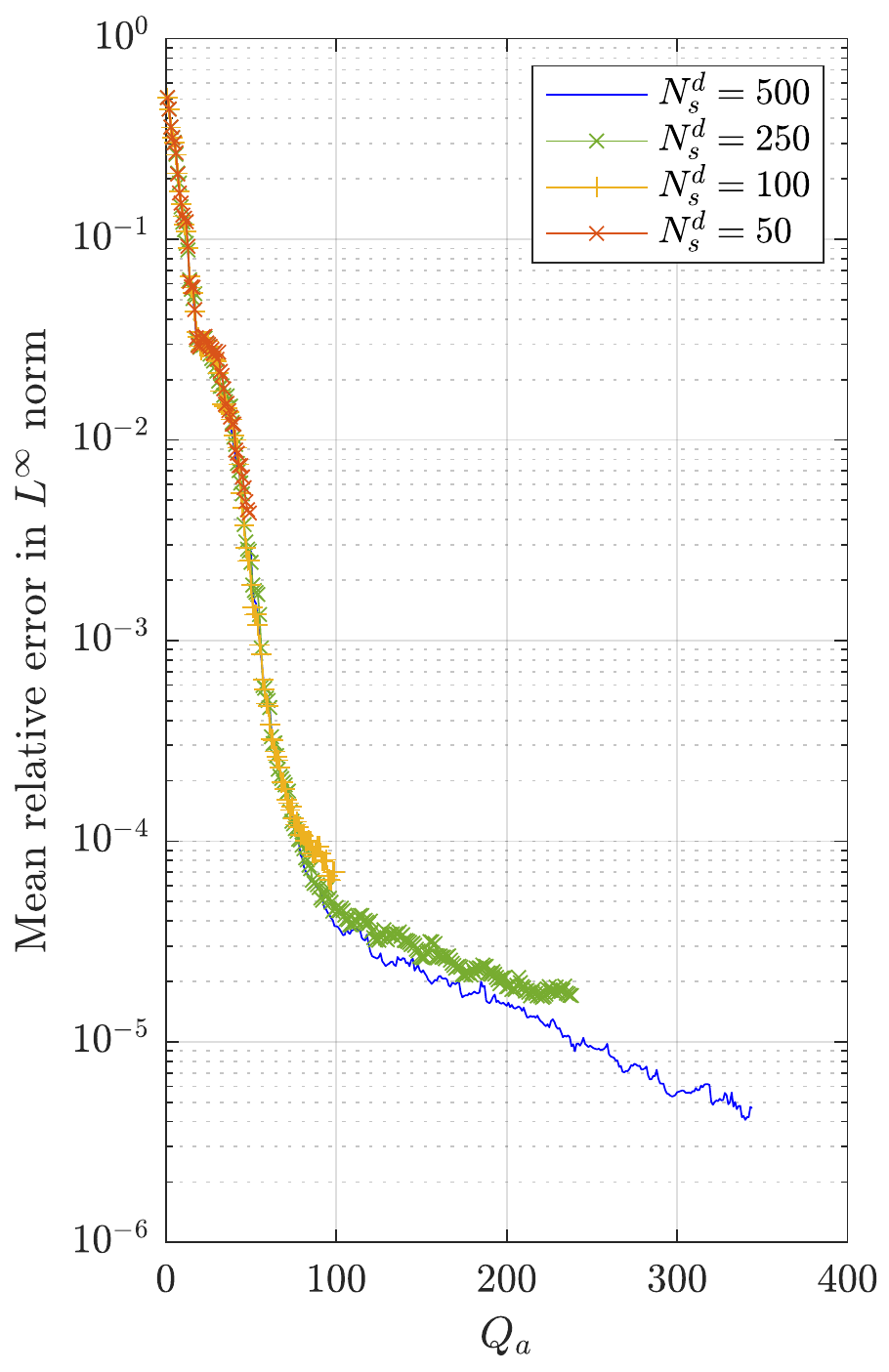}
         \caption{Stiffness matrix}
         \label{fig:DEIM_error_2}
     \end{subfigure}
        \caption{Example 6.1.1: Error decay of global DEIM approximations in $L^{\infty}$-error norm for right-hand side vector (a) and  matrix (b) based on POD tolerance of $\epsilon_{POD}^d= 10^{-7}$.}\label{fig:DEIM_error}
\end{figure}

Let us first consider the standard case where a global ROM is constructed, in order to show that such an approach is not feasible for our problem. Given the above parameterization, the extended stiffness matrix $\widehat{\mathbf{A}}(\mu)$ and right-hand side vector $\widehat{\mathbf{f}}(\mu)$ depend on the geometric parameter in a nonaffine way. Therefore, they can be approximated by DEIM to obtain an affine expansion of the given matrix and vector, respectively, as discussed in Section \ref{sec:deim}. Figure \ref{fig:DEIM_error} depicts the error decay of the DEIM approximations for varying dimension of the training set used to compute the POD basis, namely $N_{s}^d=[50,100,250,500]$. The error analysis is performed based on a test set of dimension $N_{t}=100$ by computing the mean relative error in the $L^{\infty}$ norm between the full order operators and the DEIM approximations, while the coefficients $\theta_q^{a} (\bm{\mu})$, $q=1,\dots,Q_a$ and $\theta_q^{f} (\bm{\mu})$, $q=1,\dots,Q_f$ are computed exactly using Equations \eqref{eq20}-\eqref{eq20_f}. Observing the results in Figure \ref{fig:DEIM_error}, it is evident that the training set needs to be sufficiently rich (i.e. $N_{s}^d \ge 100$) to achieve an accuracy of the order $10^{-5}$. The results indicate that a large number of basis functions $Q_a$ and $Q_f$ needs to be selected to achieve accuracy of the ROM that is not impaired by the error of the DEIM approximation. 
In fact, a large number of DEIM terms is known to reduce significantly the efficiency of the ROM, that is, the online cost within the RB framework \cite{Eftang2012}. For the problem at hand, this motivates the localized strategy introduced in Section \ref{sec:local} to contain the number of selected basis functions.

\begin{figure}[!h]
     \begin{subfigure}[b]{0.49\textwidth}
         \centering
         \includegraphics[width=0.9\textwidth]{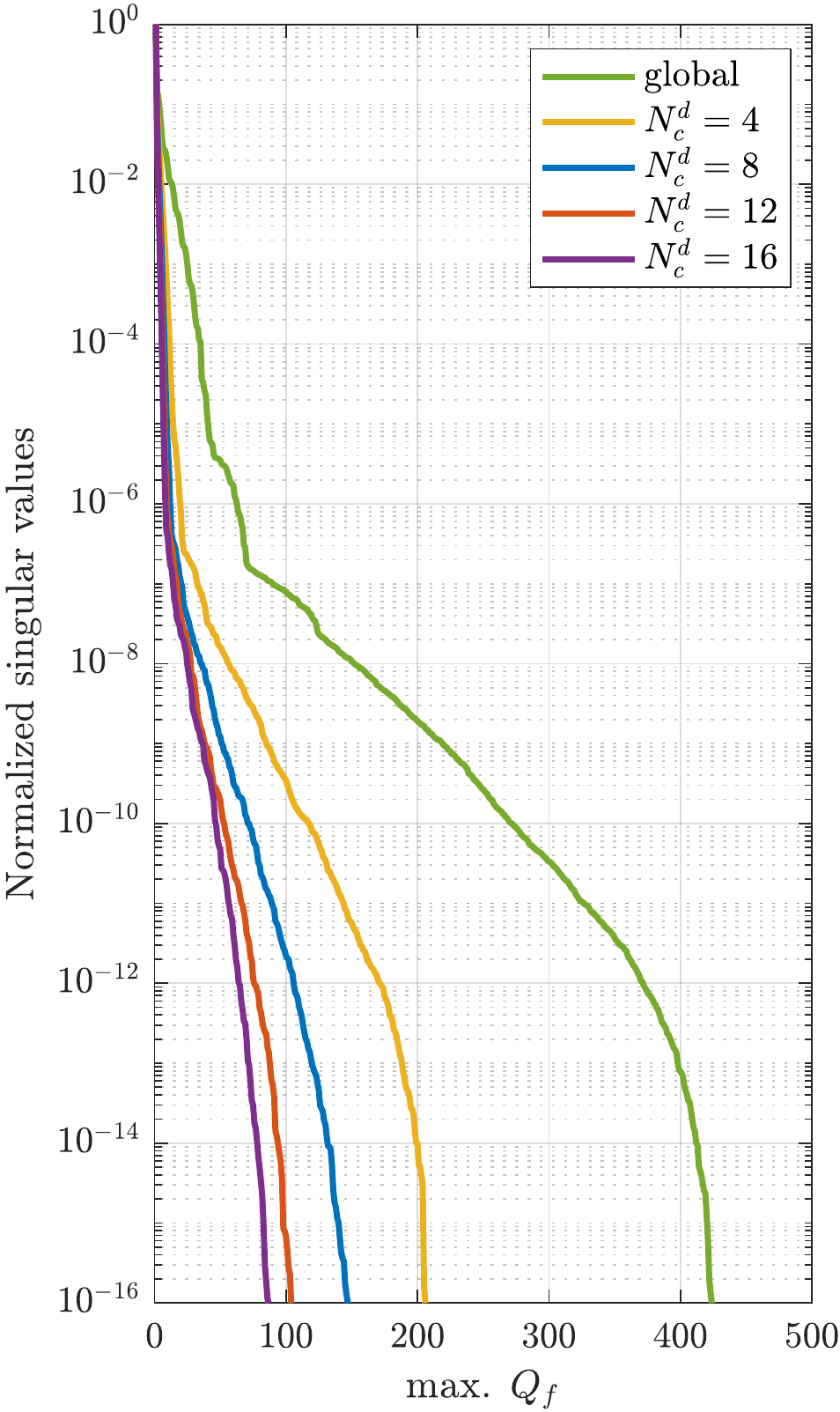}
         \caption{Right-hand side}
         \label{fig:PODDEIM_1}
     \end{subfigure}
     \begin{subfigure}[b]{0.49\textwidth}
         \centering
         \includegraphics[width=0.9\textwidth]{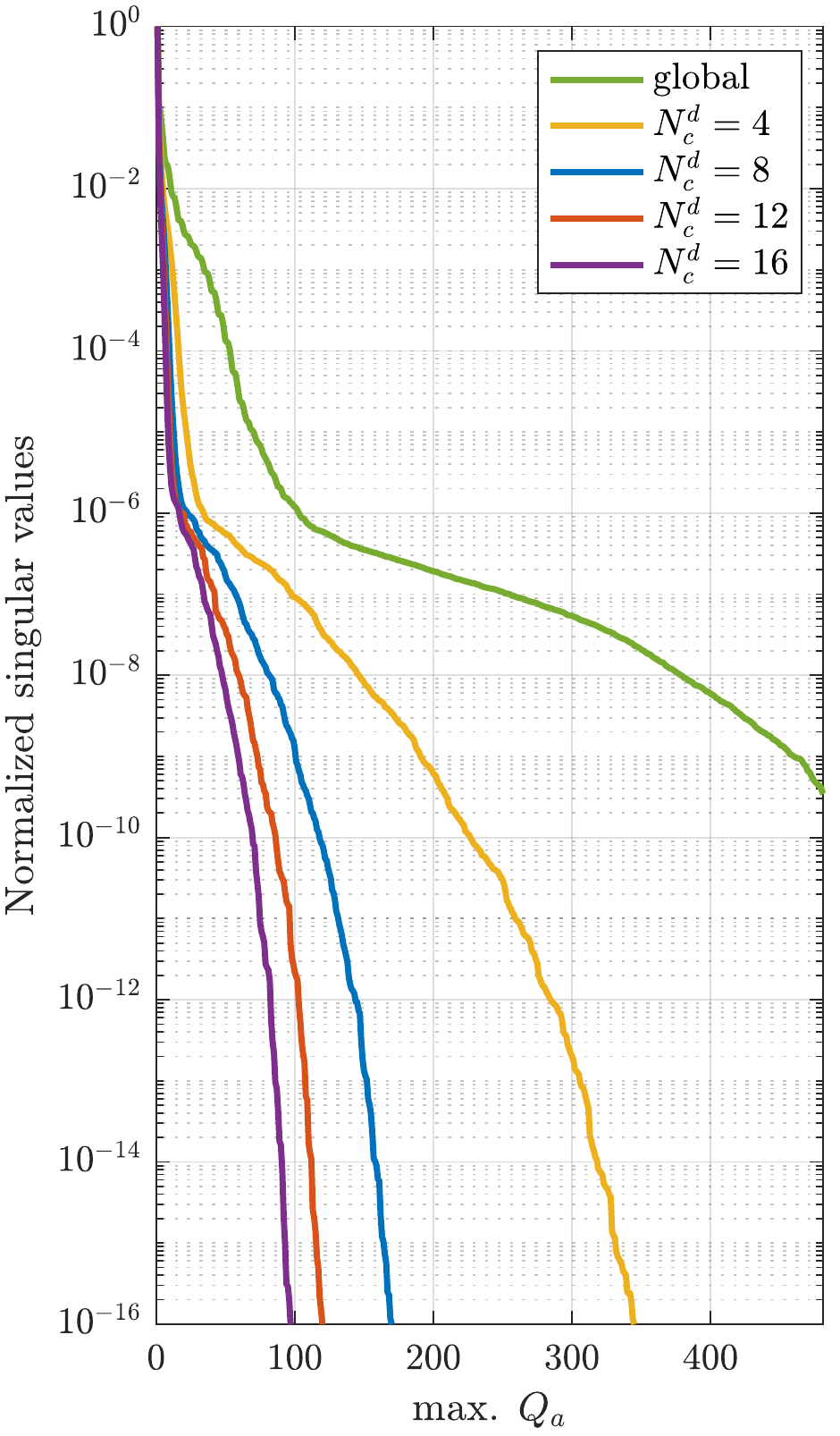}
         \caption{Stiffness matrix}
         \label{fig:PODDEIM_2}
     \end{subfigure}
        \caption{Example 6.1.1: Comparison of singular values decay between global and local DEIM approximations for right-hand side vector (a) and matrix (b)  using different number of clusters.}\label{fig:DEIM_1}
\end{figure}

Now we consider the strategy presented in Section \ref{sec:local}. To perform the localized DEIM approximation, the snapshots matrices ${\bf{S}}_{a}$ and ${\bf{S}}_{f}$ are subdivided considering the vector of parameters as a cluster indicator. We recall that the dimensionality of the training set $\mathcal{P}_{train}^d$ should be chosen sufficiently high for each cluster. Figure \ref{fig:DEIM_1} shows the decay of the singular values of the POD for the DEIM approximations with respect to the maximum number of selected basis functions ($Q_a$, $Q_f$) over all clusters for a given number of clusters $N_{c}^d$. It is evident that the number of affine terms is significantly reduced by using local subspaces.  In Table \ref{tab:DEIM_functions} we compare the selected number of affine terms for different number of clusters. Here, we use the same number of clusters for the matrix and right-hand side, although in principle this could be different. As discussed in Section \ref{sec:local}, the number of terms may differ between clusters. Thus, we only depict the minimum and maximum number of terms over all clusters.

\begin{table}[h!]
 	\centering
 	\caption{Example 6.1.1: DEIM approximations for different number of clusters. Comparison in terms of minimum and maximum number of basis functions over all clusters for the matrix ($Q_a$) and right-hand side vector ($Q_f$) based on POD tolerance $\epsilon_{POD}^d= 10^{-7}$.}\label{tab:DEIM_functions}
 	\begin{tabular}{ccccc} \hline 
 		$N_c^d$ & min. $Q_a$ & max. $Q_a$ & min. $Q_f$ & max. $Q_f$ \\\hline
 		1  & 349 & 349 & 124 & 124 \\ 
 		4  & 118 & 122 & 38 & 39 \\ 
 		8  & 62 & 66 & 21 & 23 \\ 
 		12  & 42 & 47 & 15 & 18 \\ 
 		16  & 31 & 38 & 11 & 15 \\\hline
 	\end{tabular}
 \end{table}
 
\begin{figure}[!h]
%     \centering
     \begin{subfigure}[b]{0.49\textwidth}
         \centering
         \includegraphics[width=\textwidth]{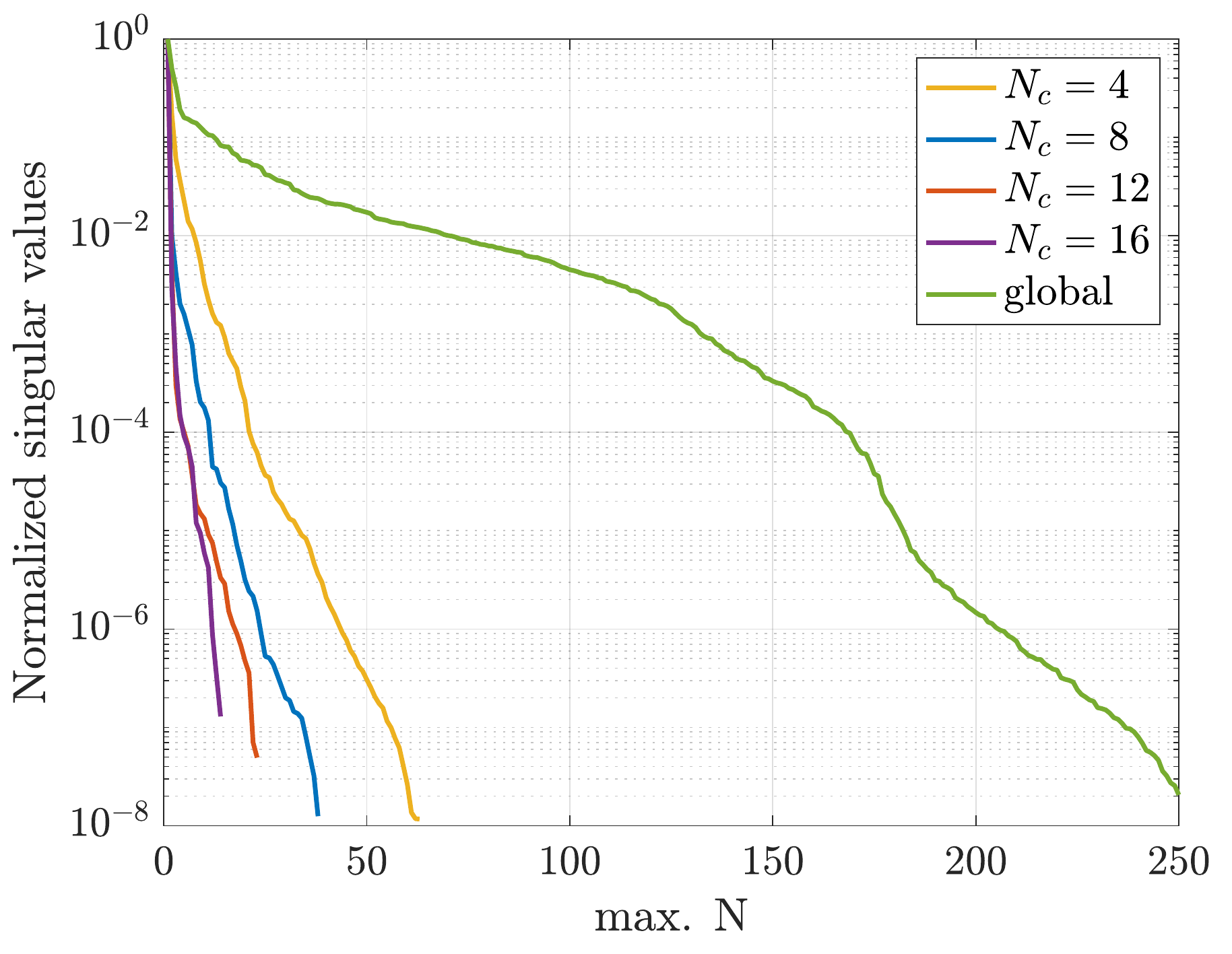}
         \caption{$N_{c}^d=16$}
         \label{fig:POD_1}
     \end{subfigure}
%     \hfill
     \begin{subfigure}[b]{0.46\textwidth}
         \centering
         \includegraphics[width=\textwidth]{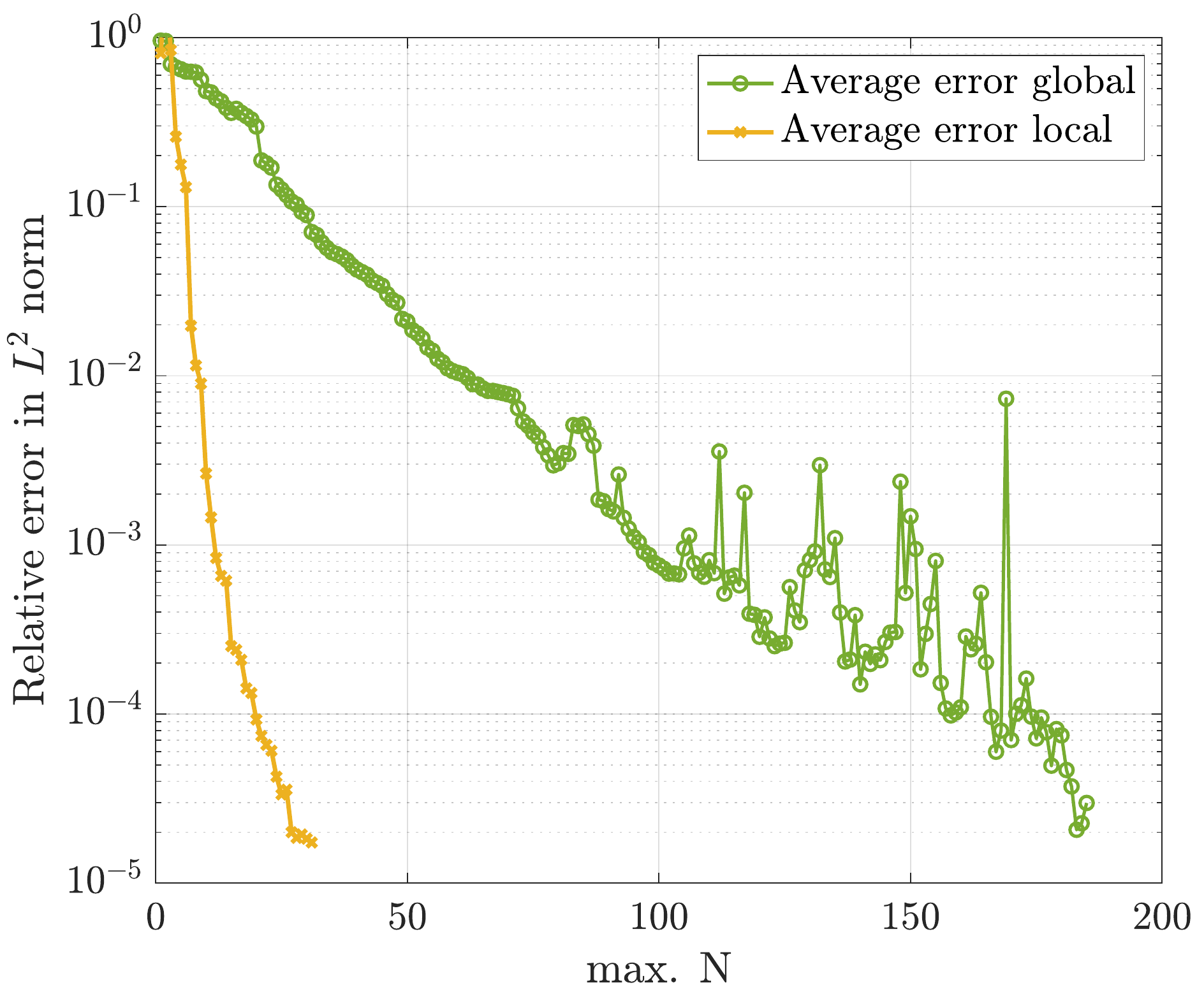}
         \caption{$N_{c}^d=16$, $N_c=4$}
         \label{fig:error_1}
     \end{subfigure}
        \caption{Example 6.1.1: Singular values decay for different numbers of clusters (a) and relative error vs. maximum number of reduced basis functions ($N$) over all clusters (b).}
      \label{fig:POD_error1}
\end{figure}

\begin{figure}[!h]
%     \centering
     \begin{subfigure}[b]{0.49\textwidth}
         \centering
         \includegraphics[width=\textwidth]{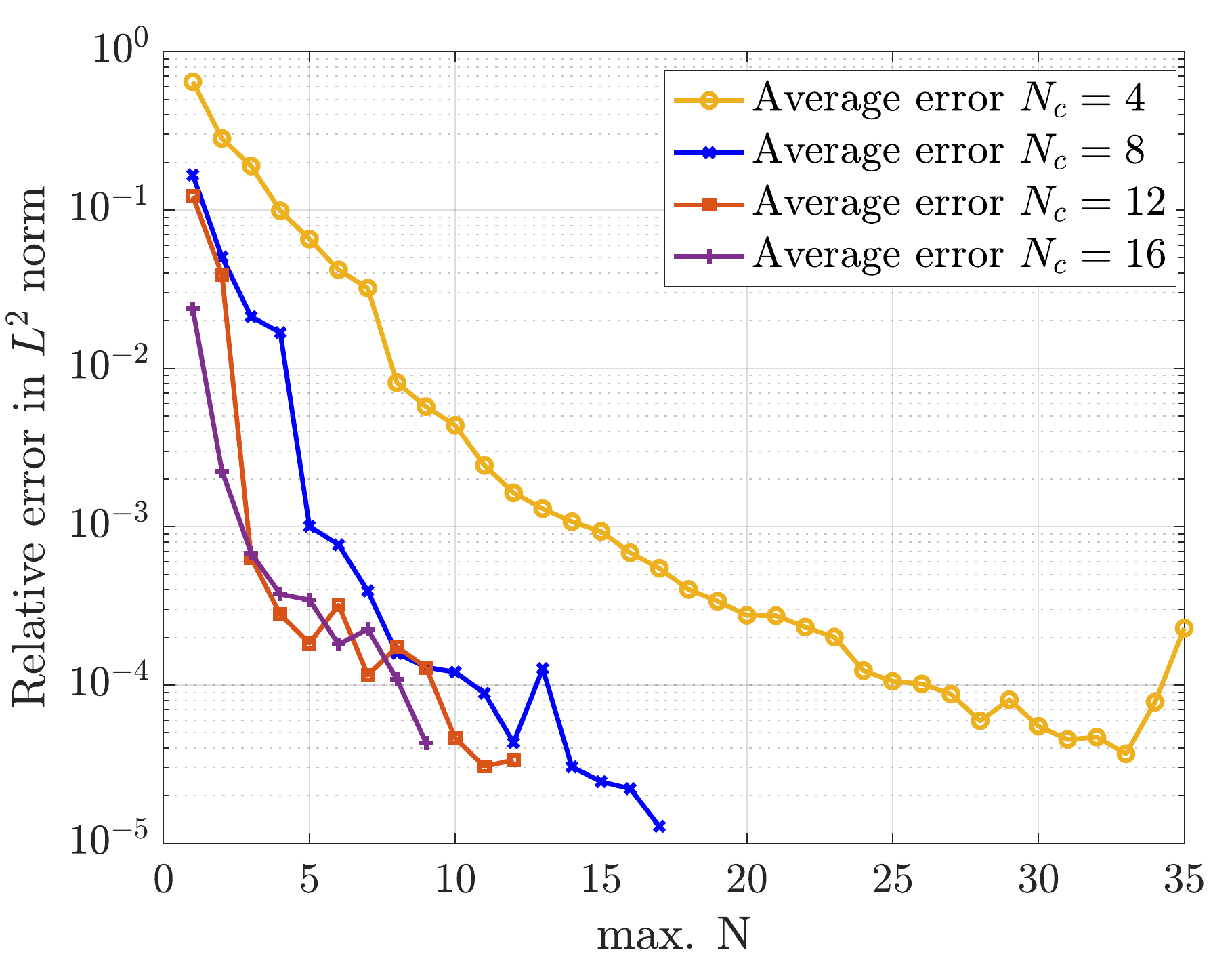}
         \caption{$N_{c}^d=8$}
         \label{fig:error_k1}
     \end{subfigure}
%     \hfill
     \begin{subfigure}[b]{0.49\textwidth}
         \centering
         \includegraphics[width=\textwidth]{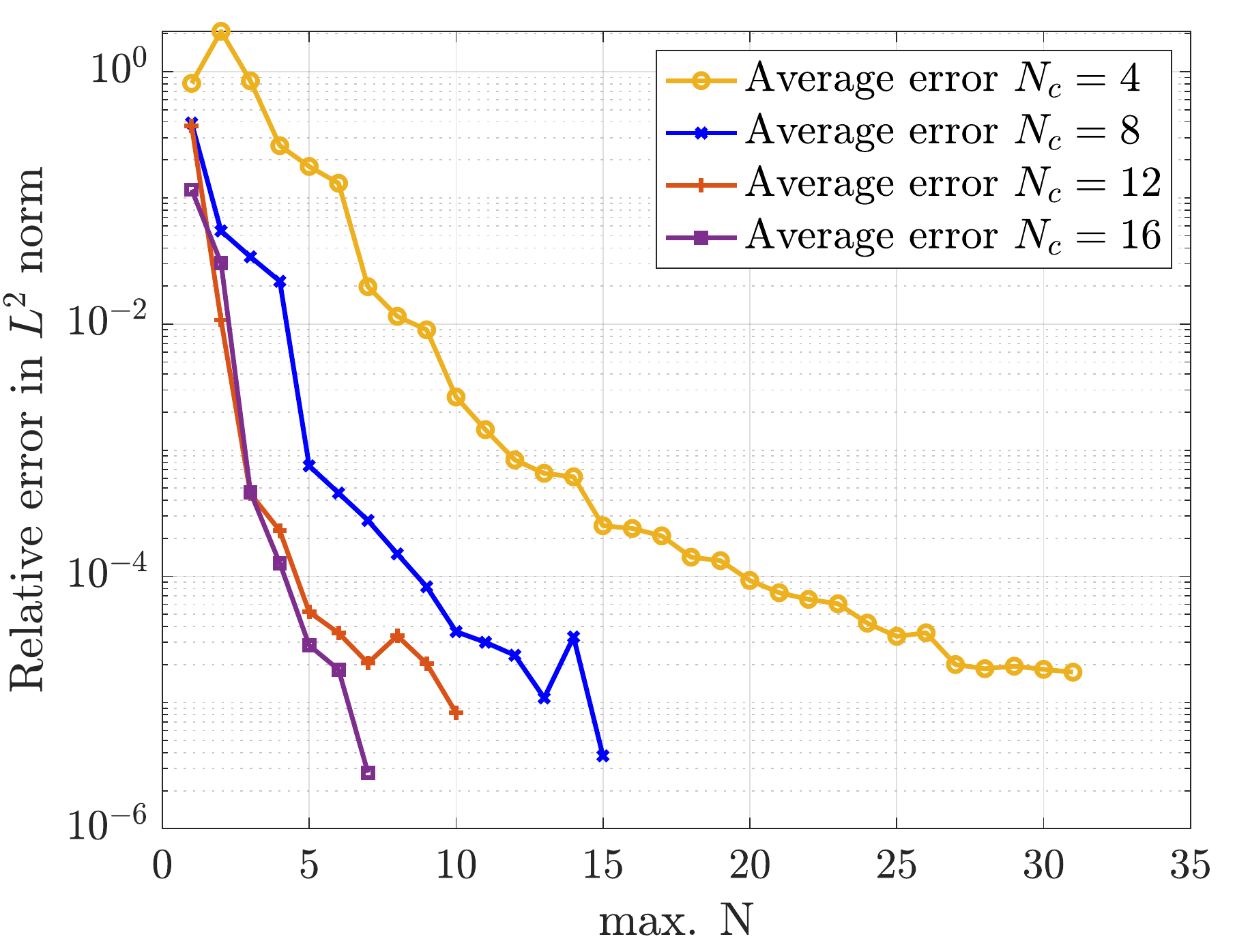}
         \caption{$N_{c}^d=16$}
         \label{fig:error_k2}
     \end{subfigure}
        \caption{Example 6.1.1: Decay of relative error vs. maximum number of reduced basis functions (N) for different numbers of clusters.}
      \label{fig:error_k}
\end{figure}

 Let us now assess the performance of the localization strategy in constructing a reduced basis with the POD. The solution snapshot matrix ${\bf{S}}_{u}$ is subdivided into clusters considering a training set of dimension $N_{s}=250$ as indicator, that is, the solutions are assigned to the same cluster as their respective parameters. We consider the number of clusters for the DEIM approximations fixed to $N_{c}^d=16$ and show the singular values decay in Figure \ref{fig:POD_1} for different number of clusters $N_c$ with respect to the maximum number of RB functions over all clusters. We observe that the clustering leads to a significant reduction of the RB functions $N$. Moreover, we perform an error analysis of the problem solution $u_h$ on a test sample of dimension $N_{t}=100$ constructed based on uniformly distributed random points in the parameter space. Figure \ref{fig:error_1} shows that a small number of clusters ($N_c=4$) is sufficient to achieve a ROM with accuracy of $10^{-5}$ and a local reduced basis with maximum dimension of $N = 35$ over all clusters, compared to the global ROM that requires $N=182$ functions. Moreover, increasing the number of clusters further improves the accuracy, while the dimension of the basis is reduced as shown in Figure \ref{fig:error_k}. 

\begin{figure}[!htb]
	\centering
	\includegraphics[width=1.0\textwidth]{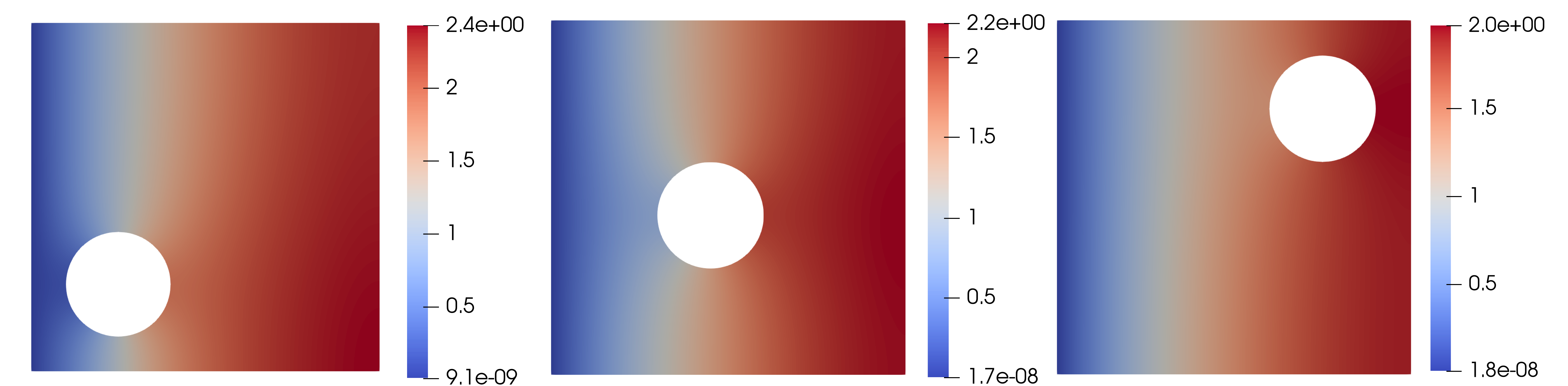} \\
	\caption{Example 6.1.1: ROM solutions for $\mu=[0.5,0.9,1.5]$  and $N_{c}^d=N_c=16$.}
	\label{fig:ROM_sol}
\end{figure}

Table \ref{tab:POD_functions} compares the efficiency of the ROMs. We observe that increasing the number of clusters $N_c$  reduces significantly the number of RB functions while it does not impair dramatically the online cost. Furthermore, the solution obtained with the local ROMs is illustrated in Figure \ref{fig:ROM_sol} for the same values of the parameter $\mu$ as in Figure \ref{fig:FOM_sol}.

 \begin{table}[!htb]
 	\centering
 	\caption{Example 6.1.1: Comparison of ROMs in terms of number of reduced basis functions (N) and computational cost for different numbers of clusters $N_c$ and $N_{c}^d=8$ with POD tolerance $\epsilon_{POD}= 10^{-5}$.}\label{tab:POD_functions}
 	\begin{tabular}{cccc} \hline 
 		$N_c$ & min. $N$ & max. $N$ & online CPU time [ms] \\\hline
 		1 & 182 & 182 & 879.0 \\
 		4  & 33 & 35 & 94.0 \\ 
 		8  & 10 & 17 & 85.7 \\ 
 		12  & 6 & 12 & 86.6 \\ 
 		16  & 4 & 9 & 93.9 \\\hline
 	\end{tabular}
 \end{table}

%\subsubsection{Two-parameter case study}
\subsubsection{Square with circular hole: 2D geometrical parameterization}
To illustrate the applicability of the methodology for multiple parameters, we consider a 2D geometrical parameterization for our problem. The parameter vector ${\bm{\mu}} = [\mu_1, \mu_2]$ represents the position of the center ($\mu_1$) and the radius of the hole ($\mu_2$) in the range $(\mu_1, \mu_2) \in [0.5,  1.5] \times [0.25,  0.35]$, while the circular hole is centered at $(\mu_1,\mu_1)$. This causes the number of active degrees of freedom to vary significantly between snapshots, which introduces a strong complexity to the solution's manifold. Note that the 2D parameter vector acts as indicator for clustering the snapshots. 

\begin{figure}[!h]
	\begin{subfigure}[b]{0.49\textwidth}
		\centering
		\includegraphics[width=0.9\textwidth]{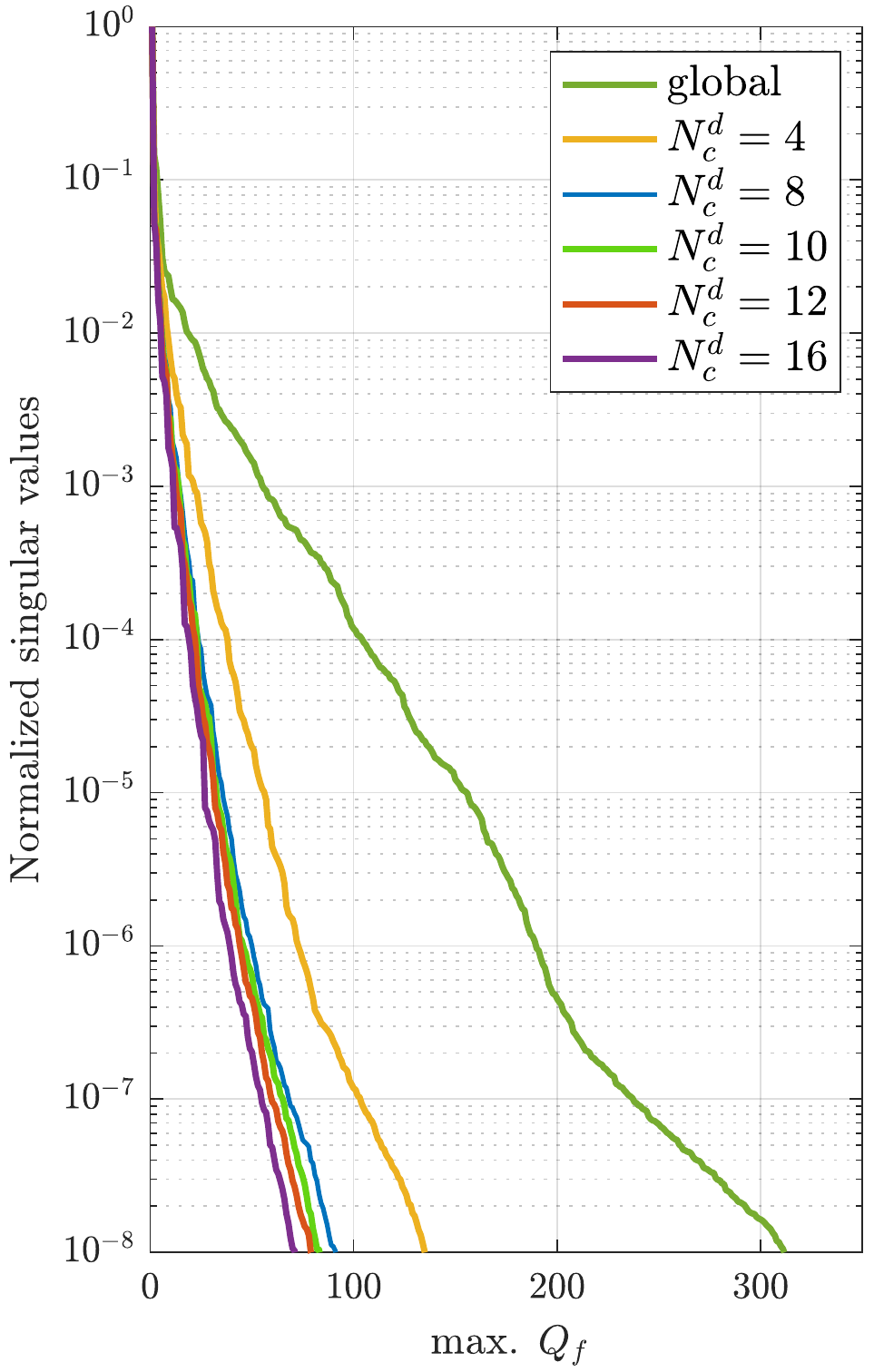}
		\caption{Right-hand side}
		\label{fig:RHS_2}
	\end{subfigure}
	\begin{subfigure}[b]{0.49\textwidth}
		\centering
		\includegraphics[width=0.9\textwidth]{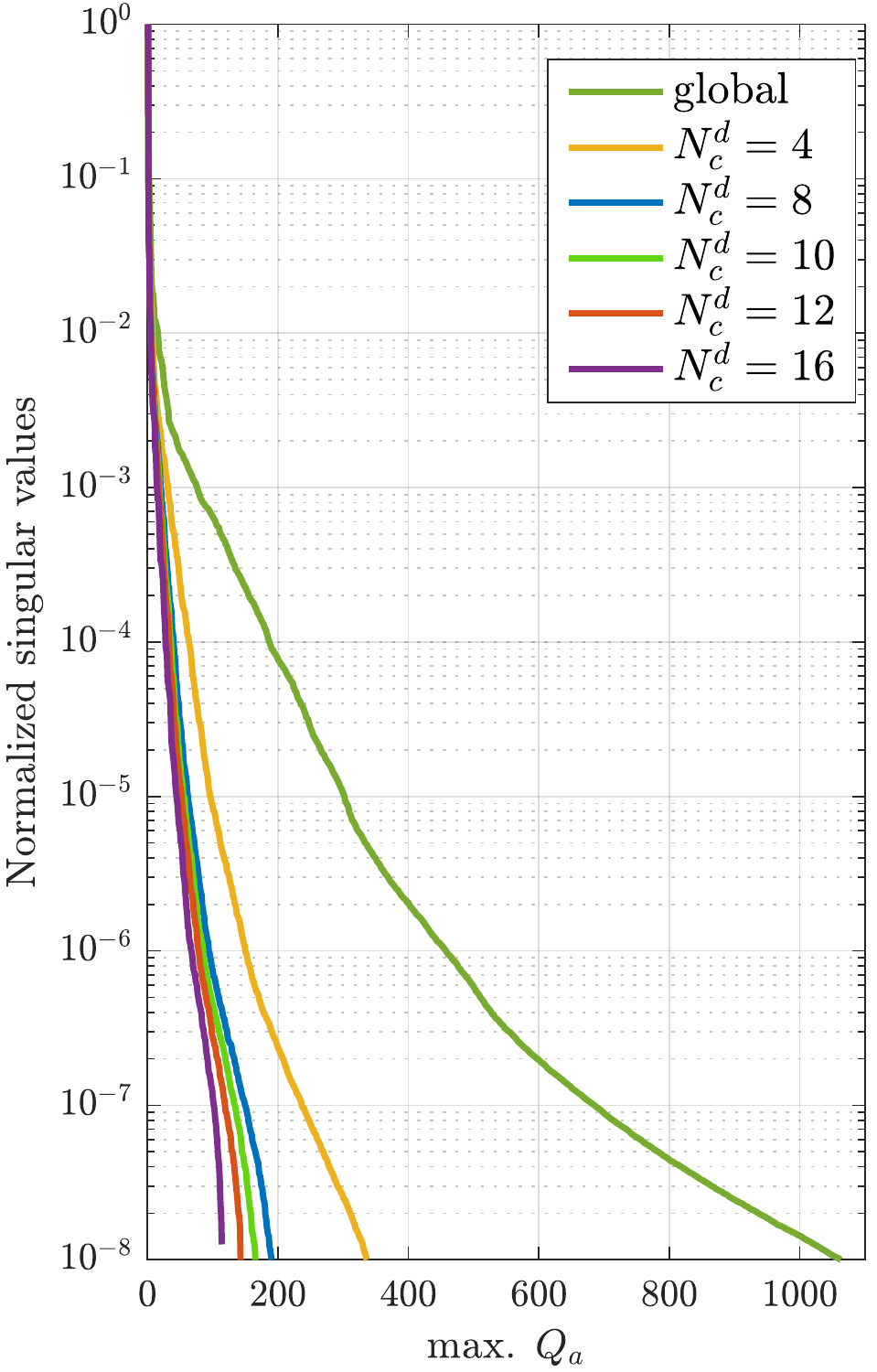}
		\caption{Stiffness matrix}
		\label{fig:matrix_2}
	\end{subfigure}
	\caption{Example 6.1.2: Comparison of singular values decay between global and local DEIM approximations for right-hand side vector (a) and matrix (b) using different numbers of clusters.}\label{fig:DEIM_2}
\end{figure}

\begin{figure}[!h]
     \begin{subfigure}[b]{0.49\textwidth}
         \centering
         \includegraphics[width=\textwidth]{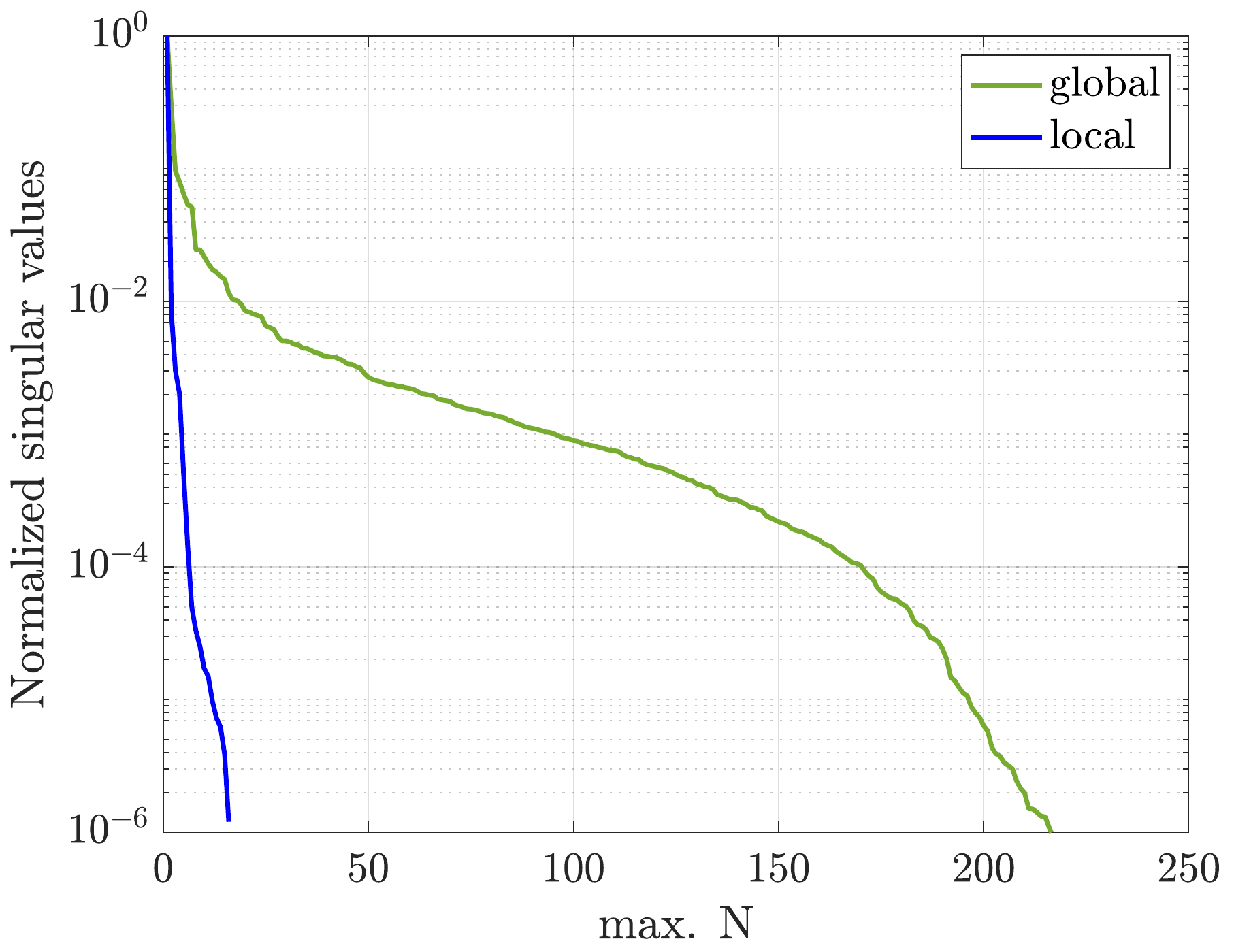}
         \caption{Singular values decay}
         \label{fig:POD_2}
     \end{subfigure}
     \begin{subfigure}[b]{0.49\textwidth}
         \centering
         \includegraphics[width=\textwidth]{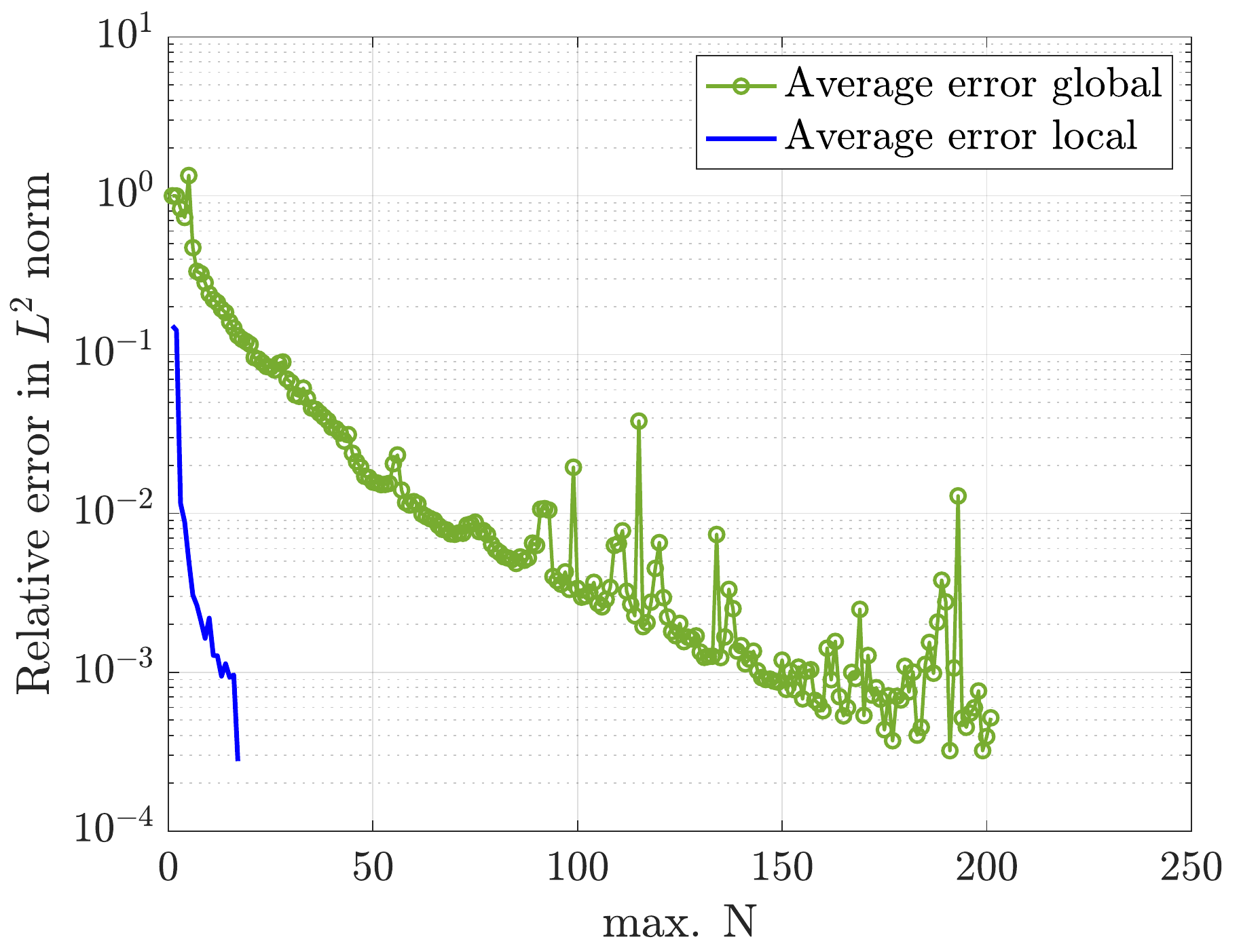}
         \caption{Error analysis}
         \label{fig:error_2}
     \end{subfigure}
        \caption{Example 6.1.2: Singular values decay for different number of clusters (a) and relative error vs. maximum number of reduced basis functions ($N$) over all clusters (b) for $N_{c}^d=N_c=16$.}
      \label{fig:POD_error2}
\end{figure}

\begin{figure}[!htb]
	\centering
	\includegraphics[width=1.0\textwidth]{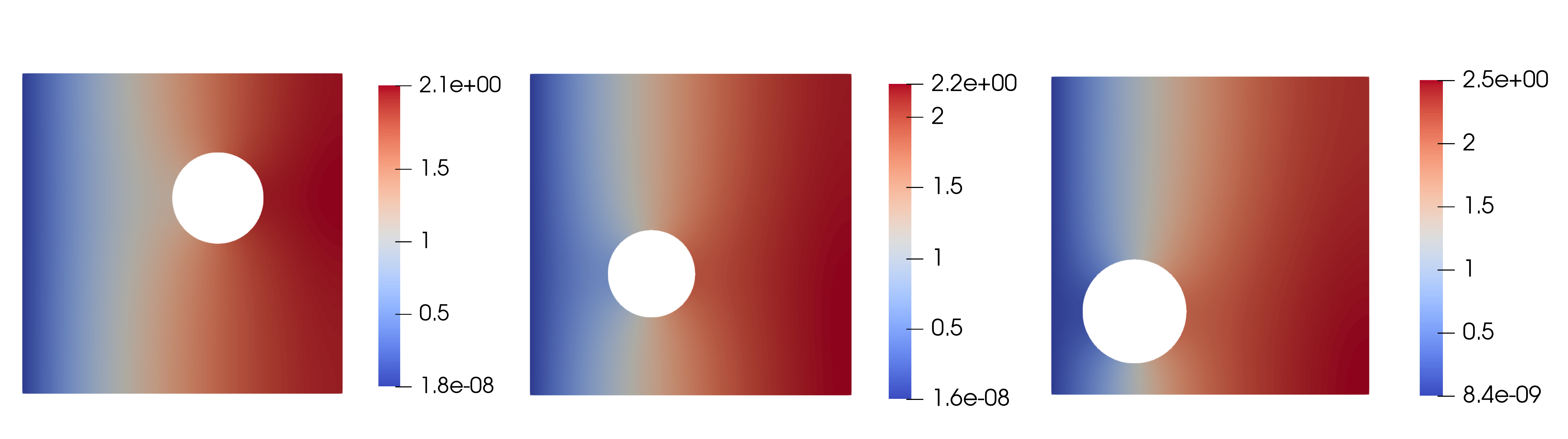} \\
	\caption{Example 6.1.2: Solution snapshots corresponding to selected cluster centroids for $\mu_1 = [1.2214, 0.7562, 0.5237]$ and $\mu_2 = [0.2851, 0.2715, 0.3268]$.}\label{fig:centroid}
\end{figure}

In fact, the number of affine terms selected by the global DEIM is much higher than the 1D parameterization as depicted in Figure \ref{fig:DEIM_2}. The use of local subspaces significantly reduces the maximum number of affine terms over all clusters. Figure \ref{fig:POD_error2} shows that the number of RB functions is effectively reduced without compromising the accuracy. The error analysis is performed on a test sample of dimension $N_{t}=100$ for a fixed number of clusters $N_{c}^d=N_c=16$. Moreover, the solution snapshots for three exemplary cluster centroids computed by \emph{k-means} are depicted in Figure \ref{fig:centroid}. In Table \ref{tab:FOM_comparison}, we compare the performance between the local and global ROMs: a 17.6x speedup is achieved with the local ROM, versus a speedup of 9.1x for the global ROM, both with respect to the FOM.

\begin{table}[!htb]
	\centering
	\caption{Example 6.1.2: Comparison of global and local ROM with $N_{c}^d=N_c=16$ in terms of maximum number of basis functions over all clusters and computational cost. POD tolerance set to $\epsilon_{POD}^d= 10^{-7}$ for the DEIM approximation and $\epsilon_{POD}= 10^{-5}$ for the reduced basis.}\label{tab:FOM_comparison}
	\begin{tabular}{cccccc} \hline 
		& max. $Q_a$ & max. $Q_f$ & max. $N$ & Online CPU time [ms] & speedup \\\hline
		local  & 107 & 59 & 17 & 122 & 17.6x \\ 
		global & 1024 & 282 & 201 & 251 & 9.1x \\\hline
	\end{tabular}
\end{table}

\begin{figure}[!h]
	\begin{subfigure}[b]{0.49\textwidth}
		\centering
		\includegraphics[width=\textwidth]{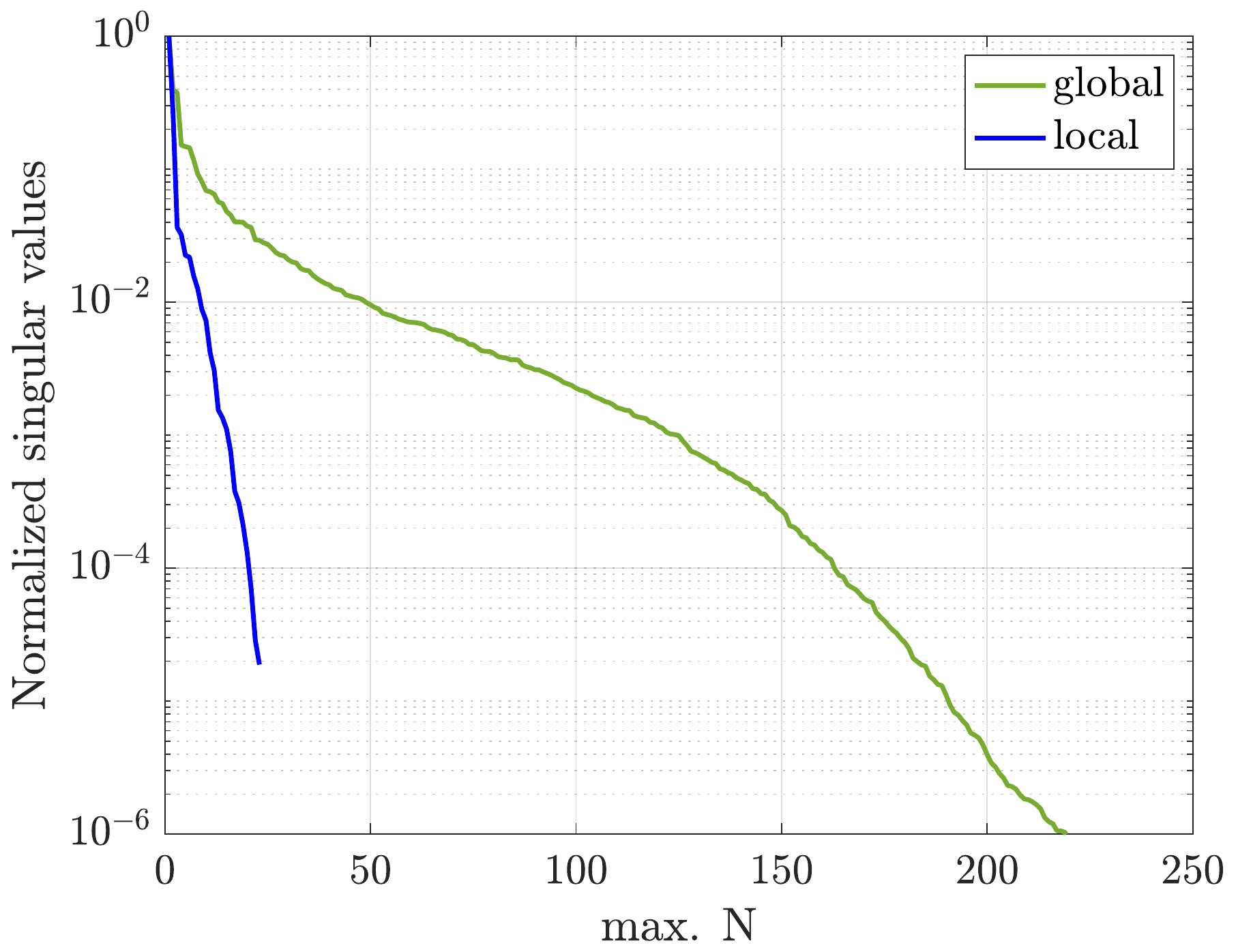}
		\caption{Singular values decay}
		\label{fig:POD_2_2}
	\end{subfigure}
	\begin{subfigure}[b]{0.49\textwidth}
		\centering
		\includegraphics[width=\textwidth]{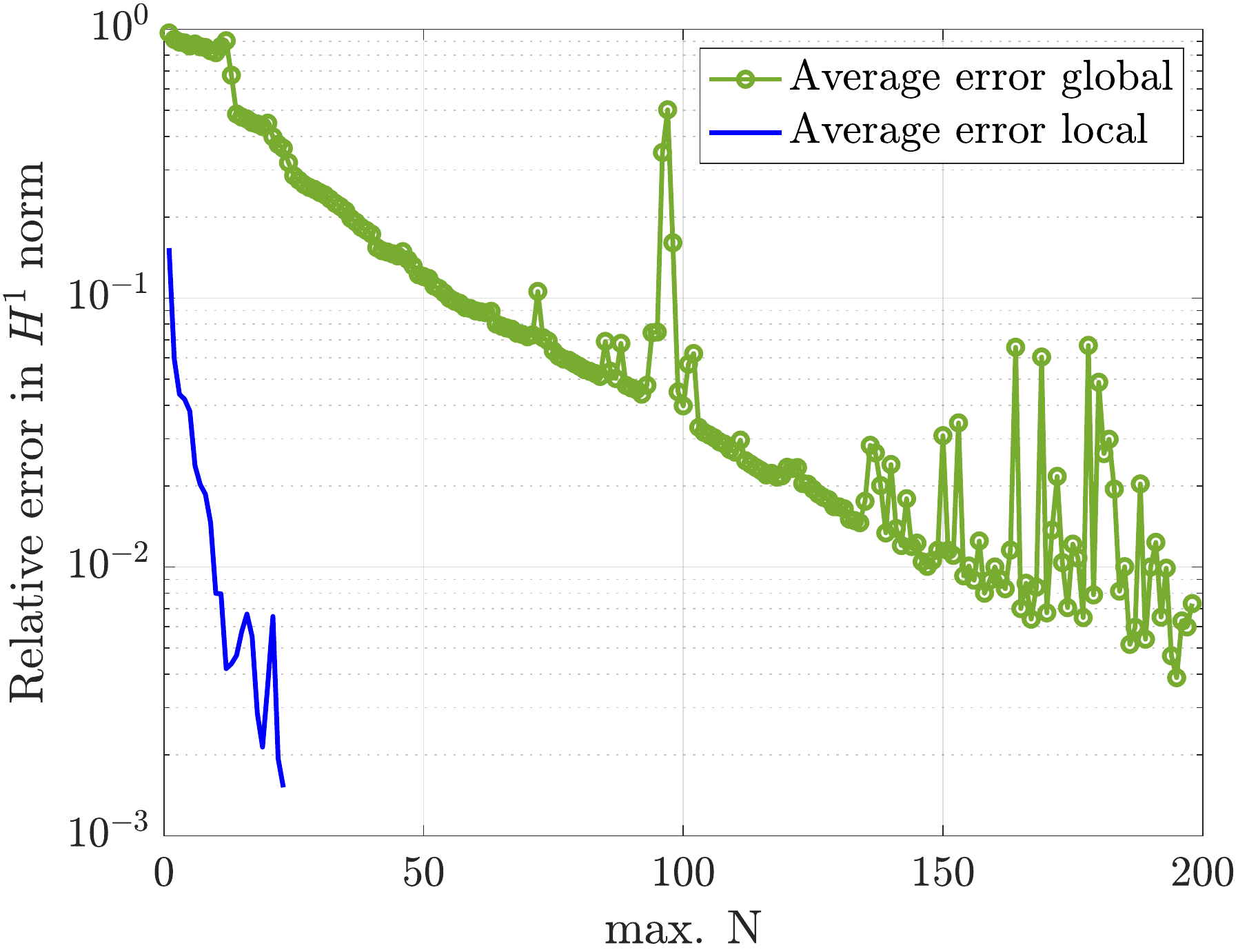}
		\caption{Error analysis}
		\label{fig:error_2_2}
	\end{subfigure}
	\caption{Example 6.1.2: Singular values decay for different number of clusters (a) and relative error in $H^1$ norm vs. maximum number of reduced basis functions ($N$) over all clusters (b) for $N_{c}^d=N_c=16$.}
	\label{fig:POD_error2_2}
\end{figure}

 Finally, in Figure \ref{fig:POD_error2_2} we show the error analysis of the problem solution $u_h$ with respect to the $H^1$ norm using again a test sample of dimension $N_{t}=100$ and a fixed number of clusters $N_{c}^d=N_c=16$. We remark that the POD basis is constructed such that it minimizes the squared projection error with respect to the algebraic counterpart of the $H^1$ norm \cite[Proposition~6.2]{QMN_RBspringer}. Similarly to the previous test cases, the singular values decay rapidly and the reduction is more effective in the case of localization.

\subsection{Linear elasticity}
Let us now briefly recall the equations of linear elasticity on a parameterized domain. We consider an isotropic parameterized solid $\Omega(\bm{\mu}) \subset \mathbb{R}^d$ with elastic deformations described in terms of a stress tensor $\boldsymbol{\sigma}$, a small strain tensor $\boldsymbol{\varepsilon}$, the body force vector ${\bm{f}}$ and the unknown displacement field $\boldsymbol{u}$. The Dirichlet and Neumann part of the boundary of the domain $\partial{\Omega}(\bm{\mu})$ are denoted by $\Gamma_{D}(\bm{\mu})$ and $\Gamma_{N}(\bm{\mu})$, respectively, while $\boldsymbol{n}$ is the outward unit normal to the boundary. Similarly to the previous case of the Poisson problem, homogeneous Dirichlet and Neumann boundary conditions are assumed, with $\Gamma_D({\bm{\mu}})\subset \partial{\Omega({\bm{\mu}})}\cap\partial{\Omega_0}$. The continuous formulation of the problem in strong form reads: for any $\bm{\mu} \in \mathcal{P}$, find $\boldsymbol{u} \in [H_{0,{\Gamma_D}}^1(\Omega(\bm{\mu}))]^d$ such that
\begin{equation}\label{eq24}
\begin{cases}
- \text{div} (\boldsymbol{\sigma}(\boldsymbol{u})) &= {\bm{f}} \qquad \qquad \qquad \qquad \qquad \ \text{in} \ \Omega({\bm{\mu}}) \\
\qquad \ \ \boldsymbol{\sigma}(\boldsymbol{u}) &= 2\tilde{\mu}\boldsymbol{\varepsilon}(\boldsymbol{u}) +  \tilde{\lambda}(\text{div} (\boldsymbol{u}))\boldsymbol{I} \qquad \text{in} \ \Omega({\bm{\mu}}) \\
\qquad \ \ \boldsymbol{\varepsilon}(\boldsymbol{u}) &= \displaystyle \frac{1}{2}(\nabla{\boldsymbol{u}} + (\nabla{\boldsymbol{u}})^T) \ \qquad \ \ \ \ \ \text{in} \ \Omega({\bm{\mu}}) \\
\qquad \ \  \boldsymbol{u} &= 0 \qquad \qquad \qquad \qquad \qquad \ \ \text{on} \ \Gamma_D(\bm{\mu}) \\
 \qquad \ \ \boldsymbol{\sigma}(\boldsymbol{u})\cdot{\boldsymbol{n}} &= 0 \qquad \qquad \qquad \qquad \qquad \ \  \text{on} \ \Gamma_{N}(\bm{\mu}). 
\end{cases}
\end{equation}
Here, the Lam\'e coefficients $\tilde{\mu}$ and $\tilde{\lambda}$ can be expressed with respect to the Young modulus $E$ and Poisson coefficient $\nu$ as 
\begin{equation}\label{eq25}
\tilde{\mu}= \frac{E}{2(1+\nu)}, \qquad \tilde{\lambda}= \frac{E\nu}{(1+\nu)(1-2\nu)}.
\end{equation}
The discrete weak formulation of the parameterized problem in Equation \eqref{eq24} can be expressed in a similar manner to the Poisson problem as: find $\boldsymbol{u}_h \in V_h$ such that 
\begin{equation}\label{eq26}
a(\boldsymbol{u}_h,\boldsymbol{v}_h;\bm{\mu}) = f(\boldsymbol{v}_h;\bm{\mu}), \qquad \forall \boldsymbol{v}_h \in V_h,
\end{equation}
where $V_h \subset [H_{0,\Gamma_D}^1(\Omega(\bm{\mu}))]^d$ is a vector subspace spanned by a B-spline basis. Then the parameterized bilinear form $a(\cdot,\cdot;\bm{\mu})$ is given as:
\begin{equation}\label{eq27a}
%\begin{aligned}
a(\boldsymbol{u}_h,\boldsymbol{v}_h;\bm{\mu}) = \int_{{\Omega}(\bm{\mu})} 2\tilde{\mu} \boldsymbol{\varepsilon}(\boldsymbol{u}_h) : \boldsymbol{\varepsilon}(\boldsymbol{v}_h) \, \textrm{d}\Omega + \int_{{\Omega}(\bm{\mu})} \tilde{\lambda} \text{div}(\boldsymbol{u}_h) \text{div}(\boldsymbol{v}_h) \, \textrm{d}\Omega, 
%\end{aligned}
\end{equation}
and the linear functional $f(\cdot;\bm{\mu})$ reads:
\begin{equation}\label{eq27b}
%\begin{aligned}
f(\boldsymbol{v}_h;\bm{\mu}) = \int_{{\Omega}(\bm{\mu})} {\bm{f}}\cdot \boldsymbol{v}_h \, \textrm{d}\Omega.
%\end{aligned}
\end{equation}

\subsubsection{Multi-perforated quarter cylinder }\label{sec:ring}
\begin{figure}[!htb]
	\centering
	\includegraphics[width=1.0\textwidth]{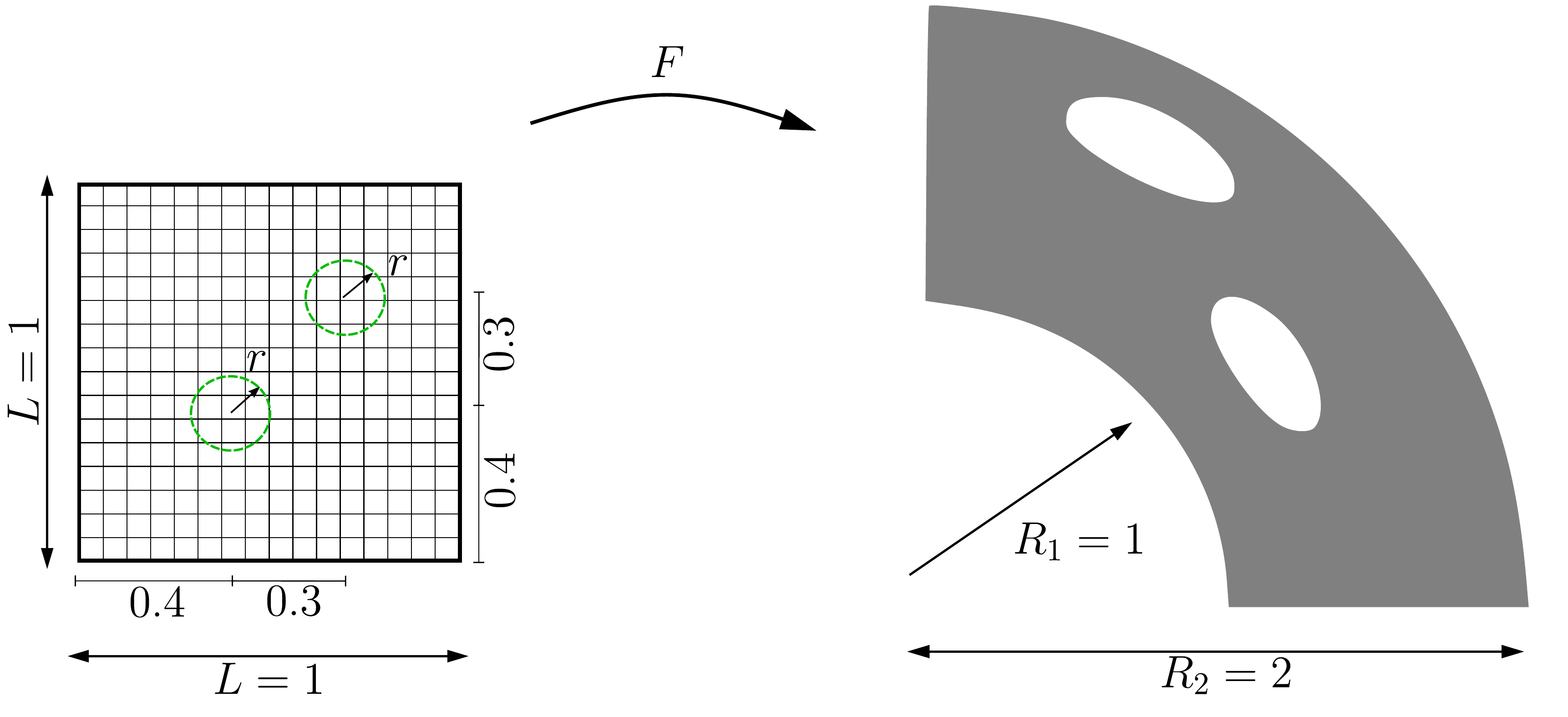} \\
	\caption{Example 6.2.1: Geometry of the multi-perforated quarter cylinder.}\label{fig:system_ring}
\end{figure}

In this example we assess the performance of the proposed procedure for geometries with multiple trimmed regions. For this purpose, we consider a two-dimensional geometry of one-quarter cylindrical ring with multiple holes. The model is defined in Figure \ref{fig:system_ring}. The trimmed domain $\Omega(\mu)$ is parameter-dependent, where $\mu = r \in [0.1, 0.2]$ is a parameter representing the radius of the circular holes in the pre-image domain, that is a Cartesian grid as discussed in Section \ref{sec:unfitted}, Remark \ref{remark1}. In this case, we make use of a spline mapping $F$ to obtain the actual geometry (see Remark \ref{remark2}).

Homogeneous Dirichlet boundary conditions are imposed on all four boundaries of the domain and the body force is set to ${\bm{f}}=[{f}_x,{f}_y]=[2xy, 2xy]$. The Young modulus and Poisson coefficient are given as $E=1.0$ and $\nu=0.3$, respectively. The geometry is discretized with quadratic $C^1$-continuous B-splines employing a mesh with 32 elements per direction of the Cartesian grid, resulting in $\mathcal{N}_{h,0}=2312$ degrees of freedom. The number of active degrees of freedom changes significantly for different values of the parameter and our problem is characterized by strong solution variations in different regions of the parameter space. 

\begin{figure}[!htb]
	\centering
	\includegraphics[width=0.6\textwidth]{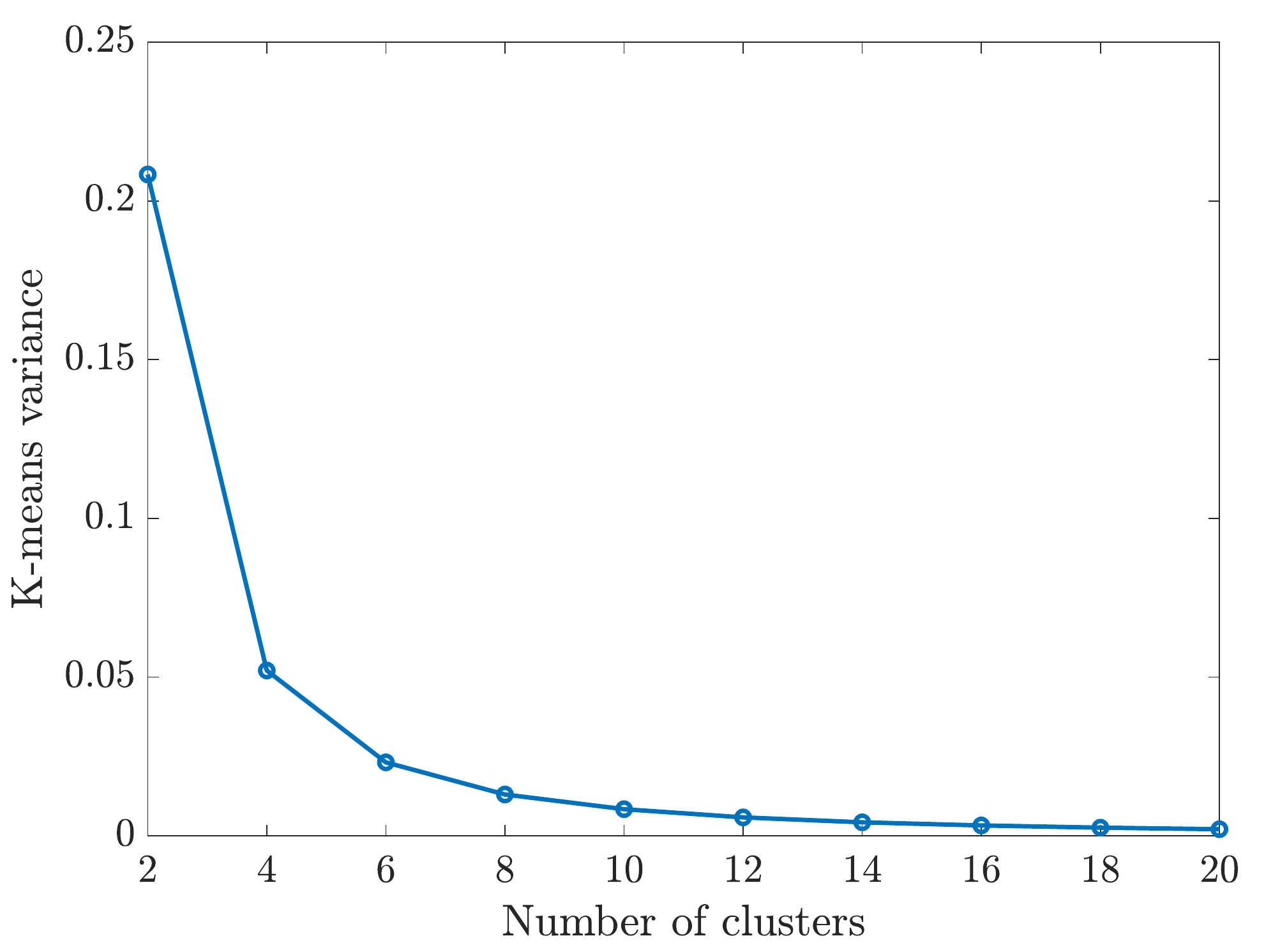} \\
	\caption{Example 6.2.1: K-means variance versus number of clusters $N_c^d$.}\label{fig:kmeans_ring}
\end{figure}

In order to construct the affine approximations and reduced bases, we consider a sufficiently rich training set of dimension $N_{s}^d=1000$ and $N_s=500$, accordingly. It should be noted that these dimensions refer to the global snapshot matrices, that is the number of snapshots in each cluster after partitioning should be sufficiently high to obtain accurate approximations. Figure \ref{fig:kmeans_ring} shows the \emph{k-means} variance \eqref{eq23} computed for different numbers of clusters $N_c^d$. It is observed that the variance does not decrease significantly after 10 clusters, therefore this number is chosen as the optimal one.

Table \ref{tab:FOM_comparison_ring} summarizes the results of the comparison between the local and global ROM. The number of basis functions is reduced significantly for both the affine approximations and the reduced basis. Regarding the performance of the ROMs, the results indicate that a 130x speedup is achieved with the local ROM, versus a speedup of 16.5x for the global ROM, both with respect to the FOM. Note that as expected, the overall efficiency of the ROM depends highly on the number of affine terms $Q_a,Q_f$. Compared to the results of the previous example in Table \ref{tab:FOM_comparison}, the number of affine terms is much lower here, which leads to a significantly higher speedup of the local ROM. 

\begin{table}[!htb]
 	\centering
 	\caption{Example 6.2.1: Comparison of global and local ROM with $N_{c}^d=N_c=10$ in terms of maximum number of basis functions over all clusters and computational cost. POD tolerance set to $\epsilon_{POD}^d= 10^{-7}$ for the DEIM approximation and $\epsilon_{POD}= 10^{-5}$ for the reduced basis.}\label{tab:FOM_comparison_ring}
 	\begin{tabular}{cccccc} \hline 
 		 & max. $Q_a$ & max. $Q_f$ & max. $N$ & Online CPU time [ms] & speedup \\\hline
 		local  & 36 & 20 & 16 & 36.7 & 130x \\ 
 		global & 260 & 117 & 198 & 291 & 16.5x \\\hline
 	\end{tabular}
 \end{table}

\begin{figure}[!h]
     \begin{subfigure}[b]{0.49\textwidth}
         \centering
         \includegraphics[width=\textwidth]{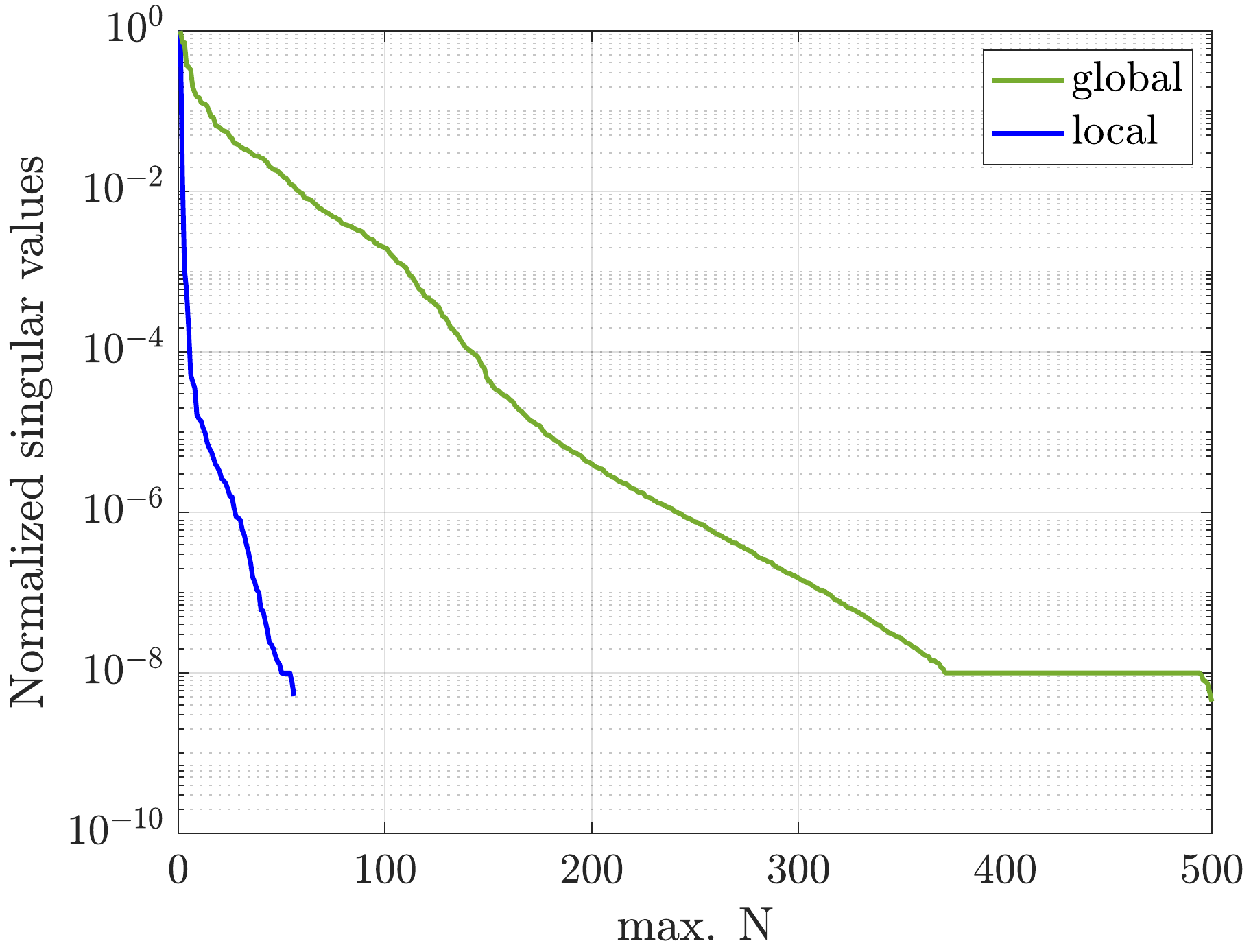}
         \caption{Singular values decay}
         \label{fig:POD_error_ring1}
     \end{subfigure}
     \begin{subfigure}[b]{0.49\textwidth}
         \centering
         \includegraphics[width=\textwidth]{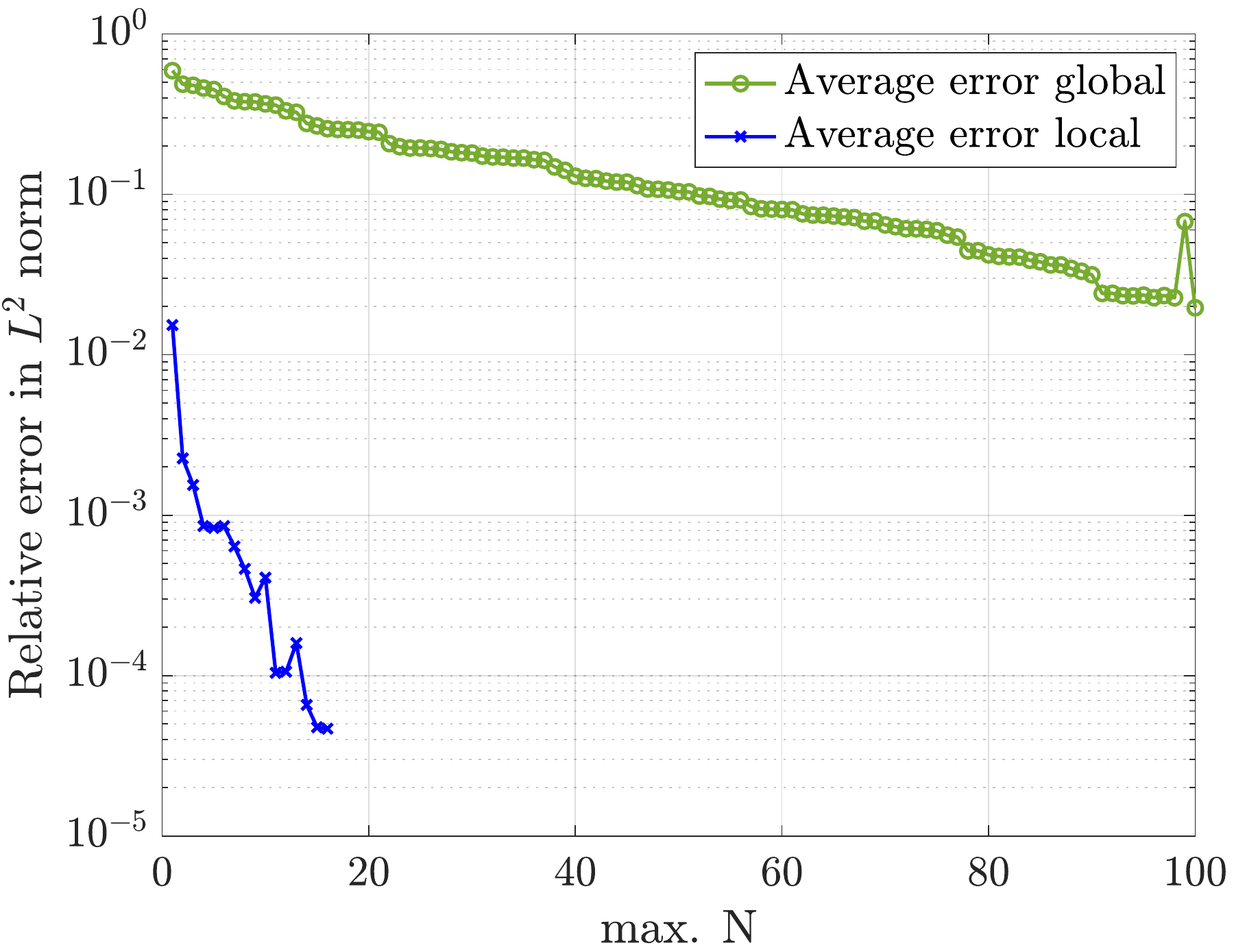}
         \caption{Error analysis}
         \label{fig:POD_error_ring2}
     \end{subfigure}
        \caption{Example 6.2.1: Decay of singular values (a) and relative error vs. maximum number of reduced basis functions ($N$) over all clusters (b) for $N_{c}^d=N_c=10$.}
      \label{fig:POD_error_ring}
\end{figure}

\begin{figure}[!h]
\begin{subfigure}[b]{0.49\textwidth}
	\centering
	\includegraphics[width=\textwidth]{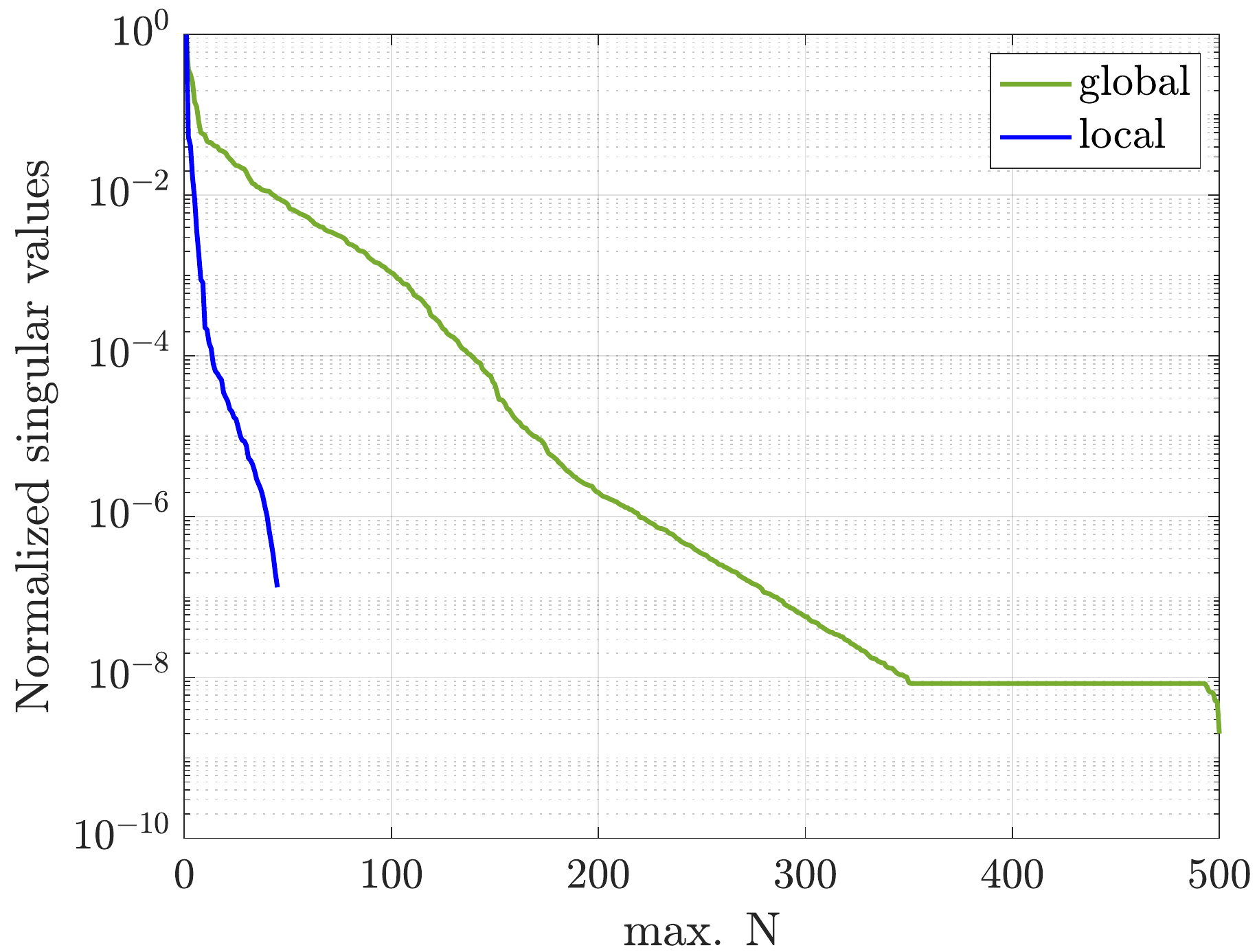}
	\caption{Singular values decay}
	\label{fig:POD_error_ring1_2}
\end{subfigure}
\begin{subfigure}[b]{0.49\textwidth}
	\centering
	\includegraphics[width=\textwidth]{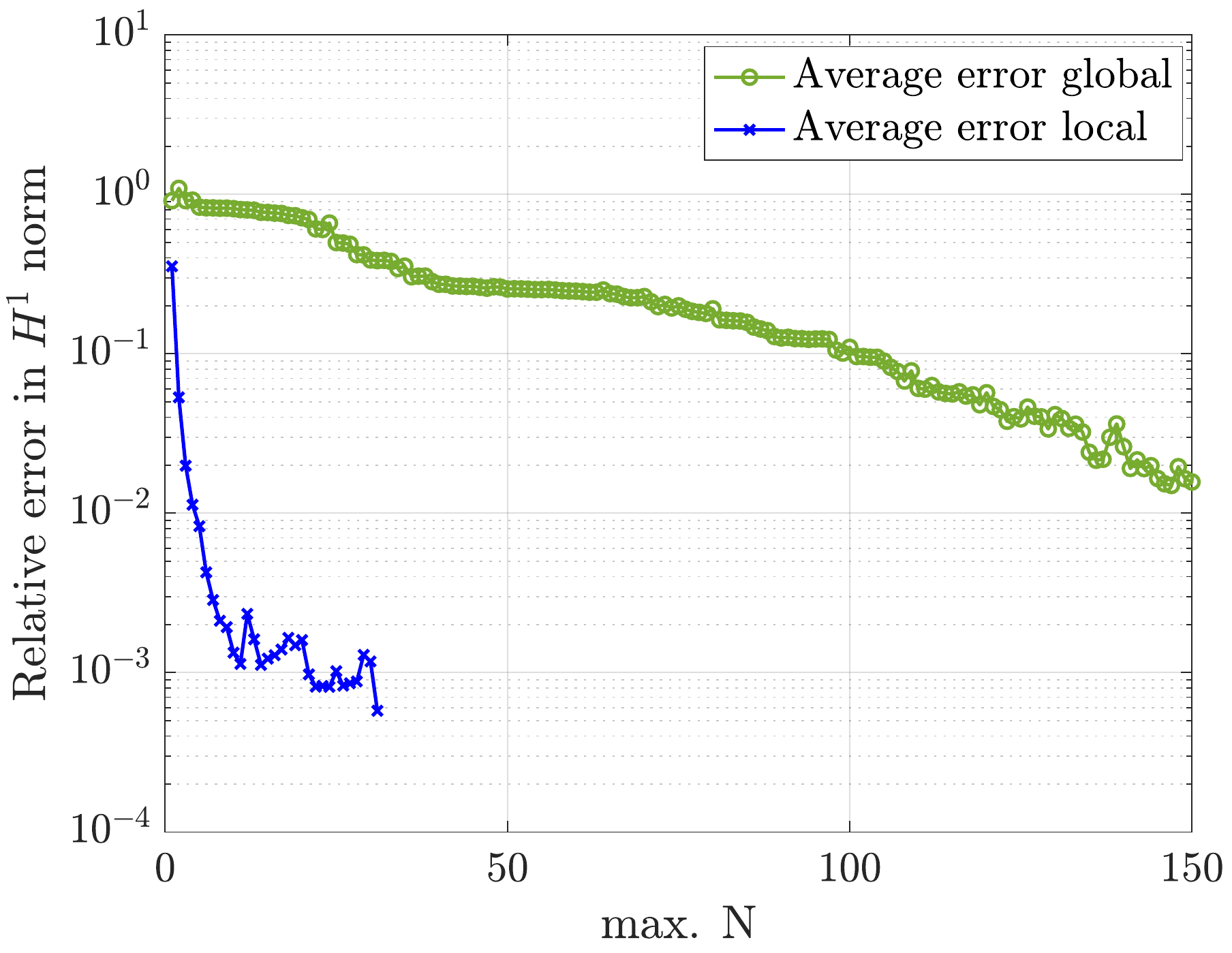}
	\caption{Error analysis}
	\label{fig:POD_error_ring2_2}
\end{subfigure}
\caption{ Example 6.2.1: Decay of singular values (a) and relative error in $H^1$ norm vs. maximum number of reduced basis functions ($N$) over all clusters (b) for $N_{c}^d=N_c=10$.}
\label{fig:POD_error_ring_2}
\end{figure}

\begin{figure}[!htb]
	\centering
	\includegraphics[width=1.0\textwidth]{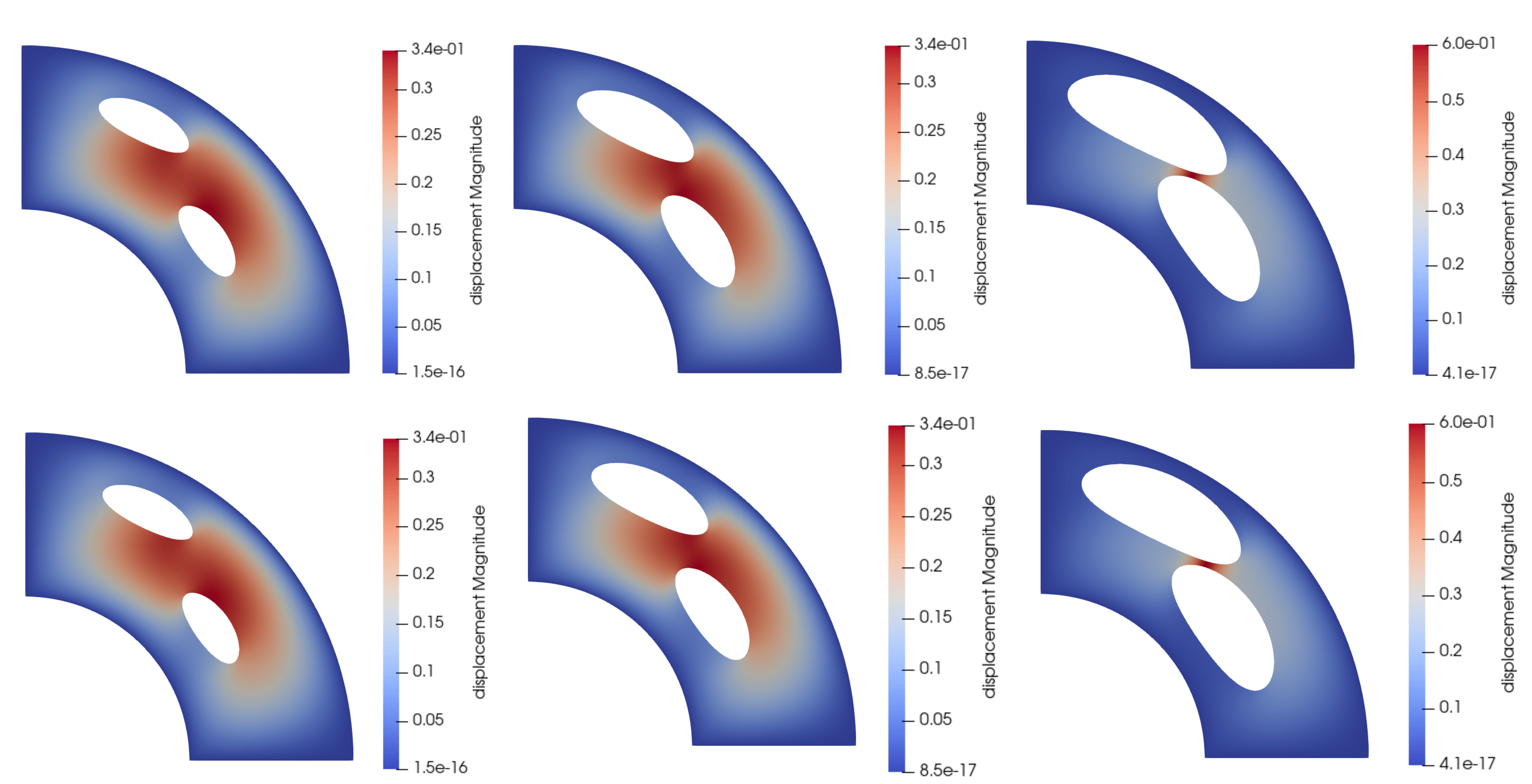} \\
	\caption{Example 6.2.1: Solution computed with the FOM (top) and local ROM (bottom) with $N_{c}^d=N_c=10$ clusters for three parameter values from the test sample $\mu = [0.1106, 0.1440, 0.1981]$.}\label{fig:solution_plots}
\end{figure}

Moreover, Figure \ref{fig:POD_error_ring} depicts the singular values and error decay for the ROM with global and local reduced basis based on a test sample of dimension $N_t=30$. Similarly to the previous test cases for the Poisson problem, we observe a rapid decay for the local ROM. Moreover, we obtain an accuracy of the order $10^{-5}$ for the local reduced basis with maximum dimension of $N=16$ over all clusters. In Figure \ref{fig:POD_error_ring_2}, the error analysis with respect to the $H^1$ norm results in accuracy of the order $10^{-4}$ for the local reduced basis with maximum dimension of $N=31$ over all clusters.  As in the previous example, the reduced basis is constructed  also here such that it minimizes the squared projection error with respect to the algebraic counterpart of the $H^1$ norm. Finally, the solutions of the local ROM are compared to the FOM for three different values of the test sample in Figure \ref{fig:solution_plots}. It can be observed that the solutions vary highly for different values of the parameter. The results indicate a good qualitative agreement with the FOM.

\subsubsection{Cube with spherical inclusion}
\begin{figure}[!h]
     \begin{subfigure}[b]{0.49\textwidth}
         \centering
         \includegraphics[width=0.9\textwidth]{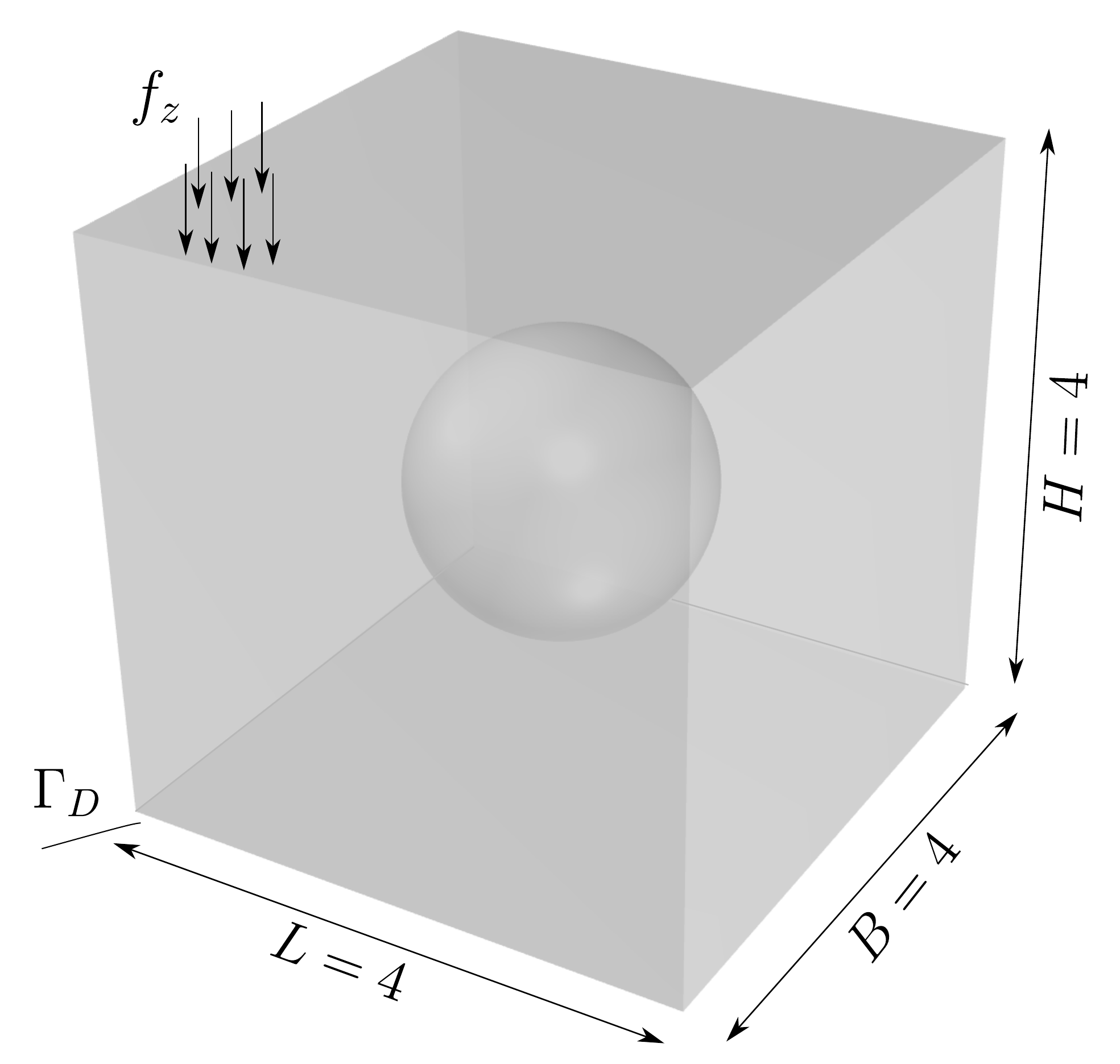}
         \caption{Geometry}
         \label{fig:cube}
     \end{subfigure}
     \begin{subfigure}[b]{0.49\textwidth}
         \centering
         \includegraphics[width=\textwidth]{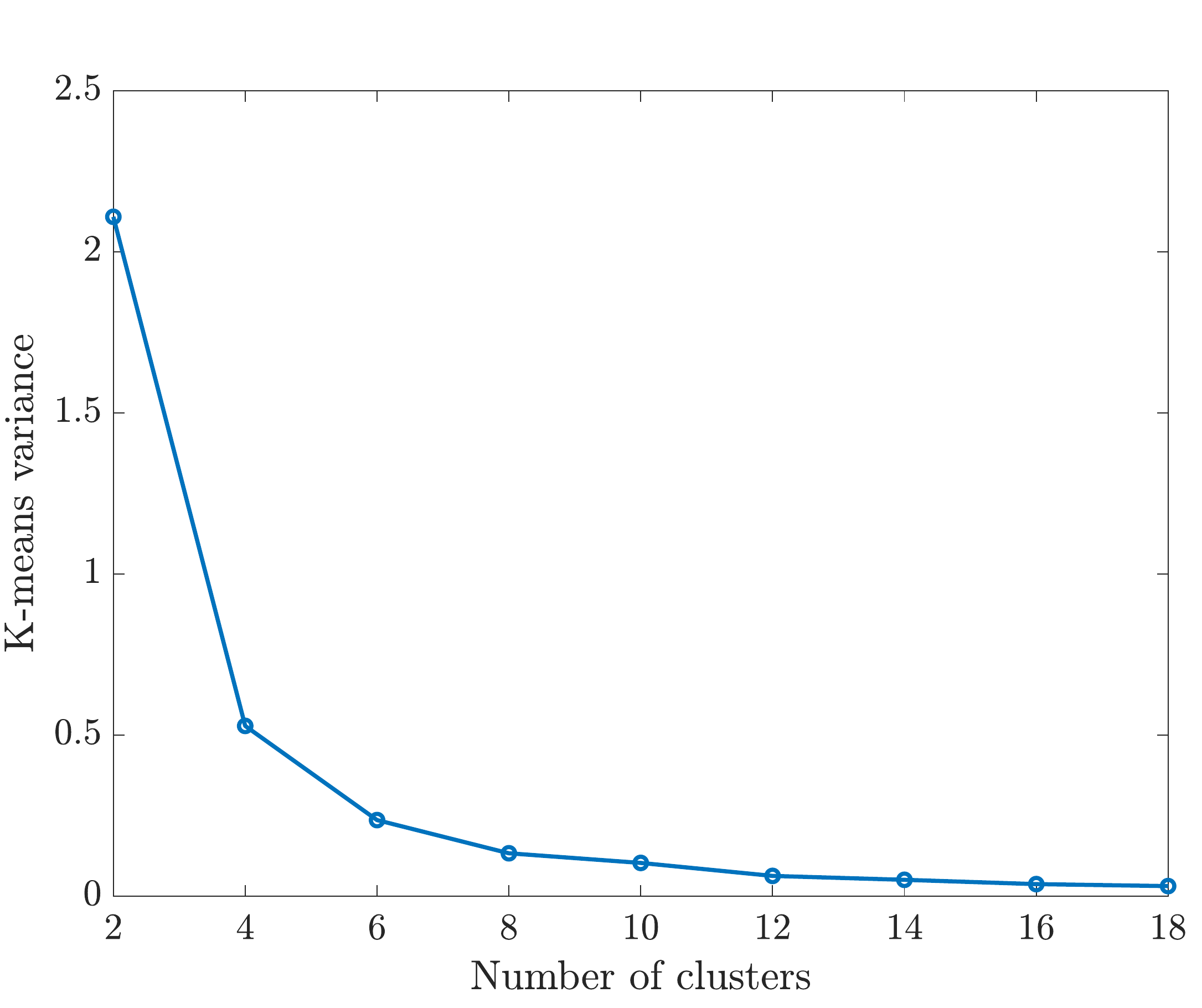}
         \caption{K-means variance}
         \label{fig:k_cube}
     \end{subfigure}
       	\caption{Example 6.2.2: Geometry of cube with spherical inclusion (a) and k-means variance versus number of clusters $N_c^d$ (b).}\label{fig:system_cube}
\end{figure}

This example aims to demonstrate the applicability of our approach to three dimensional geometries. For this purpose we consider a cube with a spherical inclusion. The model is defined in Figure \ref{fig:cube}. The geometric parameter we consider here is the radius of the sphere $\mu = R \in [0.5, 1.5]$, while the sphere is centered at the center of the cube.

We impose homogeneous Dirichlet boundary conditions on the bottom of the cube as depicted in Figure \ref{fig:cube}. The  body load is set to ${\bm{f}}=[{f}_x,{f}_y, {f}_z]=[0, 0, -10]$. The Young modulus and Poisson coefficient are given as $E=100$ and $\nu=0.3$, respectively. The geometry is discretized with quadratic $C^1$-continuous B-splines and the mesh consists of $8$ elements per direction over a Cartesian grid, resulting in $\mathcal{N}_{h,0}=3000$ degrees of freedom. Figure \ref{fig:k_cube} shows the \emph{k-means} variance \eqref{eq23} computed for different numbers of clusters $N_c^d$. Similarly to the previous test case, we observe that the variance does not decrease significantly after 10 clusters.

\begin{figure}[!h]
     \begin{subfigure}[b]{0.49\textwidth}
         \centering
         \includegraphics[width=\textwidth]{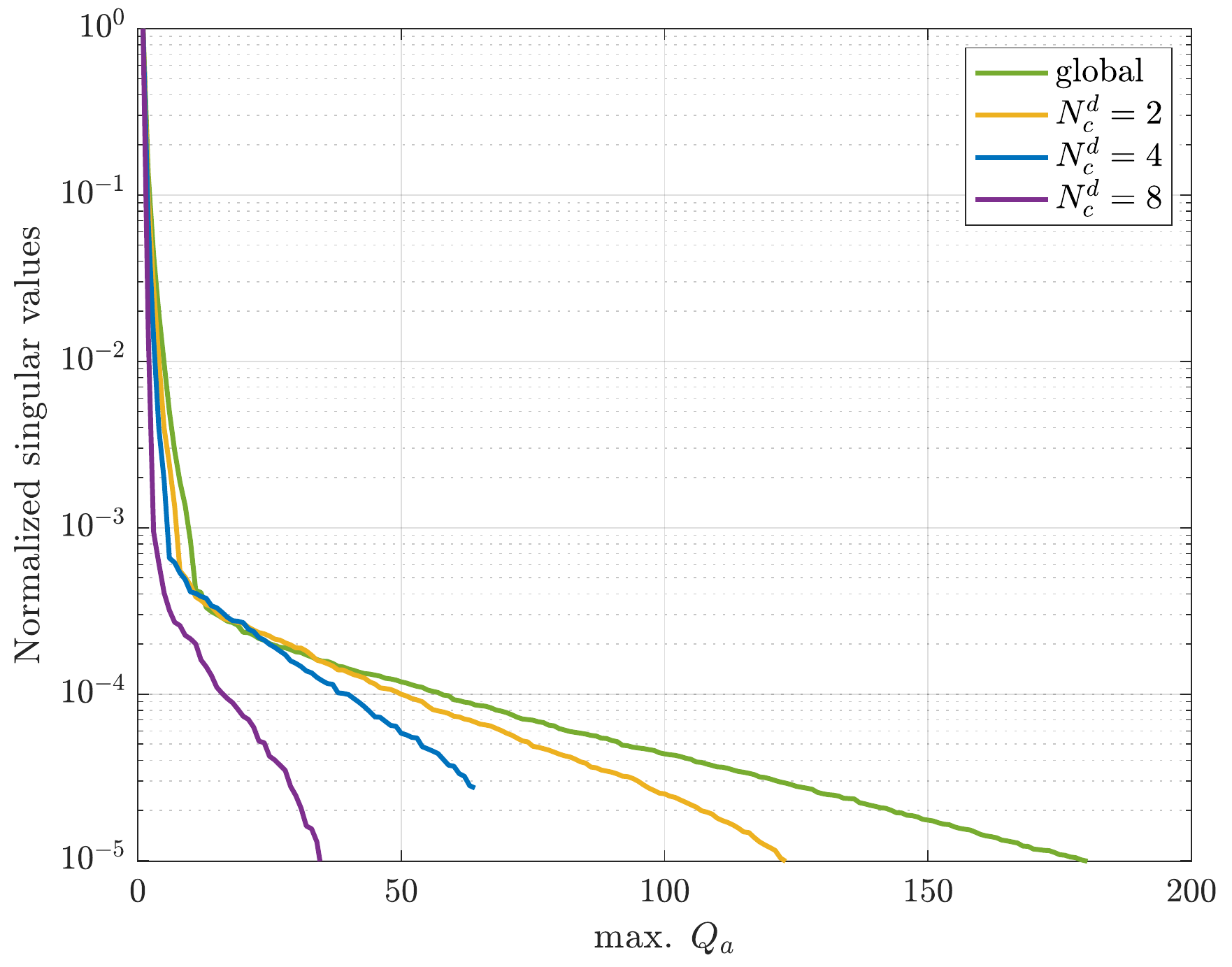}
         \caption{DEIM approximation}
         \label{fig:DEIM_3D}
     \end{subfigure}
     \begin{subfigure}[b]{0.49\textwidth}
         \centering
         \includegraphics[width=\textwidth]{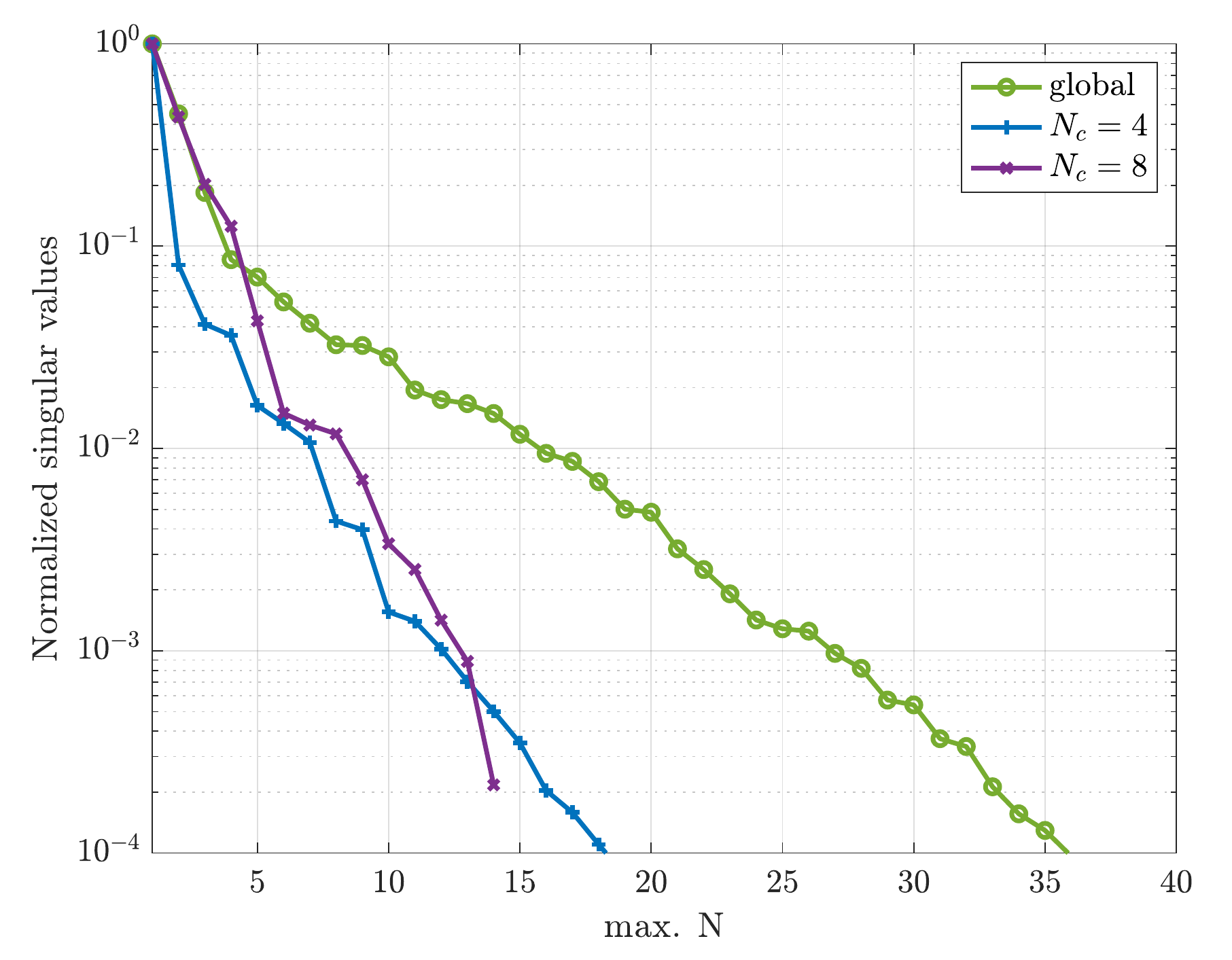}
         \caption{Reduced basis approximation}
         \label{fig:POD_3D}
     \end{subfigure}
        \caption{Example 6.2.2: Comparison of singular values decay for DEIM approximations of the stiffness matrix (a) and reduced basis approximation over all clusters (b) using different numbers of clusters.}
      \label{fig:results_3D}
\end{figure}

In the following, we will investigate the reducibility of the problem at hand. For this purpose, we consider a training set of dimension $N_{s}^d=250$ for the affine decomposition and $N_s=100$ for the reduced basis. Moreover, Figure \ref{fig:DEIM_3D} depicts the singular values decay of the DEIM approximation for the stiffness matrix versus the maximum number of basis functions $Q_a$ over all clusters. It is evident, that the maximum number of the DEIM basis functions over all clusters is reduced effectively to $Q_a=35$ with $N_{c}^d=8$ clusters, while the global approach requires $Q_a=180$. A similar behavior can be observed for the reduced basis in Figure \ref{fig:POD_3D}. The dimension of the reduced basis with $N_c=8$ clusters is reduced from $N=35$ to $N=14$ basis functions. It should be noted that the depicted decay corresponds to the cluster with the maximum number of basis functions for all cases. The decay is more rapid for the local ROMs, which implies that the solution can be captured with less basis functions and the problem at hand is more effectively reducible with the proposed localization strategy.

\section{Conclusions}
We have presented a novel reduced basis framework in the context of second-order linear elliptic PDEs defined on parameterized unfitted domains. Our approach is based on projection-based ROMs and techniques such as the reduced basis method and discrete empirical interpolation \cite{QMN_RBspringer,Negri2015}. The latter ensures an efficient offline/online procedure for problems formulated in parameterized geometries. To construct efficient ROMs for PDEs formulated on parameterized unfitted geometries, we proposed a methodology based on extension of snapshots within the cut regions and a localization strategy that reduces the dimension of the reduced basis. The presented framework allows an efficient offline/online decomposition with low online cost, while it is perfectly suitable for any discretization choice within an unfitted framework. 

 We have studied numerically the performance of the proposed methodology using the Poisson and linear elasticity problems. For this purpose, we considered trimmed spline discretizations by exploiting the re-parameterization tool for integration of cut elements in \cite{Antolin2019,Wei2021}. We observed a significant reduction of the computational cost in the online phase compared to standard ROMs, while we obtained accurate reduced basis approximations for problems distinguished by parameterized cut regions and strong variability of the solutions. Finally, we have applied the proposed strategy to a three-dimensional geometry in order to investigate the potential of our framework to achieve effective reduction. 
 
From the model reduction point of view, an interesting research direction for the future is the application of Greedy algorithms to construct localized reduced bases and error certification driven by a posteriori error estimators. Moreover, the extension to more involved problems, such as fourth-order PDEs, and complex geometrical representations is a further topic of interest. To the best of the authors' knowledge, this work comprises the first general methodology allowing reduced order modeling in the context of parameterized trimmed domains in isogeometric analysis. The proposed strategy paves the way for several applications involving complex shapes within a parametric framework, such as design, shape and topology optimization.

\section*{Acknowledgments}
The financial support of the Swiss Innovation Agency (Innosuisse) under Grant No. 46684.1 IP-EE is gratefully acknowledged. We would also like to thank Prof. Andrea Manzoni (Politecnico di Milano, Italy) and Dr. David Knezevic (Akselos S.A.) for fruitful discussions.

\bibliography{mybibfile}

\begin{thebibliography}{10}
\expandafter\ifx\csname url\endcsname\relax
  \def\url#1{\texttt{#1}}\fi
\expandafter\ifx\csname urlprefix\endcsname\relax\def\urlprefix{URL }\fi
\expandafter\ifx\csname href\endcsname\relax
  \def\href#1#2{#2} \def\path#1{#1}\fi

\bibitem{Peskin2002}
C.~S. Peskin, The immersed boundary method, Acta Numerica 11 (2002) 479–517.

\bibitem{Li2006}
Z.~Li, K.~Ito, The immersed interface method: numerical solutions of {PDE}s
  involving interfaces and irregular domains, {SIAM},Philadelphia, 2006.

\bibitem{parvizian2007finite}
J.~Parvizian, A.~D{\"u}ster, E.~Rank, Finite cell method, Computational
  Mechanics 41~(1) (2007) 121--133.

\bibitem{burman2015cutfem}
E.~Burman, S.~Claus, P.~Hansbo, M.~G. Larson, A.~Massing, Cut{FEM}:
  discretizing geometry and partial differential equations, Int. J. Numer.
  Meth. Engng. 104~(7) (2015) 472--501.

\bibitem{Main2018}
A.~Main, G.~Scovazzi, The shifted boundary method for embedded domain
  computations. {P}art {I}: {P}oisson and {S}tokes problems, J. Comput. Phys.
  372 (2018) 972--995.

\bibitem{Main2018b}
A.~Main, G.~Scovazzi, The shifted boundary method for embedded domain
  computations. {P}art {II}: Linear advection–diffusion and incompressible
  {N}avier–{S}tokes equations, J. Comput. Phys. 372 (2018) 996--1026.

\bibitem{Mittal2005}
R.~Mittal, G.~Iaccarino, Immersed boundary methods, Annu. Rev. Fluid Mech. 37
  (2005) 239 -- 261.

\bibitem{Hughes2005}
T.~J.~R. Hughes, J.~A. Cottrell, Y.~Bazilevs, Isogeometric analysis: {C}{A}{D},
  finite elements, {N}{U}{R}{B}{S}, exact geometry and mesh refinement, Comput.
  Methods Appl. Mech. Engrg. 194 (2005) 4135--4195.

\bibitem{Cottrell2007}
J.~A. Cottrell, T.~J.~R. Hughes, A.~Reali, Studies of refinement and continuity
  in isogeometric structural analysis, Comput. Methods Appl. Mech. Engrg. 196
  (2007) 4160--4183.

\bibitem{Akhras2016}
H.~A. Akhras, T.~Eljuedj, A.~Gravouil, M.~Rochette, Isogeometric
  analysis-suitable trivariate {NURBS} models from standard {B}-{R}ep models,
  Comput. Methods Appl. Mech. Eng. 307 (2016) 256--274.

\bibitem{Klinkel2020}
S.~Klinkel, M.~Chasapi, Isogeometric Analysis of Solids in Boundary
  Representation, Novel Finite Element Technologies for Solids and Structures,
  Springer, Cham, 2020, pp. 153--197.

\bibitem{Chasapi2021}
M.~Chasapi, L.~Mester, B.~Simeon, S.~Klinkel, Isogeometric analysis of {3D}
  solids in boundary representation for problems in nonlinear solid mechanics
  and structural dynamics, Int. J. Numer. Meth. Engng. 123 (2021) 1228--1252.

\bibitem{Hollig2001}
K.~H\"{o}llig, U.~Reif, J.~Wipper, Weighted extended {B}-spline approximation
  of {D}irichlet problems, SIAM Journal on Numerical Analysis 39~(2) (2001)
  442--462.

\bibitem{Schillinger2012}
D.~Schillinger, L.~Ded{\'e}, M.~Scott, J.~Evans, M.~Borden, E.~Rank, T.~J.~R.
  Hughes, An isogeometric design-through-analysis methodology based on adaptive
  hierarchical refinement of {NURBS}, immersed boundary methods, and {T}-spline
  {CAD} surfaces, Comput. Methods Appl. Mech. Engrg. 249 (2012) 116--150.

\bibitem{Elfeverson2018}
D.~Elfverson, M.~G. Larson, K.~Larsson, Cut{IGA} with basis function removal,
  Advanced Modeling and Simulation in Engineering Sciences 5 (2018) 1--19.

\bibitem{Rank2012}
E.~Rank, M.~Ruess, S.~Kollmannsberer, D.~Schillinger, A.~D{\"u}ster, Geometric
  modeling, isogeometric analysis and the finite cell method, Comput. Methods
  Appl. Mech. Eng. 249-252 (2012) 104–115.

\bibitem{Messmer2022}
M.~Me{\ss}mer, T.~Teschemacher, L.~F. Leidinger, R.~W{\"u}chner, K.-U.
  Bletzinger, Efficient {CAD}-integrated isogeometric analysis of trimmed
  solids, Comput. Methods Appl. Mech. Eng. 400 (2022) 115584.

\bibitem{Marussig2018}
B.~Marussig, T.~J.~R. Hughes, A review of trimming in isogeometric analysis:
  challenges, data exchange and simulation aspects, Arch. Computat. Methods
  Eng. 25 (2018) 1059--1127.

\bibitem{Nagy2015a}
A.~Nagy, D.~Benson, On the numerical integration of trimmed isogeometric
  elements, Comput. Methods Appl. Mech. Eng. 284 (2015) 165--185.

\bibitem{Kudela2015}
L.~Kudela, N.~Zander, T.~Bog, S.~Kollmannsberger, E.~Rank, Efficient and
  accurate numerical quadrature for immersed boundary methods, Adv. Model.
  Simul. Sci. 2 (2015) 1--22.

\bibitem{Antolin2019}
P.~Antolin, A.~Buffa, M.~Martinelli, Isogeometric analysis on {V}-reps: First
  results, Comput. Methods Appl. Mech. Engrg. 355 (2019) 976--1002.

\bibitem{Divi2020}
S.~Divi, C.~Verhoosel, F.~Aurrichio, A.~Reali, E.~van Brummelen,
  Error-estimate-based adaptive integration for immersed isogeometric analysis,
  Comput. Math. Appl. 80 (2020) 2481--2516.

\bibitem{Antolin2022}
P.~Antolin, X.~Wei, A.~Buffa, Robust numerical integration on curved polyhedra
  based on folded decompositions, Comput. Methods Appl. Mech. Eng. 395 (2022)
  114948.

\bibitem{Antolin2022b}
P.~Antolin, T.~Hirschler, Quadrature-free immersed isogeometric analysis,
  Engineering with Computers 38 (2022) 4475–4499.

\bibitem{Marussig2018b}
B.~Marussig, R.~Hiemstra, T.~J.~R. Hughes, Improved conditioning of
  isogeometric analysis matrices for trimmed geometries, Comput. Methods Appl.
  Mech. Eng. 334 (2018) 79--110.

\bibitem{Elfverson2019}
D.~Elfverson, M.~G. Larson, K.~Larsson, A new least squares stabilized
  {N}itsche method for cut isogeometric analysis, Comput. Methods Appl. Mech.
  Engrg. 349 (2019) 1--16.

\bibitem{Buffa2020}
A.~Buffa, R.~Puppi, R.~V\'azquez, A minimal stabilization procedure for
  isogeometric methods on trimmed geometries, SIAM J. Numer. Anal. 58 (2020)
  2711--2735.

\bibitem{Garotta2020}
F.~Garotta, N.~Demo, M.~Tezzele, M.~Carraturo, A.~Reali, G.~Rozza, Reduced
  order isogeometric analysis approach for {P}{D}{E}s, In Lecture Notes in
  Computational Science and Engineering - Quantification of Uncertainty:
  Improving Efficiency and Technology 137 (2020) 153--170.

\bibitem{Manzoni2015}
A.~Manzoni, F.~Salmoiraghi, L.~Heltai, Reduced basis isogeometric methods
  ({R}{B}-{I}{G}{A}) for the real-time simulation of potential flows about
  parametrized {N}{A}{S}{A} airfoils, Comput. Methods Appl. Mech. Engrg. 284
  (2015) 1147--1180.

\bibitem{Salmoiraghi2016}
F.~Salmoiraghi, F.~Ballarin, L.~Heltai, G.~Rozza, Isogeometric analysis-based
  reduced order modelling for incompressible linear viscous flows in
  parametrized shapes, Adv. Model. and Simul. in Eng. Sci. 3 (2016) 21.

\bibitem{Zhu2017}
Z.~Zhu, L.~Ded\'e, A.~Quarteroni, Isogeometric analysis and proper orthogonal
  decomposition for parabolic problems, Numer. Math. 135 (2017) 333--370.

\bibitem{Fresca2017}
S.~Fresca, A.~Manzoni, L.~Ded\'e, A.~Quarteroni, {P}{O}{D}-enhanced deep
  learning-based reduced order models for the real-time simulation of cardiac
  electrophysiology in the left atrium, Frontiers in Physiology 12 (2021)
  679076.

\bibitem{Maquart2020}
T.~Maquart, W.~Wenfeng, T.~Elguedj, A.~Gravouil, M.~Rochette, 3{D} volumetric
  isotopological meshing for finite element and isogeometric based reduced
  order modeling, Comput. Methods Appl. Mech. Engrg. 362 (2020) 112809.

\bibitem{Devaud2017}
D.~Devaud, G.~Rozza, Certified reduced basis method for affinely parametric
  isogeometric analysis {N}{U}{R}{B}{S} approximation, In Lecture Notes in
  Computational Science and Engineering: Spectral and Higher Order Methods for
  Partial Differential Equations 119 (2017) 41--62.

\bibitem{Hesthaven2016}
J.~S. Hesthaven, G.~Rozza, Certified reduced basis methods for parametrized
  partial differential equations, Springer, 2016.

\bibitem{QMN_RBspringer}
A.~Quarteroni, A.~Manzoni, F.~Negri, Reduced Basis Methods for Partial
  Differential Equations. An Introduction, Vol.~92 of Unitext, Springer, 2016.

\bibitem{Chasapi2022}
M.~Chasapi, P.~Antolin, A.~Buffa, Reduced order modelling of nonaffine problems
  on parameterized {N}{U}{R}{B}{S} multipatch geometries, arXiv (2022)
  arXiv:2211.07348.

\bibitem{Gabriel2022}
H.~D. I.~Gabriel, D.~Loukrezis, Tensor train based isogeometric analysis for
  pde approximation on parameter dependent geometries, Comput. Methods Appl.
  Mech. Eng. 401 (2022) 115593.

\bibitem{Nouy2011}
A.~Nouy, M.~Chevreuil, E.~Safatly, Fictitious domain method and separated
  representations for the solution of boundary value problems on uncertain
  parameterized domains, Comput. Methods Appl. Mech. Eng. 200~(45) (2011)
  3066--3082.

\bibitem{Balajewicz2014}
M.~Balajewicz, C.~Farhat, Reduction of nonlinear embedded boundary models for
  problems with evolving interfaces, J. Comput. Phys. 274 (2014) 489--504.

\bibitem{Karatzas2020}
E.~N. Karatzas, F.~Ballarin, G.~Rozza, Projection-based reduced order models
  for a cut finite element method in parametrized domains, Comput. Math. Appl.
  79 (2020) 833--851.

\bibitem{Karatzas2021}
E.~N. Karatzas, G.~Rozza, A reduced order model for a stable embedded boundary
  parametrized {C}ahn–{H}illiard phase-field system based on cut finite
  elements, J. Sci. Comput. 89:9.

\bibitem{Karatzas2019}
E.~N. Karatzas, G.~Stabile, L.~Nouveau, G.~Scovazzi, G.~Rozza, A reduced basis
  approach for {PDE}s on parametrized geometries based on the shifted boundary
  finite element method and application to a {S}tokes flow, Comput. Methods
  Appl. Mech. Engrg. 347 (2019) 568 -- 587.

\bibitem{Karatzas2020b}
E.~N. Karatzas, G.~Stabile, L.~Nouveau, G.~Scovazzi, G.~Rozza, A reduced-order
  shifted boundary method for parametrized incompressible {N}avier–{S}tokes
  equations, Comput. Methods Appl. Mech. Engrg. 370 (2020) 113273.

\bibitem{Zeng2022}
X.~Zeng, G.~Stabile, E.~N. Karatzas, G.~Scovazzi, G.~Rozza, Embedded domain
  {R}educed {B}asis {M}odels for the shallow water hyperbolic equations with
  the {S}hifted {B}oundary {M}ethod, Comput. Methods Appl. Mech. Eng. 398
  (2022) 115143.

\bibitem{Karatzas2022}
G.~Katsouleas, E.~N. Karatzas, F.~Travlopanos, Discrete empirical interpolation
  and unfitted mesh {FEM}s: application in {PDE}-constrained optimization,
  Optimization (2022) 1--34.

\bibitem{Barrault2004}
M.~Barrault, Y.~Maday, N.~C. Nguyen, A.~T. Patera, An ‘empirical
  interpolation’ method: application to efficient reduced-basis
  discretization of partial differential equations, C. R. Acad. Sci. Paris,
  Ser. I 339 (2004) 667--672.

\bibitem{Chaturantabut2010}
S.~Chaturantabut, D.~C. Sorensen, Nonlinear model reduction via discrete
  empirical interpolation, SIAM J. Sci. Comput. 32 (2010) 2737--2764.

\bibitem{Negri2015}
F.~Negri, A.~Manzoni, D.~Amsallem, Efficient model reduction of parametrized
  systems by matrix discrete empirical interpolation, J. Comput. Phys. 303
  (2015) 431--454.

\bibitem{Benner2021}
N.~Sarnaa, P.~Benner, Data-driven model order reduction for problems with
  parameter-dependent jump-discontinuities, Comput. Methods Appl. Mech. Engrg.
  387 (2021) 114168.

\bibitem{Amsallem2012}
D.~Amsallem, M.~J. Zahr, C.~Farhat, Nonlinear model order reduction based on
  local reduced-order bases, Int. J. Numer. Meth. Engng. 92 (2012) 891--916.

\bibitem{Peherstorfer2014}
B.~Peherstorfer, D.~Butnaru, K.~Willcox, H.-J. Bungartz, Localized discrete
  empirical interpolation method, SIAM J. Sci. Comput. 36~(1) (2014)
  A168--A192.

\bibitem{Pagani2018}
S.~Pagani, A.~Manzoni, A.~Quarteroni, Numerical approximation of parametrized
  problems in cardiac electrophysiology by a local reduced basis method,
  Comput. Methods Appl. Mech. Engrg. 340 (2018) 530--558.

\bibitem{Hess2019}
M.~Hess, A.~Alla, A.~Quaini, G.~Rozza, M.~Gunzburger, A localized reduced-order
  modeling approach for {PDE}s with bifurcating solutions, Comput. Methods
  Appl. Mech. Engrg. 351 (2019) 379--408.

\bibitem{Eftang2010a}
J.~L. Eftang, A.~T. Patera, E.~M. R{\o}nquist, An hp certified reduced basis
  method for parametrized elliptic partial differential equations, SIAM J. Sci.
  Comput. 32 (2010) 3170--3200.

\bibitem{Eftang2011}
J.~L. Eftang, D.~J. Knezevic, A.~T. Patera, A hp certified reduced basis method
  for parametrized parabolic partial differential equations, Mathematical and
  Computer Modelling of Dynamical Systems 17:4 (2011) 395--422.

\bibitem{Eftang2012}
J.~L. Eftang, B.~Stamm, Parameter multi-domain hp empirical interpolation, Int.
  J. Numer. Meth. Engng. 90 (2012) 412--428.

\bibitem{Haasdonk2011}
B.~Haasdonk, M.~Dihlmann, M.~Ohlberger, A training set and multiple bases
  generation approach for parameterized model reduction based on adaptive grids
  in parameter space, Mathematical and Computer Modelling of Dynamical Systems
  17:4 (2011) 423--442.

\bibitem{Maday2013}
Y.~Maday, B.~Stamm, Locally adaptive {G}reedy approximations for anisotropic
  parameter reduced basis spaces, SIAM J. Sci. Comput. 35 (2013) A2417--A2441.

\bibitem{Wei2021}
X.~Wei, B.~Marussig, P.~Antolin, A.~Buffa, Immersed boundary-conformal
  isogeometric method for linear elliptic problems, Computational Mechanics 68
  (2021) 1385 -- 1405.

\bibitem{BADIA2018533}
S.~Badia, F.~Verdugo, A.~F. Mart{\'\i}n, {The aggregated unfitted finite
  element method for elliptic problems}, Comput. Methods Appl. Mech. Engrg. 336
  (2018) 533 -- 553.

\bibitem{Maday2009}
Y.~Maday, N.~C. Nguyen, A.~T. Patera, S.~H. Pau, A general multipurpose
  interpolation procedure: the magic points, Communications on Pure and Applied
  Analysis 8~(1) (2009) 383--404.

\bibitem{Powell1992}
M.~J.~D. Powell, The theory of radial basis functions approximation in 1990,
  Advances in Numerical Analysis II: Wavelets, Subdivision Algorithms and
  Radial Functions, Oxford University Press, Oxford, 1992.

\bibitem{Wirtz2014}
D.~Wirtz, D.~C. Sorensen, B.~Haasdonk, A posteriori error estimation for {DEIM}
  reduced nonlinear dynamical systems, SIAM Journal on Scientific Computing
  36~(2) (2014) A311--A338.

\bibitem{Stewart1990}
G.~Stewart, J.~Sun, Matrix Perturbation Theory, Academic Press, New York, 1990.

\bibitem{Georgaka2020}
S.~Georgaka, G.~Stabile, K.~Star, G.~Rozza, M.~J. Bluck, A hybrid reduced order
  method for modelling turbulent heat transfer problems, Computer and Fluids
  208 (2020) 104615.

\bibitem{Buhmann2000}
M.~D. Buhmann, Radial basis functions, Acta Numerica (2000) 1 -- 38.

\bibitem{Beatson1999}
R.~Beatson, J.~Cherrie, C.~Mouat, Fast fitting of radial basis functions:
  Methods based on preconditioned {GMRES} iteration, Advances in Computational
  Mathematics 11 (1999) 253--270.

\bibitem{Fornberg2007}
B.~Fornberg, J.~Zuev, The {R}unge phenomenon and spatially variable shape
  parameters in {R}{B}{F} interpolation, Computers {\&} Mathematics with
  Applications 54~(3) (2007) 379--398.

\bibitem{Likas2003}
A.~Likas, N.~Vlassis, J.~Verbeek, The global k-means clustering algorithm,
  Pattern Recognition 36.2 (2003) 452 -- 461.

\bibitem{Coradello2020}
L.~Coradello, P.~Antolin, R.~V\'azquez, A.~Buffa, Adaptive isogeometric
  analysis on two-dimensional trimmed domains based on a hierarchical approach,
  Comput. Methods Appl. Mech. Engrg. 264 (2020) 112925.

\bibitem{Coradello2021}
L.~Coradello, J.~Kiendl, A.~Buffa, Coupling of non-conforming trimmed
  isogeometric {K}irchhoff-{L}ove shells via a projected super-penalty
  approach, Comput. Methods Appl. Mech. Engrg. 387 (2021) 114187.

\bibitem{Buffa2022}
A.~Buffa, O.~Chanon, R.~V\'azquez, An a posteriori error estimator for
  isogeometric analysis on trimmed geometries, IMA Journal of Numerical
  Analysis 00 (2022) 1--29.

\bibitem{Vazquez2016}
R.~V\'azquez, A new design for the implementation of isogeometric analysis in
  {O}ctave and {M}atlab: {G}eo{P}{D}{E}s 3.0, Comput. Math. Appl. 72 (2016)
  523--554.

\bibitem{redbKIT}
F.~Negri, redb{KIT} {V}ersion 2.2, \url{http://redbkit.github.io/redbKIT/}
  (2016).

\bibitem{Deprenter2017}
F.~{de Prenter}, C.~Verhoosel, G.~{van Zwieten}, E.~{van Brummelen}, Condition
  number analysis and preconditioning of the finite cell method, Comput.
  Methods Appl. Mech. Eng. 316 (2017) 297--327.

\bibitem{McKay1979}
M.~McKay, R.~Beckman, W.~Conover, Comparison of three methods for selecting
  values of input variables in the analysis of output from a computer code,
  Technometrics 21~(2) (1979) 239 – 245.

\end{thebibliography}

\end{document}